\documentclass[11pt, oneside]{article}
\usepackage{amsfonts}
\usepackage{amsfonts}
\usepackage{mathrsfs}
\usepackage{amsfonts}
\usepackage{float}
\usepackage{graphicx}
\usepackage{epstopdf}
\usepackage{caption} 
\usepackage{array}
\usepackage{amsmath}
\usepackage{bm}
\usepackage{booktabs}
\usepackage{multirow}

\usepackage{latexsym,amssymb,amsmath,amstext}
\usepackage{mathtools}
\usepackage{amsthm}
\usepackage{indentfirst}
\usepackage{color}
\usepackage{graphicx}
\usepackage{multirow,booktabs}
\usepackage{tabularx}
\usepackage{subfigure}
\usepackage{graphicx}
\usepackage[ruled,vlined,linesnumbered]{algorithm2e}
\usepackage{enumitem} 
\usepackage[colorlinks,linkcolor=blue,anchorcolor=blue,citecolor=blue]{hyperref}
\usepackage[authoryear,round,comma,sort&compress]{natbib}
\usepackage{xcolor}

\numberwithin{equation}{section}
\allowdisplaybreaks

\newtheorem{theorem}{\color{black}\indent Theorem}[section]
\newtheorem{lemma}{\color{black}\indent Lemma}[section]

\newtheorem{proposition}{\color{black}\indent Proposition}[section]
\newtheorem{definition}{\color{black}\indent Definition}[section]
\makeatletter
\def\th@remark{%
  \normalfont 
  \let\@begintheorem\@afterindenttrue\@afterheading\flushleft 
  \let\@endtheorem\endflushleft 
}
\makeatother

\theoremstyle{plain} 
\newtheorem{remark}{\color{black}Remark}

\newtheorem{corollary}{\color{black}\indent Corollary}[section]

\usepackage{amssymb,amsmath}
\begin{document} 

\title{Reduce-Rank Matrix Integer-Valued Autoregressive Model}
\author{\\Kaiyan Cui$^{1}$, Tianyun Guo$^{1}$, Suping Wang$^{2}$\\
{$^{1}$\small School of Mathematics and Statistics, Shanxi University, Taiyuan 030006, China}\\
{$^{2}$\small School of Applied Sciences, Taiyuan University of Science and Technology, Taiyuan, 030024, China.}\\
{\small (E-mail: guilinwangsuping@126.com )}
}

\date{}
 \maketitle
\noindent {\bf Abstract}:
Integer-valued time series are widely present in many fields, such as finance, economics, disease transmission, and traffic flow. With data dimensions surging, the traditional multivariate generalized integer autoregressive (MGINAR) model faces parameter overload, poor interpretability, and structural information loss. Matrix integer-valued autoregression (MINAR) model captures row-column cross-correlations and reduces the number of parameters to be estimated. However, further growth in dimensionality causes data redundancy, which degrades the MINAR model's performance and increases the number of parameters. To solve the limitations of the MINAR model described above, this paper proposes the reduced-rank matrix integer-valued autoregression (RRMINAR) model. Reducing rank is achieved by adding low-rank constraints to the coefficient matrices in the MINAR model, leading to RRMINAR reducing parameter quantity while incorporating matrix structure information. We develop an iterative conditional least squares estimation and analyze its asymptotic properties. Simulation results demonstrate that the proposed RRMINAR model exhibits more robust parameter estimation and higher prediction accuracy than MGINAR and MINAR models when the data structure is low-rank. Empirical analysis using criminal data validates the proposed RRMINAR model's effectiveness and uncovers structural temporal-spatial information in criminal behavior.

\thispagestyle{empty}
\section{Introduction}
Integer-time series data are closely related to people's daily production and life, and frequently appear in many fields such as epidemiology, biology, engineering, and finance. For example, the number of fires occurring in an area within a  period of time, the monthly sales volume of a product in a shopping mall, the number of COVID-19 patients in a certain area, etc., are all specific manifestations of integer-valued time series data. As described in \cite{Wei2021}, there are mainly two types of methods for modeling integer-valued time series data currently: (1) Conditional regression model. (2) Model based on thinning operator. An introduction to the former can be found in \cite{Ferland2006} and \cite{EncisoMora2009}, while the latter represents the main current research direction.

\cite{Steutel1979} first proposed the thinning operator, denoted by $\circ$, to construct integer-valued autoregressive (INAR) models. This operator solves the problem that multiplying an integer random variable by a real number typically results in a non-integer value. Specifically, the thinning operator $\circ$ is defined as follows: 
\[\alpha\circ x=\sum_{i=1}^x \xi_i,\]
where \(x\) is a non-negative integer, and \(\{\mathbf{\xi}_i\}\) is an i.i.d. (independent and identically distributed) Bernoulli random sequence with parameter \(\alpha\) and independent of \(x\). 

Building upon the thinning operator $\circ$, \cite{Al-Osh1987} introduced the first-order integer autoregressive (INAR(1)) model, defined as follows:
\[x_t =\alpha \circ x_{t-1} + {e}_t,~t\in\mathbb{Z},\]
where \(\{e_{t}\}\) is an i.i.d. Bernoulli random sequence with the same parameter \(\lambda\). \cite{Alzaid1988} and \cite{Al-Osh1992} developed INAR(1) models whose marginal distributions follow the geometric, generalized Poisson, and negative binomial distributions, respectively. \cite{Li1991} extended the INAR(1) model to a higher-order INAR(p) process and investigated its robustness and parameter estimation methods. \cite{Bu2008} studied the estimation problem for the INAR(p) model and proposed the maximum likelihood estimator. The study compared the performance of this estimator with the conditional least squares (CLS) and conditional maximum likelihood (CML) estimators for the Poisson INAR(2) model. To better analyze real-world data, the INAR model has been expanded in more directions. For example, \cite{Jazi2012} proposed the first-order zero-inflation integer autoregressive (ZI-INAR(1)) model to process zero-inflation data, and the experimental effect on the Pittsburgh drug crime data was good. \cite{Pang2023} proposed a time-varying parameter integer-valued autoregressive (TV-INAR) model that allows the parameter $\alpha$ to vary over time. The model demonstrated a better fit than conventional INAR models when applied to sexual offense data from Barina City, Australia.

The INAR model is unable to adapt to the diversified development of data; therefore \cite{Franke1993} introduced a significant multivariate INAR(1) model constructed using a matrix of univariate independent binomial thinning operators. Subsequently, \cite{Latour1997} proposed another important extension: the multivariate generalized integer autoregressive (MGINAR) model. This model replaces the binomial thinning operators with generalized thinning operators based on counting random variables that are not necessarily Bernoulli-distributed. \cite{Kirchner2017} further extended the multivariate INAR(1) to a $p$-th order model, estimating parameters via CLS and establishing the consistency and asymptotic normality of the estimators. \cite{Pedeli2011} and \cite{Pedeli2013a} investigated the properties and estimation of multivariate INAR models with bivariate and multivariate Poisson-distributed error terms, respectively. To address over-dispersion in multivariate integer-valued data, \cite{Yu2024} proposed a multivariate INAR(1) model with a mixed innovation term, featuring multivariate log-normal and multivariate log-Poisson distributions. To effectively model multivariate integer-valued time series with periodic characteristics, \cite{Santos2021} proposed a periodic multiple integer autoregressive (PMINAR(1)) model, incorporating periodic time-varying parameters and multiple negative binomial innovations. These studies have further promoted the development of integer-valued time series analysis.

Due to the increasing complexity and dynamic nature of data, the scope of analysis has broadened to incorporate multiple variables observed across various entities. These data can be structured in matrix form, as exemplified by economic indicators such as the gross domestic product (GDP) of multiple countries, bank interest rates, and total industrial production. Over time, the continuous observation of such matrix-form data gives rise to a matrix time series. This concept was first elaborated in the signal processing literature by \cite{Walden2002}. In recent years, the academic community has demonstrated growing interest in the modeling and analysis of this type of time series and has achieved a series of innovative results. 

The modeling and analysis of matrix time series can be broadly categorized into two approaches. The first category comprises methods based on vectorizing the matrix data, such as \cite{Basu2015} and \cite{Bauer2011}. These methods reshape the data to fit a vector autoregressive (VAR) model, but this leads to both a loss of structural information and over-parameterization. The second category of methods preserves the original data structure. These approaches utilize matrix factorization techniques to construct models directly from the matrix observations (\cite{Wang2019, Chen2019a}). The first-order matrix autoregressive (MAR(1)) model proposed by \cite{Chen2021} deeply describes the dynamic interaction between rows and columns in the matrix time series. Compared with VAR model ((\cite{Hannan1970}, \cite{Lutkepohl2005}), MAR method significantly reduces the number of the model parameters, and enhances the structural information. \cite{Han2024} further proposed a model based on the bilinear transformation, which simplifies the model structure. All the above studies are conducted within the real number domain. Even so, the above-mentioned progresses have opened up a new path for the theoretical and application research of matrix integer-valued time series. Recently, \cite{Xu2024} introduced a matrix integer-valued autoregressive (MINAR) model and employed both the projection estimation method and the iterative conditional least squares estimation (ICLSE) method for parameter estimation.

As data dimensions increase, redundant information impairs modeling performance and estimation efficiency, making it necessary to employ alternative methods to mitigate the curse of dimensionality. Although the matrix sequence can be stacked column by column and transformed into a vector sequence for dimensionality reduction or thinning processing, e.g. the regularization and thinning regression techniques(\cite{Basu2015} and \cite{Bauer2011}), and dimensionality reduction techniques, such as the reduce-rank vector autoregressive (RRVAR) model (\cite{Cubadda2021}), blind stacking will destroy the matrix data structure and lose the mutual influence of specific dimensions. To reduce the number of parameters without destroying the data structure, \cite{Xiao2022}, \cite{Liu2023}, and others proposed a reduced-rank matrix autoregressive (RRMAR) model based on the MAR model. They applied rank reduction to its bilinear coefficient matrix to mitigate the effects of high dimensionality and noise. Regrettably, existing methods in integer-valued time series remain underdeveloped, restricting their applicability in high-dimensional contexts.

To extend the theoretical results and application scope of matrix integer-valued time series, this study aims to conduct and analyze a reduce-rank autoregressive model for matrix integer-valued time series data. In particular, we propose the reduced-rank matrix integer-valued autoregressive (RRMINAR) model. Rank reduction is achieved by imposing a low-rank constraint on the dual coefficients of the MINAR model. It is worth noting that the low-rank structure of the proposed RRMINAR model does not affect the interpretability of its parameters, and it requires fewer parameters than the MINAR model proposed by \cite{Xu2024}. The estimation process for this model is based on reduced-rank regression, and we employ an iterative conditional least squares estimator. The asymptotic properties of the model parameters are derived. Experimental results demonstrate that, compared with the MGINAR and MINAR models, the RRMINAR model exhibits faster parameter convergence efficiency, achieves greater dimensionality reduction and noise removal, and reduces the  parameter estimation errors. For real-world data, the RRMINAR model demonstrates superior fitting performance in both in-sample and out-of-sample predictions compared with the other six benchmark models.

The rest of this article is organized as follows. Section 2 introduces the definition and probabilistic properties of the MINAR(1) model. Section 3 provides the definition and parameter estimation method of the RRMINAR(1) model. Section 4 presents the statistical inferences and asymptotic properties of the RRMINAR(1) model. Section 5 reports the relevant experimental results. Section 6 contains a brief summary. All the proofs are collected in the Appendix.

\noindent{\it Notation.} 
We use $\otimes$ to denote the Kronecker product, $\|\cdot\|_F$ the Frobenius norm of a matrix, and $\|\cdot\|_2$ the $L^{2}$-norm of a matrix. The bold capital letters, e.g., $\mathbf{A}$, denote matrices, and the bold lowercase letters, e.g., $\mathbf{v}$, denote vectors. $\mathbf{A}(i,j)$ denotes element in the \(i\)th row and \(j\)th column of matrix $\mathbf{A}$. Let $\mathrm{vec}(\cdot)$ be the vectorization of a matrix by stacking its columns. Define $\rho(\cdot)$ as the spectral radius of a matrix, i.e., the maximum modulus of its eigenvalues. $\operatorname{diag}(\mathbf{v})$ denotes the diagonal matrix formed with the elements of vector $\mathbf{v}$ as its main diagonal entries; $\operatorname{tr}(\mathbf{A})$ and $\operatorname{rank}(\mathbf{A})$ denote the trace and rank of the matrix $\mathbf{A}$, respectively. $\mathbf{I}_m$ and $\mathbf{0}_{m \times n}$ denote the $m$-dimensional identity matrix and the $(m, n)$-dimensional all-zero matrix, respectively.

\section{Matrix integer-valued autoregressive model}
\subsection{Definition and model}
In this section, we introduce some integer thinning operators, i.e., the thinning operator \(\circ\), vector thinning operator \(\ast\), and matrix thinning operator \(\circledast\), in Definitions \ref{def:2.1}-\ref{def:2.3}, respectively.

{\definition\label{def:2.1} For an $\mathbb{N}_0$-valued random variable $x$ and a constant $\alpha\geq0$, define the thinning operator $\circ$ by \[\alpha\circ x\coloneqq \sum_{k=1}^{x}\xi_{k}^{(\alpha)},\]
where $\xi_{1}^{(\alpha)}, \xi_{2}^{(\alpha)}, \ldots$ are {\rm i.i.d.} and independent of $Y$ with $\xi_{1}^{(\alpha)}\sim\text{Poisson}(\alpha)$. We use the convention that $\sum_{k=1}^{0}\xi_{k}^{(\alpha)}=0$.
}

{\definition \label{def:2.2}For a $m\times m$ matrix $\mathbf{A}=(\alpha_{i,j})_{m\times m}\in \mathbb{R}_{\geq0}^{m\times m}$ and an $\mathbb{N}_0^{m}$-valued random vector $\mathbf{v}=(v_1,\ldots,v_d)^\top$, define the vector (multivariate) thinning operator $\ast$ by
\[\mathbf{A}\ast\mathbf{v}\coloneqq\left(\sum_{k=1}^{m}\alpha_{i,k}\circ v_{k}\right)_{1\leq i\leq m},\]
where the thinning operator $\alpha_{i,k}\circ\cdot$ independently over $1\leq i,k\leq m$.
}

{\definition \label{def:2.3}For a $m\times m$ matrix $\mathbf{A}=(\alpha_{i,j})_{m\times m}\in \mathbb{R}_{\geq0}^{m\times m}$, an $n\times n$ matrix $\mathbf{B}=(\beta_{i,j})_{n\times n}\in \mathbb{R}_{\geq0}^{n\times n}$ and an $\mathbb{N}_0^{m\times n}$-valued random matrix $\mathbf{Y}=(y_{i,j})_{m\times n}$, define the single matrix thinning operator $\circledast$ by
\[\mathbf{A}\circledast\mathbf{Y}\coloneqq \left(\sum_{k=1}^{m}\alpha_{i,k}\circ y_{k,j}\right)_{1\leq i\leq m,1\leq j\leq n},\]
\[\mathbf{Y}\circledast \mathbf{B}\coloneqq \left(\sum_{l=1}^{n}\beta_{l,j}\circ y_{i,l}\right)_{1\leq i\leq m,1\leq j\leq n}.\]
}

{\definition \label{def:2.4}For a $m\times m$ matrix $\mathbf{A}=(\alpha_{i,j})_{m\times m}\in \mathbb{R}_{\geq0}^{m\times m}$, an $n\times n$ matrix $\mathbf{B}=(\beta_{i,j})_{n\times n}\in \mathbb{R}_{\geq0}^{n\times n}$ and an $\mathbb{N}_0^{m\times n}$-valued random matrix $\mathbf{Y}=(y_{i,j})_{m\times n}$, define the multiple matrix thinning operator $\circledast$ by
\[\mathbf{A}\circledast\mathbf{Y}\circledast \mathbf{B}^{\top}\coloneqq (\mathbf{A}\circledast\mathbf{Y})\circledast \mathbf{B}^{\top}=\mathbf{A}\circledast(\mathbf{Y}\circledast \mathbf{B}^{\top})=\left(\sum_{l=1}^{n}\sum_{k=1}^{m}(\beta_{j,l}\alpha_{i,k})\circ y_{k,l}\right)_{1\leq i\leq m,1\leq j\leq n},\]
where the thinning operator $(\beta_{j,l}\alpha_{i,k})\circ\cdot$ independently over $1\leq i,k\leq m$, $1\leq j,l\leq n$.
}

\begin{figure}[htbp]
  \centering
  \includegraphics[width=1\textwidth]{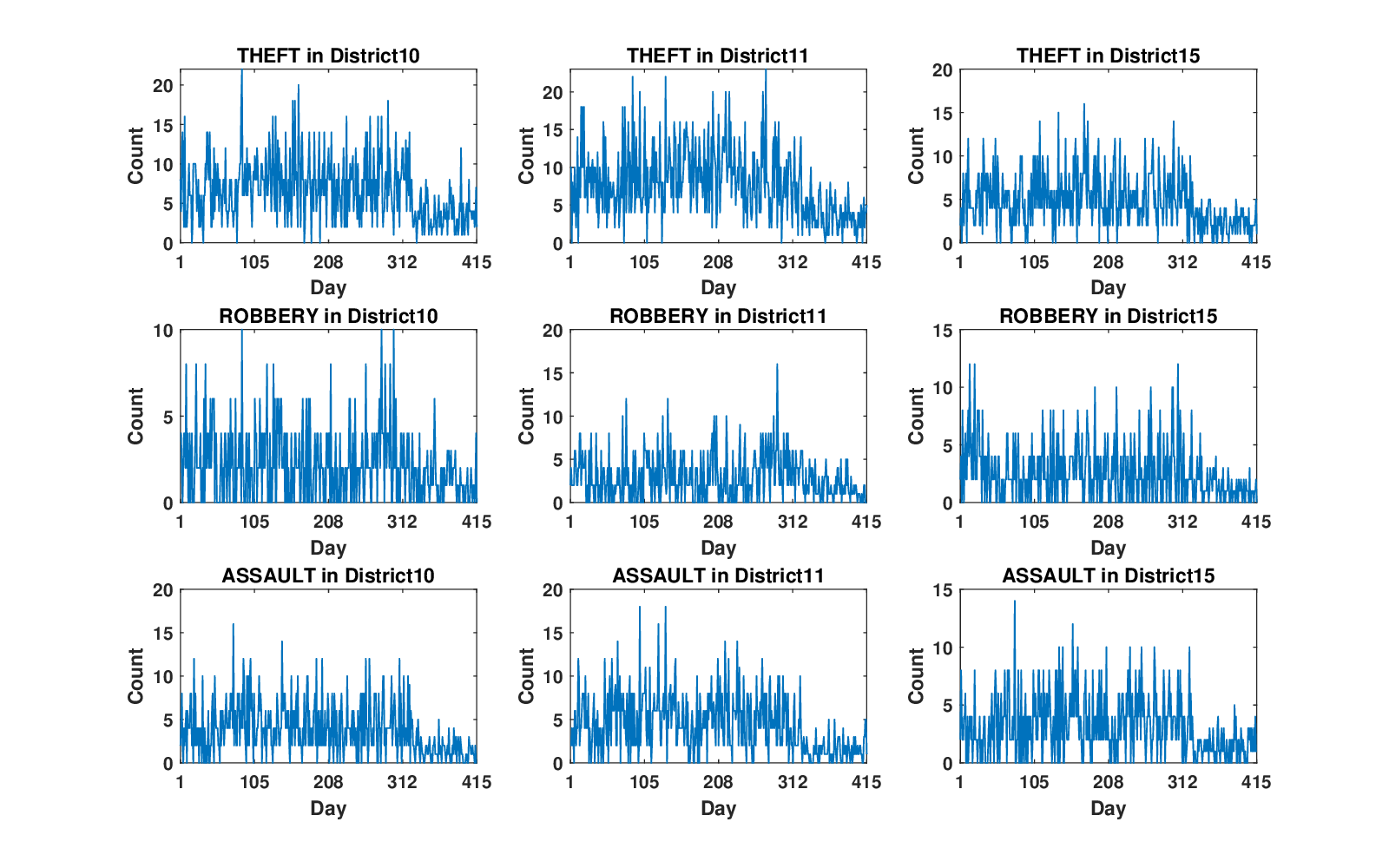}
  \caption{Time series of three types of crime (THEFT, ROBBERY and ASSAULT) in the three districts of Chicago.} 
  \label{fig:DQG}
  \end{figure}

Let $\mathbf{X}_t$ be a non-negative integer-valued random matrix of dimension $m \times n$, $\Phi$ be a non-negative matrix of dimension $mn \times mn$. Based on Definition \ref{def:2.2}, the first-order multivariate generalized integer autoregressive (MGINAR(1)) model proposed by \cite{Latour1997} is constructed in the form
\begin{equation}
\mathrm{vec}(\mathbf{X}_t)=\Phi\ast\mathrm{vec}(\mathbf{X}_{t-1})+\mathbf{e}_t,
\end{equation}
where \(\{\mathbf{e}_t\}\) is an i.i.d. random vector sequence. The matrix data is stacked and mixed together column-wise without taking into account the structural information of the matrix, ignoring the interactions between rows and columns in the data. A detailed explanation is provided below using an example of crime data. In Figure \ref{fig:DQG}, we plot data on three crime types from three districts of Chicago and obtain a $3 \times 3$ matrix observed at each time point. The rows and columns of the observed matrix correspond to different crime types and districts, respectively. In the MGINAR(1) model, \(\mathrm{vec}(\mathbf{X}_t)\) denotes the vectorization of matrix time series data representing different crime types across different districts. However, the coefficient matrix \(\Phi\) in the MGINAR(1) model simply mixes districts and crime types together to explore their interactions, without considering the strong interactions between different districts and different crime types. This implies that the coefficient matrix \(\Phi\) has limitations in interpreting the relationship between crime types and districts. Moreover, when handling high-dimensional data, it results in the number of model parameters multiplying. To solve the above issues, the first-order matrix integer-valued autoregressive (MINAR(1)) model is proposed as shown in Definition \ref{def:2.5}.

{\definition \label{def:2.5} Let $\mathbf{A}\in \mathbb{R}_{\geq0}^{m\times m}$, $\mathbf{B}\in \mathbb{R}_{\geq0}^{n\times n}$, $\mathbf{C}=(c_{i,j})_{m\times n}\in \mathbb{R}_{\geq0}^{m\times n}$, and $\{\mathbf{E}_t\}_{t\in\mathbb{Z}}$ be an {\rm i.i.d.} sequence of matrices in $\mathbb{N}_0^{m\times n}$ with mutually independent components $e_{t,i,j}\sim Poisson(c_{i,j})$, $1\leq i\leq m,1\leq j\leq n$. Assume that,  $\forall t_1>t_2$,  $\{\mathbf{E}_{t_1}\}$ and $\{\mathbf{X}_{t_2}\}$  are independent of each other. Then the MINAR(1) model is denoted in the form
\begin{equation}\label{equ:2.2}
  \mathbf{X}_t = \mathbf{A}\circledast\mathbf{X}_{t-1}\circledast \mathbf{B}^{\top} + \mathbf{E}_t,~t\in\mathbb{Z},
  \end{equation}
where $\circledast$ operator independently over $t\in\mathbb{Z}$ and also independently of $\{\mathbf{E}_t\}_{t\in\mathbb{Z}}$.
}

Based on Definition \ref{def:2.4}, it can be obtained
\[\operatorname{vec}\left(\mathbf{A} \circledast\mathbf{Y} \circledast \mathbf{B}^{\top}\right)=(\mathbf{B} \otimes \mathbf{A}) * \operatorname{vec}(\mathbf{Y}),\]
where $\otimes$ denotes the matrix Kronecker product. Therefore, there are interchangeable relations between MINAR(1) and MGINAR(1) models:
\begin{equation}\label{equ:2.3}
\mathrm{vec}(\mathbf{X}_t)=(\mathbf{B}\otimes \mathbf{A})\ast \mathrm{vec}(\mathbf{X}_{t-1})+ \mathrm{vec}(\mathbf{E}_t).
\end{equation}

\begin{remark}\label{remark:1}Relying on model \eqref{equ:2.2}, there are also many extended models. For example,
\[
\mathbf{X}_t = \mathbf{A}_1\circledast\mathbf{X}_{t - 1}\circledast\mathbf{B}_1^\top+\mathbf{A}_2\circledast\mathbf{X}_{t - 1}\circledast\mathbf{B}_2^\top+\cdots+\mathbf{A}_d\circledast\mathbf{X}_{t - 1}\circledast\mathbf{B}_d^\top + \mathbf{E}_t
\]
and\[\mathbf{X}_{t}=\mathbf{A}_{1}\circledast\mathbf{X}_{t-1}\circledast\mathbf{B}_{1}^{\top}+\cdots+\mathbf{A}_{p}\circledast \mathbf{X}_{t-p}\circledast\mathbf{B}_{p}^{\top}+\mathbf{E}_{t}.\]
The first extended model is still a MINAR(1), but it involves more parameters, which helps to more fully capture the influences between rows and columns. However, there are more difficulties in terms of parameter estimation and parameter identifiability problems. The second extended model is a \(p\)th order MINAR model, which is a more general model. The proposed parameter estimation in subsequent section can be extended to the \(p\)th order MINAR model.
\end{remark}

\subsection{Model interpretation}\label{subsec:2.2}
This section provides a detailed further discussion of the coefficient matrices $\mathbf{A}$ and $\mathbf{B}$ in the MINAR(1) model. First, the left-hand matrix $\mathbf{A}$ in the model represents row-wise interactions, while the right-hand matrix $\mathbf{B}$ represents column-wise interactions. In order to specifically understand the action mechanism of matrices $\mathbf{A} $ and $\mathbf{B}$, we first ignore the information reflected by matrix $\mathbf{A}$, and assume that $\mathbf{A}=\mathbf{I}$, then \eqref{equ:2.2} becomes
\[\mathbf{X}_t=\mathbf{X}_{t-1}\circledast\mathbf{B}^{\top}+\mathbf{E}_t.\]

Consider the example in Figure \ref{fig:DQG}. The row of $\mathbf{X}_{t}$ represents different crime types at time $t$, and column represents different districts. Then the conditional expectations for the first column of $ \mathbf{X}_t$ can be expressed as
\[
\resizebox{1\hsize}{!}{$
\begin{pmatrix}
\text{THEFT} \\
\text{ROBBERY} \\
\text{ASSAULT}
\end{pmatrix}_{t}^{(10)} = 
\beta_{11}
\begin{pmatrix}
\text{THEFT} \\
\text{ROBBERY} \\
\text{ASSAULT}
\end{pmatrix}_{t-1}^{(10)} + 
\beta_{12}
\begin{pmatrix}
\text{THEFT} \\
\text{ROBBERY} \\
\text{ASSAULT}
\end{pmatrix}_{t-1}^{(11)} + 
\beta_{13}
\begin{pmatrix}
\text{THEFT} \\
\text{ROBBERY} \\
\text{ASSAULT}
\end{pmatrix}_{t-1}^{(15)} +
\begin{pmatrix}
\text{C\_THEFT} \\
\text{C\_ROBBERY} \\
\text{C\_ASSAULT}
\end{pmatrix}_{t}^{(10)},
$}
\]
where (10) denotes district 10, THEFT denotes the expected number of THEFT crime, \(\text{C\_THEFT}\) denotes the expected error in the number of THEFT crime, and the remaining entries are defined analogously. This means that the number of crime types occurring in one district at time $t$ is a linear combination of the number of the same crime type in all districts at time $t-1$, and this linear combination is identical for different crime types. Therefore, the coefficient matrix $\mathbf{B}$ represents the interactions between columns in the matrix data.

When $\mathbf{B}=\mathbf{I}$, the model becomes
\[\mathbf{X}_t=\mathbf{A}\circledast\mathbf{X}_{t-1}+\mathbf{E}_t.\]
Then the conditional expectation of the first row of the above model is
\[\resizebox{1\hsize}{!}{$
\begin{pmatrix}\text{THEFT}^{(10)}\\\text{THEFT}^{(11)}\\\text{THEFT}^{(15)}\end{pmatrix}_t^{\top}=
\alpha_{11}\begin{pmatrix}\text{THEFT}^{(10)}\\\text{THEFT}^{(11)}\\\text{THEFT}^{(15)}\end{pmatrix}_{t-1}^{\top}+
\alpha_{12}\begin{pmatrix}\text{ROBBERY}^{(10)}\\\text{ROBBERY}^{(11)}\\\text{ROBBERY}^{(15)}\end{pmatrix}_{t-1}^{\top}+
\alpha_{13}\begin{pmatrix}\text{ASSAULT}^{(10)}\\\text{ASSAULT}^{(11)}\\\text{ASSAULT}^{(15)}\end{pmatrix}_{t-1}^{\top}+
\begin{pmatrix}
\text{C\_THEFT}^{(10)} \\
\text{C\_THEFT}^{(11)} \\
\text{C\_THEFT}^{(15)}
\end{pmatrix}_{t}^{\top}
$},
\]
where \(\text{THEFT}^{(10)}\) represents expected number of THEFT crime in district 10, \(\text{C\_THEFT}^{(10)}\) represents the expected error in the number of THEFT crime in district 10, and the remaining entries follow analogously. It follows that the number of a certain crime type at time \( t \) is linearly related to the number of all crime types in the current district at time \( t-1 \). 

From another perspective, the first column of $\mathbf{X}_{t}$ can be expessed as
\[
\resizebox{0.8\hsize}{!}{$
\begin{pmatrix}
\text{THEFT\_obs}^{(10)} \\
\text{ROBBERY\_obs}^{(10)} \\
\text{ASSAULT\_obs}^{(10)}
\end{pmatrix}_t =
\mathbf{A}\circledast\begin{pmatrix}
\text{THEFT\_obs}^{(10)} \\
\text{ROBBERY\_obs}^{(10)} \\
\text{ASSAULT\_obs}^{(10)}
\end{pmatrix}_{t-1} +
\begin{pmatrix}
\text{E\_THEFT\_obs}^{(10)} \\
\text{E\_ROBBERY\_obs}^{(10)} \\
\text{E\_ASSAULT\_obs}^{(10)}
\end{pmatrix}_{t-1}
$},
\]
where \(\text{THEFT\_obs}^{(10)}\) represents the observed number of THEFT crime in district 10, and \(\text{E\_THEFT\_obs}^{(10)} \) represents the error number of THEFT crime in district1 10. It follows that, for each district, the number of crime types follows a MGINAR(1) model.

\subsection{Probabilistic properties of MINAR(1)}

Since the MINAR(1) model can be transformed into the MGINAR(1) model as shown in \eqref{equ:2.3}. In the MGINAR(1) model, if $\rho (\Phi)<1$, then the model is stationary and causal. Therefore, when $\rho (\mathbf{B}\otimes \mathbf{A})<1$, the MINAR(1) model is considered to be stationary and causal. Further, the following Proposition \ref{pro:2.1} holds.

{\proposition\label{pro:2.1} If $\rho(\mathbf{A})\rho(\mathbf{B})<1$, then the MINAR(1) model is stationary and causal.
}

The proof of this proposition can be refered to the proof of Proposition 1 in \cite{Chen2021}, as the two proofs are similar. 

If Proposition \ref{pro:2.1} holds, then the MINAR(1) model has the following vectorized causal representation (see \cite{Silva2005}) 
\[\operatorname{vec}(\mathbf{X}_t)=\sum_{k=0}^\infty\left(\mathbf{B}^k\otimes\mathbf{A}^k\right)\ast\operatorname{vec}(\mathbf{E}_{t-k}).\]

\section{Methodology}
\subsection{Reduce-rank MINAR model}
Compared with the MGINAR model, the MINAR model has achieved significant improvements in reducing the number of parameters and integrating data structural information. However, it still encounters parameter explosion when processing high-dimensional data. When the matrix dimension is high or there is noise, the coefficient matrix of the MINAR model obtained after dimensionality reduction and denoising becomes rank-deficient. In this case, fitting the data with a full-rank MINAR model would cause overfitting, leading to incorrect estimation of the coefficient matrix, and reduce  predictive performance. Drawing on the construction methodology of the RRMAR model proposed by \cite{Xiao2022}, we impose reduce-rank constraints on the coefficient matrix of the MINAR model, and propose a reduce-rank MINAR (RRMINAR) model. The key distinction of our model lies in its specialized design for integer-valued data, which is crucial for extracting key information from the data and remove noise. Specifically, the proposed RRMINAR model not only preserves the MINAR model's advantage of capturing row-column interactions but also overcomes the problems of overfitting and incorrect coefficient matrix estimation. In addition, the RRMINAR model reduces the number of parameters to be estimated, thereby lowering computational costs. This paper focuses on the RRMINAR(1) model.

{\definition\label{def:3.1} Under conditions of Definition 
\ref{def:2.4}, the RRMINAR(1) model is defined as
\begin{equation}\label{equ:3.1}
  \mathbf{X}_{t} = \mathbf{A} \circledast \mathbf{X}_{t - 1} \circledast \mathbf{B}^{\top} + \mathbf{E}_{t},
\end{equation}
where ${\rm rank}(\mathbf{A})=k_1<m$ and ${\rm rank}(\mathbf{B})=k_2<n$.
}

In the RRMINAR(1) model, to avoid unidentifiability caused by matrix scaling, we normalize matrix $\mathbf{A}$ such that $\|\mathbf{A}\|_F = 1$. Since the RRMINAR(1) model only introduces rank constraints to the MINAR(1) model without affecting parameter interpretation or stationarity conditions, Proposition \ref{pro:2.1} and the explanations of coefficient matrices for the MINAR(1) model in Section \ref{subsec:2.2}  remain equally applicable to the RRMINAR(1) model. These details will not be reiterated here. The number of model parameters in the RRMINAR(1) model depends on the rank of the coefficient matrix: The lower the rank, the fewer the parameters. When the data dimension is relatively high, the RRMINAR(1) model has fewer parameters than both the MGINAR(1) model and the MINAR(1) model.

\begin{remark}\label{remark:3.1} 
$\mathrm{(i)}$ The coefficient matrices of the MGINAR(1), MINAR(1) and RRMINAR(1) models contain \(m^2 n^2\), \(m^{2}+n^{2}\) and \(m^2+n^2-(m-k_1)^2-(n-k_2)^2\) parameters, respectively.   
$\mathrm{(ii)}$The MGINAR(1) model obscures interactions between different dimensions of matrix data, leading to loss of inherent structural information. In contrast, both the MINAR(1) and RRMINAR(1) models directly process matrix data, resulting in stronger model interpretability.
$\mathrm{(iii)}$ Similar to the MINAR model, the RRMINAR(1) model can be extended to the RRMINAR(p) model as follows
\[
\mathbf{X}_t = \sum_{j=1}^p \mathbf{A}_{j} \circledast \mathbf{X}_{t-j} \circledast \mathbf{B}_{j}^{\top} + \mathbf{E}_t,
\]
where both $\mathbf{A}_{j}$ and $\mathbf{B}_{j}$ are low-rank matrices. More lag orders can handle more information, but the parameters to be estimated also increase accordingly.
\end{remark}

\subsection{Iterative conditional least squares estimation} 
To obtain the parameter estimates of model \eqref{equ:3.1}, we need to solve the following constrained optimization problem. A direct idea is to find:
\begin{equation}\label{equ:Objective function}
\min_{\substack{
\mathbf{A}:~ \operatorname{rank}(\mathbf{A})=k_1 \\ 
\mathbf{B}:~ \operatorname{rank}(\mathbf{B})=k_2}} 
\operatorname{tr} \left\{ \sum_{t=2}^{T} (\mathbf{X}_t - \mathbf{A} \mathbf{X}_{t-1} \mathbf{B}^\top - \mathbf{C})(\mathbf{X}_t - \mathbf{A} \mathbf{X}_{t-1} \mathbf{B}^\top - \mathbf{C})^\top \right\},
\end{equation}
which is obviously a nonlinear solving problem, and it is difficult to find an explicit solution. Alternativel, we proposed the iterative conditional least squares estimation (ICLSE) for the RRMINAR(1) model, along with the corresponding algorithmic implementation procedure.
In this section, \(\mathbf{\hat{A}}_{\rm{RR.LS}}\) and \(\mathbf{\hat{B}}_{\rm{RR.LS}}\) are used to represent the ICLSE of $\mathbf{A}$ and $\mathbf{B}$ in \eqref{equ:3.1}. To obtain \(\mathbf{\hat{A}}_{\rm{RR.LS}}\), we first suppose \(\mathbf{B}\) and \(\mathbf{C}\) are given. Then we obtain \(\mathbf{\hat{A}}_{\rm{RR.LS}}\) by minimizing the trace of the sample conditional covariance matrix of the residuals:
\begin{equation}\label{equ:3.2}
\begin{aligned}
&\min_{\mathbf{A}:~ {\rm rank}(\mathbf{A})=k_1} \sum_{t=2}^{T} \left\| \mathbf{X}_t - \mathbb{E}(\mathbf{X}_t \mid \mathbf{X}_{t-1}) \right\|_F^2 \\
& = \min_{\mathbf{A}:~ {\rm rank}(\mathbf{A})=k_1} \operatorname{tr} \left\{ \sum_{t=2}^{T} (\mathbf{X}_t - \mathbf{A} \mathbf{X}_{t-1} \mathbf{B}^\top - \mathbf{C})(\mathbf{X}_t - \mathbf{A} \mathbf{X}_{t-1} \mathbf{B}^\top - \mathbf{C})^\top \right\} \\
& = \min_{\mathbf{A}:~ {\rm rank}(\mathbf{A})=k_1} \operatorname{tr} \Biggl\{ \left( \sum_{t=2}^{T} \mathbf{X}_t \mathbf{X}_t^\top - \mathbf{S}_{1yx} \mathbf{S}_{1xx}^{-1} \mathbf{S}_{1yx}^\top \right) + \sum_{t=2}^{T} (\mathbf{C}\mathbf{C}^\top - \mathbf{X}_t \mathbf{C}^\top - \mathbf{C}\mathbf{X}_t^\top) \\
& \quad ~ + (\mathbf{S}_{1yx} \mathbf{S}_{1xx}^{-\frac{1}{2}} - \mathbf{A} \mathbf{S}_{1xx}^{\frac{1}{2}})(\mathbf{S}_{1yx} \mathbf{S}_{1xx}^{-\frac{1}{2}} - \mathbf{A} \mathbf{S}_{1xx}^{\frac{1}{2}})^\top \Biggr\},
\end{aligned}
\end{equation}
where \(\mathbf{S}_{1xx}=\sum_{t = 2}^{T}\mathbf{X}_{t - 1}\mathbf{B}^{\top}\mathbf{B}\mathbf{X}_{t - 1}^{\top}\) and \(\mathbf{S}_{1yx}=\sum_{t = 2}^{T}(\mathbf{X}_{t}-\mathbf{C})\mathbf{B}\mathbf{X}_{t - 1}^{\top}\). Let \(\mathbf{U}_{1:k_1} = [\mathbf{u}_1, \mathbf{u}_2,... \mathbf{u}_{k_1}]\), where \(\mathbf{u}_j\) is the normalized eigenvector of \(\mathbf{S}_{1yx}\mathbf{S}_{1xx}^{-1}\mathbf{S}_{1yx}^{\top}\) corresponding to the \(j\)th eigenvalue. For simplify, we refer to \(\mathbf{u}_j\) as the \(j\)th normalized eigenvector. By Theorem 2.1 in \cite{Reinsel1998}, we have
\begin{equation*}
  \begin{aligned}
\mathbf{\hat{A}}_{\rm{RR.LS}}=\mathbf{U}_{1:k_1}\mathbf{U}_{1:k_1}^{\top}\mathbf{S}_{1yx}\mathbf{S}_{1xx}^{-1}.
  \end{aligned}
\end{equation*}
Similarly, by transposing \eqref{equ:3.1} and following \eqref{equ:3.2} with given \(\mathbf{A}\) and \(\mathbf{C}\), we also have
\begin{equation*}
\begin{aligned}
\mathbf{\hat{B}}_{\rm{RR.LS}}=\mathbf{\tilde U}_{1:k_2}\mathbf{\tilde U}_{1:k_2}^{\top}\mathbf{S}_{2yx}\mathbf{S}_{2xx}^{-1},
\end{aligned}
\end{equation*}
where \(\mathbf{S}_{2xx}=\sum_{t = 2}^{T}\mathbf{X}_{t - 1}^{\top}\mathbf{A}^{\top}\mathbf{A}\mathbf{X}_{t - 1}\), \(\mathbf{S}_{2yx}=\sum_{t = 2}^{T}(\mathbf{X}_{t}-\mathbf{C})^{\top}\mathbf{A}\mathbf{X}_{t - 1}\), \(\mathbf{\tilde U}_{1:k_2} = [\mathbf{\tilde{u}}_1, \mathbf{\tilde{u}}_2,... \mathbf{\tilde{u}}_{k_2}]\) and \(\mathbf{\tilde{u}}_j\) is the \(j\)th normalized eigenvector of \(\mathbf{S}_{2yx}\mathbf{S}_{2xx}^{-1}\mathbf{S}_{2yx}^{\top}\). 

To obtain \(\mathbf{\hat{C}}\), we minimizes the trace of the sample conditional covariance matrix of the residuals:
\[\sum_{t=2}^{T} \left\| \mathbf{X}_t - \mathbb{E}(\mathbf{X}_t \mid \mathbf{X}_{t-1}) \right\|_F^2 .\]
Specifically, taking the partial derivative of the above objective function with respect to \(\mathbf{C}\) and letting it equals zero, it holds that
\begin{equation}\label{equ:3.3}
\begin{aligned}
\mathbf{\hat{C}}=\frac{1}{T} \sum_{t = 2}^{T}(\mathbf{X}_{t}-\mathbf{A}\mathbf{X}_{t - 1}\mathbf{B}^{\top}).
\end{aligned}
\end{equation}

The ICLSE algorithm for the RRMINAR(1) model updates parameters in each iteration, and the process is presented in Algorithm \ref{alg:rrminar1}. If the initial values $\hat{\mathbf{B}}_{0}$ and $\hat{\mathbf{C}}_{0}$  satisfy certain conditions, the ICLSE estimation obtained by the Algorithm \ref{alg:rrminar1} has asymptotic normality. During the parameter updating process, the initial values $\hat{\mathbf{B}}_{0}$ and $\hat{\mathbf{C}}_{0}$ are derived from the projection estimation for the MINAR(1) model proposed by \cite{Xu2024}. 

\begin{algorithm}[H]
\caption{ICLSE for RRMINAR(1) Model}
\label{alg:rrminar1}
\KwIn{Initial values $\hat{\mathbf{B}}_{0}, \hat{\mathbf{C}}_{0}$; $\delta_T=o(T^{-1/2})$; integer-valued time series $\{\mathbf{X}_t\}_{t=1}^{T}$.}
\KwOut{Estimates: $\hat{\mathbf{A}}_{\rm RR.LS},~ \hat{\mathbf{B}}_{\rm RR.LS}, ~\hat{\mathbf{C}}$.}
\SetAlgoLined
$\mathbf{B}_{\text{prev}} \gets \hat{\mathbf{B}}_{0}$, $\mathbf{C}_{\text{prev}} \gets \hat{\mathbf{C}}_{0}$\;
\For{$r \gets 1$ \textbf{\rm \bf to} $2000$}{
    $\hat{\mathbf{S}}_{1yx} \gets \mathbf{0}$, $\hat{\mathbf{S}}_{1xx} \gets \mathbf{0}$\;
    \For{$t \gets 2$ \textbf{\rm \bf to} ${T}$}{
        $\hat{\mathbf{S}}_{1yx} \gets \hat{\mathbf{S}}_{1yx} + (\mathbf{X}_t - \mathbf{C}_{\text{prev}}) \cdot \mathbf{B}_{\text{prev}} \cdot \mathbf{X}_{t-1}^\top$\;
        $\hat{\mathbf{S}}_{1xx} \gets \hat{\mathbf{S}}_{1xx} + (\mathbf{X}_{t-1} \mathbf{B}_{\text{prev}}^\top \mathbf{B}_{\text{prev}} \mathbf{X}_{t-1}^\top)$\;
    }
    $\hat{\mathbf{U}}_{1:k_1} \gets \text{top-$k_1$ eigenvectors of } \hat{\mathbf{S}}_{1yx} \hat{\mathbf{S}}_{1xx}^{-1} \hat{\mathbf{S}}_{1yx}^\top$\;
    $\hat{\mathbf{A}}_r \gets \hat{\mathbf{U}}_{1:k_1} \hat{\mathbf{U}}_{1:k_1}^\top \hat{\mathbf{S}}_{1yx} \hat{\mathbf{S}}_{1xx}^{-1}$\;

    $\hat{\mathbf{S}}_{2yx} \gets \mathbf{0}$, $\hat{\mathbf{S}}_{2xx} \gets \mathbf{0}$\;
    \For{$t \gets 2$ \textbf{\rm \bf to} ${T}$}{
        $\hat{\mathbf{S}}_{2yx} \gets \hat{\mathbf{S}}_{2yx} + (\mathbf{X}_t - \mathbf{C}_{\text{prev}})^\top \cdot \hat{\mathbf{A}}_r \cdot \mathbf{X}_{t-1}$\;
        $\hat{\mathbf{S}}_{2xx} \gets \hat{\mathbf{S}}_{2xx} + (\mathbf{X}_{t-1}^\top \hat{\mathbf{A}}_r^\top \hat{\mathbf{A}}_r \mathbf{X}_{t-1})$\;
    }
    $\hat{\mathbf{U}}_{1:k_2} \gets \text{top-$k_2$ eigenvectors of } \hat{\mathbf{S}}_{2yx} \hat{\mathbf{S}}_{2xx}^{-1} \mathbf{S}_{2yx}^\top$\;
    $\hat{\mathbf{B}}_r \gets \hat{\mathbf{U}}_{1:k_2} \hat{\mathbf{U}}_{1:k_2}^\top \hat{\mathbf{S}}_{2yx} \hat{\mathbf{S}}_{2xx}^{-1}$\;

    $\hat{\mathbf{C}}_r \gets \mathbf{0}$\;
    \For{$t \gets 2$ \textbf{\rm to} ${T}$}{
        $\hat{\mathbf{C}}_r \gets \hat{\mathbf{C}}_r + (\mathbf{X}_t - \hat{\mathbf{A}}_r \mathbf{X}_{t-1} \hat{\mathbf{B}}_r^\top)$\;
    }
    $\hat{\mathbf{C}}_r \gets \hat{\mathbf{C}}_r / ({T}-1)$\;

    \If{$\|\hat{\mathbf{A}}_r - \mathbf{A}_{\mathrm{prev}}\|_F < \delta_t $ \textbf{\rm \bf and} $\|\hat{\mathbf{B}}_r - \mathbf{B}_{\mathrm{prev}}\|_F < \delta_t$ \textbf{\rm \bf and} $\|\hat{\mathbf{C}}_r - \mathbf{C}_{\mathrm{prev}}\|_F < \delta_t$}{
        \textbf{break}\;
    }
    $\mathbf{A}_{\text{prev}} \gets \hat{\mathbf{A}}_r$, $\mathbf{B}_{\text{prev}} \gets \hat{\mathbf{B}}_r$, $\mathbf{C}_{\text{prev}} \gets \hat{\mathbf{C}}_r$\;
}
\Return{
  $\hat{\mathbf{A}}_{\rm RR.LS} \gets \hat{\mathbf{A}}_r$, 
  $\hat{\mathbf{B}}_{\rm RR.LS} \gets \hat{\mathbf{B}}_r$, 
  $\hat{\mathbf{C}} \gets \hat{\mathbf{C}}_r$.}
\end{algorithm}

\begin{remark}
In the optimization process of the objective function \eqref{equ:Objective function}, the value of the objective function decreases monotonically  after each update, and the optimization process gradually converges the local optimal solution. Meanwhile, as iteration proceed, the error between the estimated parameters and the true parameters decreases gradually, and the final estimated parameter values converge stable solutions. During the iteration process, we naturally expect $\hat {\mathbf{A}}_ {r + 1} - \mathbf {A}_{r} \stackrel {\mathrm {p}} {\to} 0$ fast, $\mathbf{B}$ and $\mathbf{C}$ in the same way. So we set the stop condition \(\delta_T = o(T^{-1/2}) \). The rationality of this condition will be further discussed in the next section. This condition ensures the rationality of the algorithm, effectively balances computational efficiency and estimation accuracy, and avoids excessive iteration.

However, from a theoretical perspective, this optimization problem involves interaction between nonlinear constraints (such as low-rank constraints) and alternating least squares, making its theoretical analysis rather complex. In particular, the non-convexity of the objective function and the coupling relationship among parameters make the proof of global convergence challenging. Furthermore, the strict theoretical basis of the stopping conditions (such as the convergence efficiency and the error bound) still requires further in-depth research. Future work revolves around these issues.
\end{remark}

\section{Theory}
In this section, the consistency and asymptotic normality of the proposed ICLSE of RRMINAR(1) model obtained by Algorithm \ref{alg:rrminar1} will be presented. We first present some conclusions of the MINAR(1) model and the RRMINAR(1) model. These conclusions will be helpful in proving the consistency and asymptotic normality of the ICLSE of the RRMINAR(1) model.

\begin{theorem}\label{pro:4.1} 
Let $\{\mathbf{X}_t\}_{t\in\mathbb{Z}}$ be a matrix sequence generated from a $m\times n$-dimensional MINAR(1) model with coefficient matrices $\mathbf{A}\in \mathbb{R}_{\geq0}^{m\times m}$ and $\mathbf{B}\in \mathbb{R}_{\geq0}^{n\times n}$ such that $\rho(\mathbf{A})\rho(\mathbf{B})<1$ and innovation-parameter matrix $\mathbf{C}\in \mathbb{R}_{\geq0}^{m\times n}\backslash\{\mathbf{0}_{m\times n}\}$. Then
\[\mathbf{\mathbf{\Delta}}_t\coloneqq\mathbf{X}_t - \mathbf{A}\mathbf{X}_{t-1}\mathbf{B}^{\top} - \mathbf{C},~t\in\mathbb{Z},\]
define a matrix white noise sequence with concurrent correlations among its own entries, i.e., $\{\mathbf{\Delta}_t\}_{t\in\mathbb{Z}}$ is stationary with $\mathbb{E}(\mathbf{\Delta}_t) = \mathbf{0}_{m\times n}$ for $t\in\mathbb{Z}$ and
\[\mathbb{E}\big\{{\rm vec}(\mathbf{\Delta}_{t_1}){\rm vec}(\mathbf{\Delta}_{t_2})^{\top}\big\}=
\begin{cases}
\operatorname{diag}\big\{(\mathbf{1}_{mn \times mn}- \mathbf{B}\otimes \mathbf{A})^{-1}{\rm vec}(\mathbf{C})\big\},&t_1= t_2;\\
\mathbf{0}_{m\times n},&t_1\neq t_2.
\end{cases}
\]
Further, for any $t_1<t$, $\mathbf{\Delta}_t$ is uncorrelated with $\mathbf{X}_{t_1}$. 
\end{theorem}

The proof of Theorem \ref{pro:4.1} is given in Section \ref{Proof of Proposition 4.1} of Appendix. As already noted in \cite{Latour1997}, the INAR(p) model can be represented as a AR(p) model with white noise innovation terms. According to Theorem \ref{pro:4.1}, a $m\times n$ dimensional MINAR(1) model can be transformed into a $m\times n$ dimensional MAR(1) model with white noise errors. Then we have the following corollary.

\begin{corollary} \label{pro:4.2}Let the conditions of Theorem \ref{pro:4.1} hold and $\{\mathbf{\Delta}_t\}_{t\in\mathbb{Z}}$ be a matrix white noise sequence defined as in Theorem \ref{pro:4.1}. Then $\{\mathbf{X}_{t}\}_{t\in\mathbb{Z}}$ is a solution to the following stochastic difference equation
\[\mathbf{X}_t = \mathbf{A}\mathbf{X}_{t-1}\mathbf{B}^{\top} + \mathbf{C} + \mathbf{\Delta}_t,~t\in\mathbb{Z}.\]
\end{corollary}

By Corollary \ref{pro:4.2}, we have 
$$\operatorname{tr} \bigg\{\sum_{t=2}^{T} (\mathbf{X}_t - \mathbf{A} \mathbf{X}_{t-1} \mathbf{B}^\top - \mathbf{C})(\mathbf{X}_t - \mathbf{A} \mathbf{X}_{t-1} \mathbf{B}^\top - \mathbf{C})^\top\bigg\} = \operatorname{tr} \bigg(\sum_{t=2}^{T} \mathbf{\Delta}_t\mathbf{\Delta}_t^{\top}\bigg).$$ 
Hence, the trace of sample conditional covariance matrix for minimizing the residuals in \eqref{equ:3.2} is equal to the sum of squares of the minimized innovation terms of $\mathbf{\Delta}_t$.

Due to the alternating nature of the proposed ICLSE obtained by Algorithm \ref{alg:rrminar1}, the asymptotic properties of \(\mathbf{\hat{A}}_{\rm{\rm{RR.LS}}}\) and \(\mathbf{\hat{B}}_{\rm{RR.LS}}\) are intertwined. For instance, \(\mathbf{\hat{A}}_{\rm{RR.LS}}=\mathbf{\hat{U}_1}\mathbf{\hat{U}_1}^{\top}\mathbf{\hat{S}}_{1yx}\mathbf{\hat{S}}_{1xx}^{-1}\), where \(\mathbf{\hat{U}_1}\), \(\mathbf{\hat{S}}_{1yx}\), \(\mathbf{\hat{S}}_{1xx}^{-1}\) are substitutes obtained by replacing \(\mathbf{B}\) in \(\mathbf{U}_1\), \(\mathbf{S}_{1yx}\), \(\mathbf{S}_{1xx}^{-1}\) with \(\mathbf{\hat{B}_{\rm{RR.LS}}}\), respectively. A similar estimation applies to \(\mathbf{\hat{B}}_{\rm{RR.LS}}\). The following Theorem \ref{t:Iterative convergence} guarantees the validity of the ICLSE obtained by Algorithm \ref{alg:rrminar1}.

\begin{theorem}\label{t:Iterative convergence}
Consider a $m\times n$-dimensional RRMINAR(1) model with coefficient matrices $\mathbf{A}\in \mathbb{R}_{\geq0}^{m\times m}$ and $\mathbf{B}\in \mathbb{R}_{\geq0}^{n\times n}$ such that $\rho(\mathbf{A})\rho(\mathbf{B})<1$, ${\rm rank}(\mathbf{A})=k_1\le m$, ${\rm rank}(\mathbf{B})=k_2\le n$, and innovation-parameter matrix $\mathbf{C}\in \mathbb{R}_{\geq0}^{m\times n}\backslash\{\mathbf{0}_{m\times n}\}$. Let $\hat{\mathbf{B}}_0$ and $\hat{\mathbf{C}}_0$ be the initial values of Algorithm \ref{alg:rrminar1}. Assume that $\hat{\mathbf{B}}_0 = \mathbf{B} + O_\mathrm{p}(T^{-1/2})$ and $\hat{\mathbf{C}}_0 = \mathbf{C} + O_\mathrm{p}(T^{-1/2})$. Then, for the $r$th iteration values $\hat{\mathbf{A}}_r$, $\hat{\mathbf{B}}_r$ and $\hat{\mathbf{C}}_r$ obtained by Algorithm \ref{alg:rrminar1}, we have
$$\hat{\mathbf{A}}_r = \mathbf{A} + O_\mathrm{p}(T^{-1/2}),\quad \hat{\mathbf{B}}_r = \mathbf{B} + O_\mathrm{p}(T^{-1/2}),\quad \hat{\mathbf{C}}_r = \mathbf{C} + O_\mathrm{p}(T^{-1/2}),$$
where $r\ge1$.

\end{theorem}
The proof of Theorem \ref{t:Iterative convergence} is given in Section \ref{proof:iterative consistency} of Appendix.

Theorem \ref{t:Iterative convergence} shows that the iterative values during the preform of Algorithm \ref{alg:rrminar1} will converge to the true values in probability, respectively, if the initial values are properly selected. However, obtaining the asymptotic normality of the ICLSE of the RRMINAR(1) model obtained by Algorithm \ref{alg:rrminar1} is a highly challenging task. In fact, for the unconstrained MINAR(1) model, its conditional expectation least squares objective function is
\begin{equation}\label{MINAR_mubiao}
\sum_{t=2}^{T} \left\| \mathbf{X}_t - \mathbb{E}(\mathbf{X}_t \mid \mathbf{X}_{t-1}) \right\|_F^2 = \sum_{t=2}^{T} \| \mathbf{X}_t - \mathbf{A} \mathbf{X}_{t-1} \mathbf{B}^\top - \mathbf{C} \|_F^2.
\end{equation}  
Taking partial derivatives of \eqref{MINAR_mubiao} with respect to \(\mathbf{A}\), \(\mathbf{B}\), and \(\mathbf{C}\), respectively, we obtain the gradient conditions of the ICLSE of the unconstrained MINAR(1) model:  
\begin{equation}\label{MINAR_gradient}
\begin{cases}
\displaystyle\sum_{t=2}^{T} \mathbf{A} \mathbf{X}_{t-1} \mathbf{B}^\top \mathbf{B} \mathbf{X}_{t-1}^\top - \sum_{t=2}^{T} (\mathbf{X}_t - \mathbf{C}) \mathbf{B} \mathbf{X}_{t-1}^\top = \mathbf{0}_{m\times m}, \\
\displaystyle\sum_{t=2}^{T} \mathbf{B} \mathbf{X}_{t-1}^\top \mathbf{A}^\top \mathbf{A} \mathbf{X}_{t-1} - \sum_{t=2}^{T} (\mathbf{X}_t - \mathbf{C})^\top \mathbf{A} \mathbf{X}_{t-1} =\mathbf{0}_{n\times n}, \\
\displaystyle\sum_{t=2}^{T} \mathbf{C} - \sum_{t=2}^{T} (\mathbf{X}_t - \mathbf{A} \mathbf{X}_{t-1} \mathbf{B}^\top) = \mathbf{0}_{m\times n}.
\end{cases}
\end{equation}  
\cite{Chen2021} demonstrate that objective function \eqref{MINAR_mubiao} has at least one global minimum solution in \eqref{MINAR_gradient}. Condition (R) in \cite{Chen2021} ensures that the randomness of the error terms is sufficiently dispersed, avoiding the occurrence of flat regions in the objective function that lead to non-differentiable points. Thus, it ensures that the solutions of the gradient equations \eqref{MINAR_gradient} are isolated, that is, the number of local minima is finite. Moreover, condition (R) in \cite{Chen2021} guarantees that there exists a unique global minimum for \eqref{MINAR_gradient} with probability 1, providing a theoretical basis for establishing asymptotic normality of the ICLSE of the unconstrained MINAR(1) model. Whereas the low-rank constraints of the coefficient matrices make the solution space non-convex, and there may be multiple local optimal solutions, so that the uniqueness of the solution no longer holds, and it cannot be guaranteed that the iteration converges to the true parameter at a sufficiently fast rate. In addition, the low-rank constraint may violate the absolute continuity requirement in condition (R) in \cite{Chen2021}. 

To avoid the above issue and obtain the asymptotic normality of the ICLSE of the RRMINAR(1) model obtained by Algorithm \ref{alg:rrminar1}, we require the following technical assumption: $\exists~ r \in \mathbb{N}^+$ such that $|\hat{\mathbf{B}}_{r+1}-\hat{\mathbf{B}}_{r}|=o_\mathrm{p}(T^{-1/2})$ and $|\hat{\mathbf{C}}_{r+1}-\hat{\mathbf{C}}_{r}|=o_\mathrm{p}(T^{-1/2})$, which implies that the if-statement in Line 22 of Algorithm \ref{alg:rrminar1} will be executed within a finite number of steps. This technical assumption guarantees that the iteration converges to the true parameter at a sufficiently fast rate. While it is difficult to derive under the current conditions, and we will further explore the conditions for the validity of this technical assumption in subsequent research work. Recall that $\mathbf{U}_{1:k_1}$ and $\mathbf{U}_{1:k_2}$ are consist of top-$k_1$ eigenvectors of $\mathbf{S}_{1yx} \mathbf{S}_{1xx}^{-1} \mathbf{S}_{1yx}^\top$ and top-$k_2$ eigenvectors of $\mathbf{S}_{2yx} \mathbf{S}_{2xx}^{-1} \mathbf{S}_{2yx}^\top$, respectively. Let \(\mathbf{\Gamma}_1 = \mathbb{E}(\mathbf{X}_t^{\top}\mathbf{A}^{\top}\mathbf{A}\mathbf{X}_t)\), 
 \(\mathbf{\Gamma}_2 = \mathbb{E}(\mathbf{X}_t\mathbf{B}^{\top}\mathbf{B}\mathbf{X}_t^{\top})\), \(\mathbf{P}_1 = \mathbf{U}_{1:k_1}\mathbf{U}_{1:k_1}^{\top}\), \(\mathbf{P}_2 = \mathbf{U}_{1:k_2}\mathbf{U}_{1:k_2}^{\top}\). Define
\begin{equation*}
\mathbf{Q}_t = 
\begin{pmatrix}
(\mathbf{X}_t\mathbf{B}^{\top}) \otimes \mathbf{P}_1 + [\Gamma_2\mathbf{A}^{\top}(\mathbf{A}\Gamma_2\mathbf{A}^{\top})^+\mathbf{A}\mathbf{X}_t\mathbf{B}^{\top}]\otimes (\mathbf{I} - \mathbf{P}_1)\\
\mathbf{P}_2\otimes (\mathbf{X}_t^{\top}\mathbf{A}^{\top}) + (\mathbf{I} - \mathbf{P}_2)\otimes [\Gamma_1\mathbf{B}^{\top}(\mathbf{B}\Gamma_1\mathbf{B}^{\top})^+\mathbf{B}\mathbf{X}_t^{\top}\mathbf{A}^{\top}]\\
\mathbf{I}_n\otimes \mathbf{I}_m
\end{pmatrix}
\end{equation*}
and
\begin{equation*}
\mathbf{W}_t = 
\begin{pmatrix}
(\mathbf{X}_t\mathbf{B}^{\top}\mathbf{B}\mathbf{X}_t^{\top}) \otimes \mathbf{I}_m & (\mathbf{X}_t\mathbf{B}^{\top})\otimes(\mathbf{A}\mathbf{X}_t) & (\mathbf{X}_t\mathbf{B}^{\top})\otimes \mathbf{P}_1\\
(\mathbf{B}\mathbf{X}_t)\otimes(\mathbf{X}_t\mathbf{A}^{\top}) & \mathbf{I}_n\otimes (\mathbf{X}_t\mathbf{A}^{\top}\mathbf{A}\mathbf{X}_t) & \mathbf{P}_2\otimes (\mathbf{X}_t^{\top}\mathbf{A}^{\top})\\
(\mathbf{B}\mathbf{X}_t^{\top})\otimes \mathbf{I}_m & \mathbf{I}_n\otimes (\mathbf{A}\mathbf{X}_t) & \mathbf{I}_n\otimes \mathbf{I}_m
\end{pmatrix}.
\end{equation*}
According to \cite{Chen2021} , matrix $\mathbb{E}(\mathbf{W}_t)$ may be singular. To address this, we define \(\mathbf{H}=\mathbb{E}(\mathbf{W}_t)+\gamma_1\gamma_1^{\top}\) with \(\gamma_1\coloneqq[\mathrm{vec}(\mathbf{A})^{\top},\mathbf{0}^{\top}]^{\top}\in \mathbb{R}^{m^2+n^2}\). Although $\mathbf{H}$ is generally invertible, providing a rigorous proof is highly challenging and goes beyond the core focus of this study. Therefore, we adopt the invertibility of $\mathbf{H}$ as a technical assumption and defer its analysis to future work.

Under the above technical assumptions and notation, Theorem \ref{t:4.2} gives the asymptotic normality of the ICLSE of the RRMINAR(1) model obtained by Algorithm \ref{alg:rrminar1}. 

\begin{theorem} \label{t:4.2} 
Let the conditions of Theorem \ref{t:Iterative convergence} hold. Assume that $\exists~ r \in \mathbb{N}^+$ such that $|\hat{\mathbf{B}}_{r+1}-\hat{\mathbf{B}}_{r}|=o_\mathrm{p}(T^{-1/2})$ and $|\hat{\mathbf{C}}_{r+1}-\hat{\mathbf{C}}_{r}|=o_\mathrm{p}(T^{-1/2})$, and \(\Sigma_{\mathbf{\Delta}}:={\rm Cov}\{\mathrm{vec}(\mathbf{\Delta}_t)\}\) is positive definite. Then,
\begin{equation*}
\sqrt{T}
\begin{pmatrix}
\text{\rm{vec}}(\mathbf{\hat{A}}_{\rm{RR.LS}}-\mathbf{A})\\
\text{\rm{vec}}(\mathbf{\hat{B}}_{\rm{RR.LS}}^{\top}-\mathbf{B}^{\top})\\
\text{\rm{vec}}(\mathbf{\hat{C}}_{\rm{RR.LS}}-\mathbf{C})
\end{pmatrix}
\Rightarrow N(0,\mathbf{\Xi_2}),
\end{equation*} 
where \(\mathbf{\Xi}_2 = \mathbf{H}^{-1}\mathbb{E}(\mathbf{Q}_t\Sigma_{\mathbf{\Delta}}\mathbf{Q}_t^{\top})\mathbf{H}^{-1}\), and $\Rightarrow$ denotes convergence in distribution.
\end{theorem}

The proof of Theorem \ref{t:4.2} is given in Section \ref{Proof of Theorem 4.2} of the Appendix.
\section{Experiments}
\subsection{Simulation}
In this section, we describe the performance of the least squares estimation (LSE) for the MGINAR(1) model, as well as the ICLSE for the MINAR(1) and the RRMINAR(1) models, through detailed simulation. Our simulation are carried out under various settings of matrix dimensions \(m\) and \(n\), ranks \(k_1\) and \(k_2\), and time series length \(T\).

For given dimensions \(m\) and \(n\), ranks \({k_1}\) and \(k_2\), we generate random matrices \(\mathbf{A}\) and \(\mathbf{B}\) such that $\rho(\mathbf{B})\rho(\mathbf{A})<1$ and \(\|\mathbf{A}\|_F = 1\). According to \eqref{equ:3.1}, we can generate a numerical simulation observation sequence \(\mathbf {X}_{t}\). In multiple repeated simulations under the same dimensions and ranks, the coefficient matrices \(\mathbf{A}\) and \(\mathbf{B}\) remain fixed. Specifically, we consider three distinct scenarios to examine the finite-sample performance of the estimators. The three settings are as follows:
\begin{itemize}
  \item[$\bullet$] \textbf{Setting I:} The covariance matrix $\text{Cov}(\text{vec}(\mathbf{E}_t))=\Sigma $ is set to $\Sigma = \mathbf{I}$.
  \item[$\bullet$] \textbf{Setting II:} The covariance matrix $\text{Cov}(\text{vec}(\mathbf{E}_t))=\Sigma$ is diagonal, with its elements following a uniform distribution on (0,1).
  \item[$\bullet$] \textbf{Setting III:} The covariance matrix Cov $(\operatorname{vec}(\mathbf{E}_t))=\Sigma=\Sigma_c \otimes \Sigma_r$ is randomly generated according to $\Sigma_c=\mathbf{Q}\Lambda \mathbf{Q}^{\top}$, where the eigenvalues in the diagonal matrix $\Lambda$ are the absolute values of i.i.d. standard normal random variates, and the eigenvector matrix $\mathbf{Q}$ is a random orthonormal matrix. $\Sigma_r$ is formed following an identical procedure.
  \end{itemize}
\begin{figure}[htbp]
  \centering
  \includegraphics[scale=0.5]{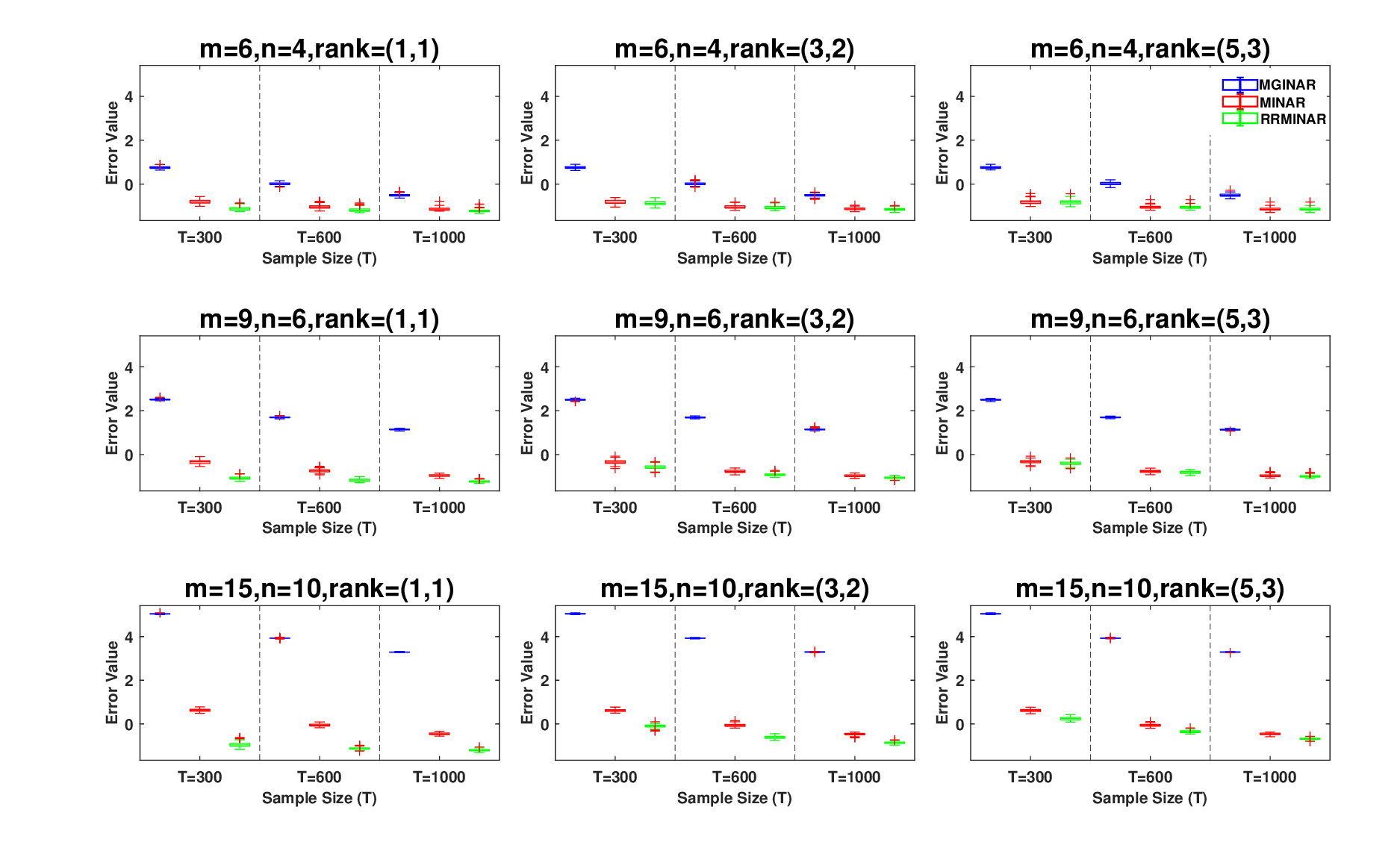}
  \caption{Comparison of the estimation errors for the MGINAR(1) (LSE), MINAR(1) (ICLSE), and RRMINAR(1) (ICLSE) models under Setting I. The three panels in each sub-figure correspond to sample sizes of 300, 600, and 1000, respectively.}
  \label{fig:set1_exp1}
  \end{figure}

For different settings, we first compare the estimation errors of the LSE for the MGINAR(1) model, the ICLSE for the MINAR(1) and RRMINAR(1) models under different time series lengths. The estimation errors are calculated using the following formulas:
\begin{equation}\label{equ:5.1}
\mathrm{MINAR}(1),~\mathrm{RRMINAR}(1){:}~\log\| \hat{\mathbf{B}}\otimes\hat{\mathbf{A}}-\mathbf{B}\otimes \mathbf{A}\| _F^2;~\mathrm{MGINAR}(1){:}~\log\|\hat{\Phi}-\mathbf{B}\otimes \mathbf{A}\|_F^2.
\end{equation} 

The time series lengths are set to \( T = 300, 600, 1000 \). The data dimensions are chosen as \( (m, n) = (6, 4), (9, 6), (15, 10) \) with the corresponding ranks \( (k_1, k_2) = (1, 1), (3, 2), (5, 3) \), representing low, medium and high rank levels, respectively. The simulation is repeated 100 times, and the distributions of the estimation errors for the coefficient matrices of the MGINAR(1), MINAR(1) and RRMINAR(1) models under Setting I are shown in Figure \ref{fig:set1_exp1}. 

As shown in Figure \ref{fig:set1_exp1}, under Setting I with the same dimension and rank constraint, the estimation errors of all three models decrease as the sample size increases, reflecting the consistency of the estimators. Further, under the same sample size and dimension, the estimation error of the RRMINAR(1) model approaches that of the MINAR(1) model as the rank increases. In contrast, the estimation error of the RRMINAR(1) model outperforms the MINAR(1) model when the rank decreases. Moreover, under the same sample size and rank constraint, the higher the dimension, the greater the errors of the MGINAR(1) and MINAR(1) models. While the proposed RRMINAR(1) model exhibits high robustness to the increase in dimensionality. In addition, the LSE of the MGINAR(1) model shows the worst performance. It should be noted that the ICLSE of the RRMINAR(1) model shows relatively better estimation performance when the ranks of the coefficient matrices are low.
\begin{figure}[htbp]
  \centering
  \includegraphics[scale=0.5]{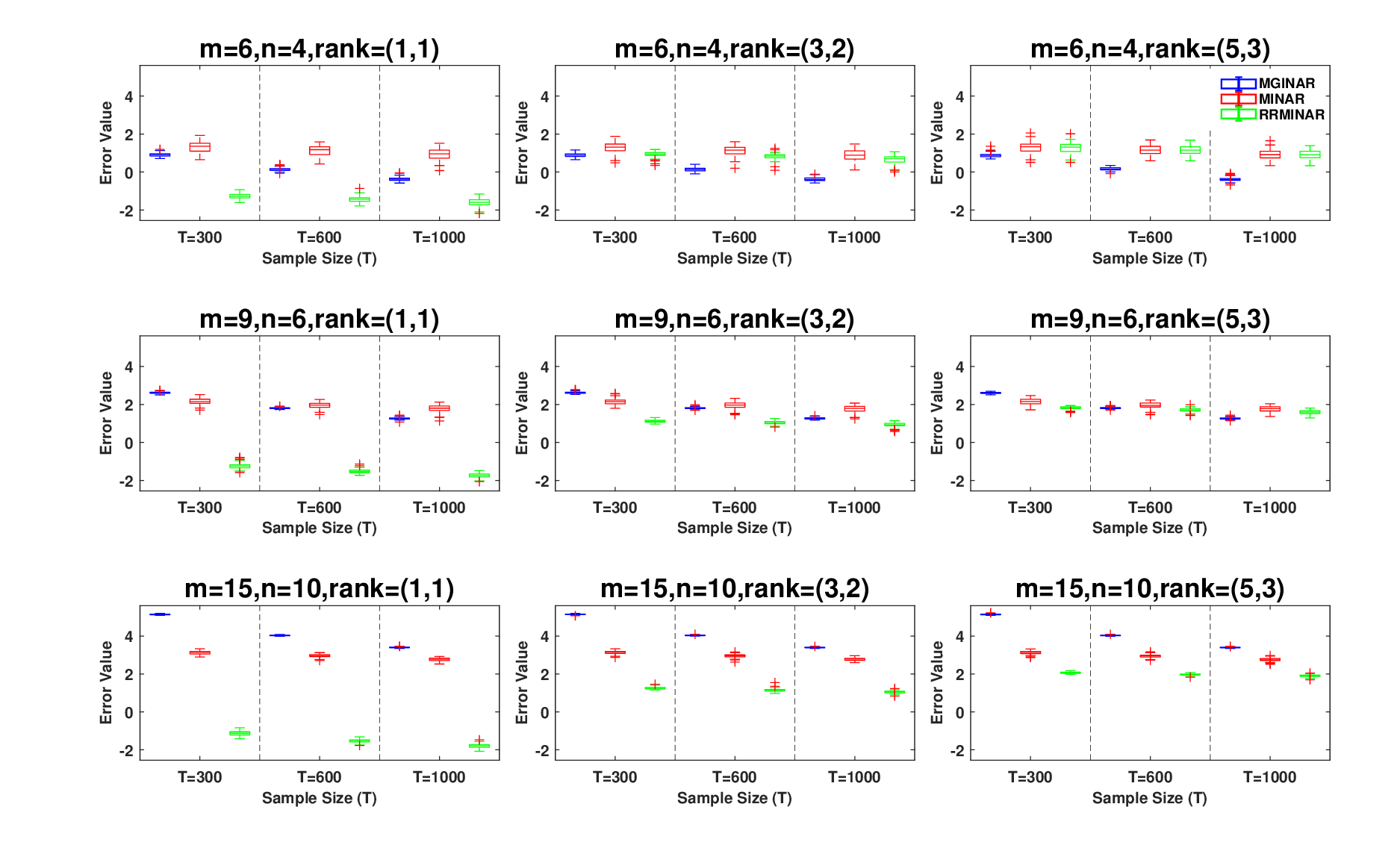}
  \caption{Comparison of the estimation errors for the MGINAR(1) (LSE), MINAR(1) (ICLSE), and RRMINAR(1) (ICLSE) models under Setting II. The three panels in each sub-figure correspond to sample sizes of 300, 600, and 1000, respectively.}
  \label{fig:set2_exp1}
  \end{figure}

Figure \ref{fig:set2_exp1} shows that under Setting II, with the same dimension and rank constraint, the estimation errors decrease as the sample size increases. Under fixed dimension and sample size, as the rank gradually increases, the RRMINAR(1) model's estimation error approaches that of the MINAR(1) model. Unlike the results of Setting I, when the covariance structure of the innovations becomes more complex and the dimension is small, the estimation error of the MGINAR(1) model is smaller than those of the RRMINAR(1) and MINAR(1) models. It is possible that the MGINAR(1) model contains more parameters and can thus capture more information when the coefficient matrix approaches full rank. While maintaining a fixed rank and increasing the dimension, this relationship reverses. Furthermore, under the same rank constraint and sample size, as the dimension increases, the estimation error of the MGINAR(1) model becomes larger and the variation range of the estimation error is faster than that of the other two models, which implies that the MGINAR(1) model is unstable and less suitable for high-dimensional data. In addition, the results under Setting III (see Appendix \ref{fig:appendix}) are similar to those under Setting II. Overall, these experiments show that under the high-rank setting, the RRMINAR(1) model and the MINAR(1) model are isomorphic, whereas under low-rank settings, the RRMINAR(1) model shows better performance.
\begin{figure}[!htbp]
  \centering
  \includegraphics[width=1\textwidth]{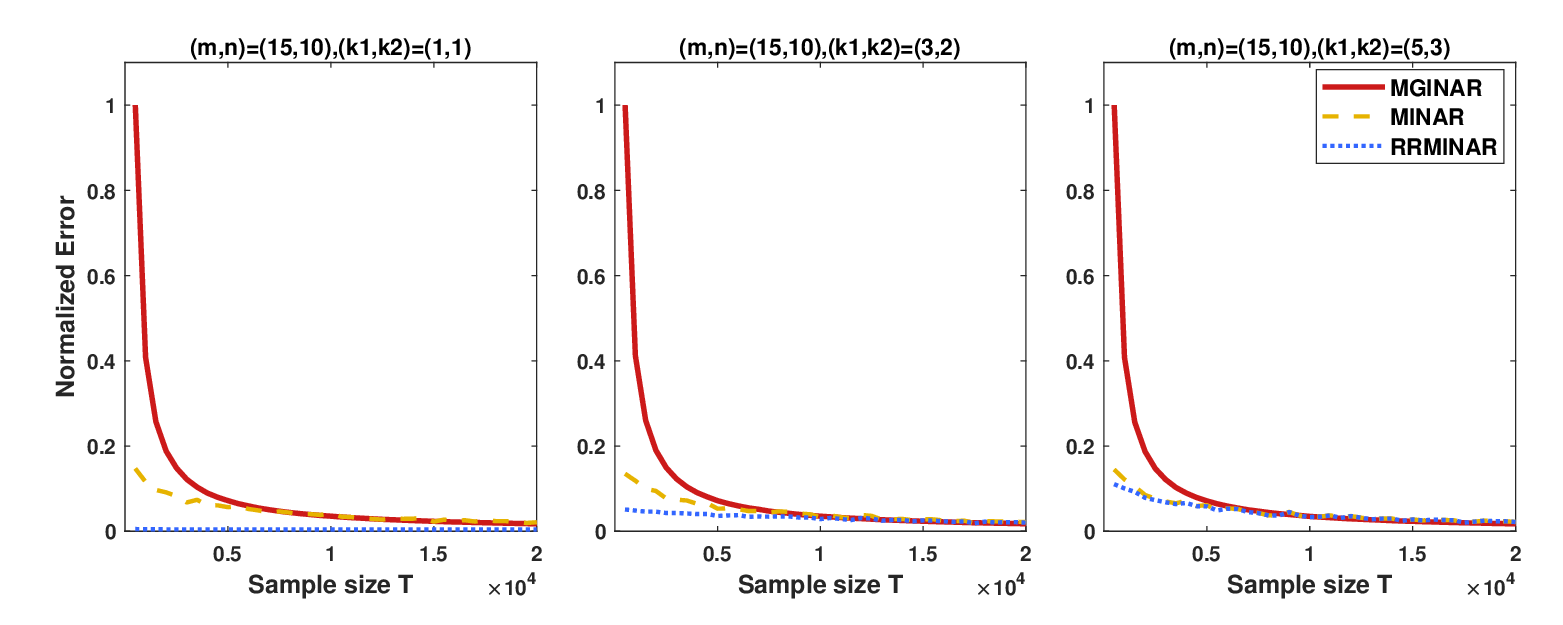}
  \caption{Convergence efficiency of coefficient matrix \(\mathbf{B}\otimes \mathbf{A}\) under Setting III.}
  \label{fig:efficiency}
\end{figure}

To study the convergence efficiency of the RRMINAR(1), MINAR(1) and MGINAR(1) models, we plot the variations in the mean values of estimation error \eqref{equ:5.1} for different models under different settings and different $T$ over 100 repetitions. Specifically, to highlight the relative performance trends of different models, the estimation errors are normalized by the global maximum across all methods and sample sizes as follows:
\[
\tilde{e}_{model}(T)=\frac{\bar{e}_{model}(T)}{\max\limits_{{model}, T}\{\bar{e}_{model}(T)\}},
\]
where \(\bar{e}_{model}(T) = \frac{1}{L} \sum_{r = 1}^{L} \exp\{e_{model}^{(r)}(T)\}\), $model\in$ \{{MGINAR(1), MINAR(1), RRMINAR(1)}\}, $L$ represents the number of repeated experiments and $e_{model}^{(r)}(T)$ represents the $r$th estimation errors under different sample size $T \in \{500,1000,...,20000\}$. The results based on Setting III with dimension $(m, n) =(15,10)$ are shown in Figure \ref{fig:efficiency}. In Figure \ref{fig:efficiency}, the red, yellow and blue lines in each panel represent the MGINAR(1), MINAR(1) and RRMINAR(1) models' results, respectively. As can be seen from Figure \ref{fig:efficiency}, the normalized estimation errors of all three models exhibit a stable decreasing trend as $T$ increases. Under all three rank constraints, the MGINAR(1) model shows the slowest convergence efficiency, while the RRMINAR(1) model converges the fastest. When the rank is (1,1), the RRMINAR(1) model has the highest convergence efficiency. As the rank increases, the convergence efficiency of the MINAR(1) model gradually approaches that of the RRMINAR(1) model, which is consistent with the previous results. Results under other settings are similar and provided in Section \ref{fig:appendix} of the Appendix.

\subsection{The $C_p$ criterion}
In the RRMINAR model, the selection of the ranks of coefficient matrices $\mathbf{A}$ and $\mathbf{B}$ is involved. This paper selects the ranks of coefficient matrices $\mathbf{A}$ and $\mathbf{B}$ based on the $C_p$ criterion proposed by \cite{Mallows2000}. The calculation formula of $C_p$ criterion is as follows:
\[
C_p=\frac{RSS_s}{\sigma^2}-(n-2k),
\]
\[
{\sigma}^2=\sum_{i=1}^n\frac{\|\mathbf{X}_t-\mathbf{\hat{X}}_t\|_F^2}{n-p},
\]
where, $RSS_s$ represents the residual sum of squares of sub-model $S$, ${\sigma}^2$ represents the mean squared error of the full model, $n$ is the number of observed data, $p$ denotes the number of parameters in the full model, and $k$ represents the number of parameters in the sub-model. For the RRMINAR model in this paper, $RSS_s$ corresponds to the sum of the squared residuals of the coefficient matrices $\mathbf{A}$ and $\mathbf{B}$ under different ranks, while ${\sigma}^2$ corresponds to the mean squared error of the unconstrained MINAR model. 

The $C_p$ criterion reflects the degree of fit of the sub-model to the data through the term ${RSS_s \sigma^{-2}}$, while incorporating a penalty term to balance goodness of fit and model complexity. To obtain a more accurate estimate of the rank, we adopt a method combining rolling fitting and the $C_p$ criterion to select the ranks of the coefficient matrices $\mathbf{A}$ and $\mathbf{B}$ in the RRMINAR model. Specifically, in each experiment, the data are divided into three segments respectively, and the $C_p$ values for different sub-models are computed for each segment. The average $C_p$ across segments is used as the final criterion for model selection. 

To evaluate the effectiveness of the proposed selection approach, simulation experiments are conducted under each of the three settings described earlier. We report the success rates of rank selection over 20 repeated experiments under varying dimensions, sample sizes $T$, and true ranks. Given that the initial estimate for the RRMINAR model is derived from the projected estimate of the MGINAR model, the experimental results are shown in Tables \ref{tab:rank_selection_6x4} and \ref{tab:rank_selection_3x3}.

\begin{table}[htbp]
\centering
\caption{Proportion of correct rank selection by $C_p$ criterion for dimension $(m,n)=(6,4)$}
\label{tab:rank_selection_6x4}
\begin{tabular}{cccccc}
\toprule
\multirow{2}{*}{Setting} & \multirow{2}{*}{Sample Size (T)} & \multicolumn{3}{c}{True Rank $(r_1,r_2)$} \\
\cmidrule(lr){3-5}
 & & $(1,1)$ & $(3,2)$ & $(4,3)$ \\
\midrule
\multirow{4}{*}{I} 
 & 600 & (0.65, 0.60) & (0.10, 0.70) & (0.00, 0.25) \\
 & 1000 & (0.75, 0.80) & (0.30, 0.80) & (0.10, 0.70) \\
 & 2000 & (0.95, 0.90) & (0.75, 0.80) & (0.35, 0.85) \\
 & 5000 & (0.85, 0.85) & (0.85, 0.85) & (0.80, 0.75) \\
\midrule
\multirow{4}{*}{II}
 & 600 & (0.60, 0.70) & (0.25, 0.95) & (0.00, 0.55) \\
 & 1000 & (0.75, 0.80) & (0.45, 0.85) & (0.25, 0.50) \\
 & 2000 & (0.65, 0.80) & (0.75, 0.90) & (0.40, 0.75) \\
 & 5000 & (0.75, 0.90) & (0.80, 0.90) & (0.75, 0.90) \\
\midrule
\multirow{4}{*}{III}
 & 600 & (0.55, 0.50) & (0.30, 0.70) & (0.05, 0.55) \\
 & 1000 & (0.55, 0.70) & (0.50, 0.80) & (0.15, 0.60) \\
 & 2000 & (0.90, 0.60) & (0.75, 0.85) & (0.50, 0.75) \\
 & 5000 & (0.75, 0.75) & (0.80, 0.95) & (0.60, 0.65) \\
\bottomrule
\end{tabular}
\end{table}

\begin{table}[htbp]
\centering
\caption{Proportion of correct rank selection by $C_p$ criterion for dimension $(m,n)=(3,3)$}
\label{tab:rank_selection_3x3}
\begin{tabular}{ccccc}
\toprule
\multirow{2}{*}{Setting} & \multirow{2}{*}{Sample Size (T)} & \multicolumn{3}{c}{True Rank $(r_1,r_2)$} \\
\cmidrule(lr){3-5}
 & & $(1,1)$ & $(1,2)$ & $(2,1)$ \\
\midrule
\multirow{3}{*}{I} 
 & $100$ & (0.85, 0.80) & (0.90, 0.25) & (0.30, 0.75) \\
 & $300$ & (0.90, 0.80) & (0.70, 0.45) & (0.20, 0.75) \\
 & $1000$ & (0.95, 0.95) & (0.80, 0.55) & (0.75, 0.90) \\
\midrule
\multirow{3}{*}{II}
 & $100$ & (0.75, 0.80) & (0.85, 0.40) & (0.30, 0.55) \\
 & $300$ & (0.70, 0.50) & (0.90, 0.60) & (0.25, 0.80) \\
 & $1000$ & (0.85, 0.75) & (0.85, 0.70) & (0.65, 0.75) \\
\midrule
\multirow{3}{*}{III}
 & $100$ & (0.70, 0.65) & (0.70, 0.25) & (0.10, 0.80) \\
 & $300$ & (0.85, 0.75) & (0.80, 0.55) & (0.40, 0.70) \\
 & $1000$ & (0.75, 0.95) & (0.85, 0.70) & (0.65, 0.90) \\
\bottomrule
\end{tabular}
\end{table}

The experimental results show that under different dimensions $(m,n)$ and rank combinations $(r_1,r_2)$, the rank selection success rate of the $C_p$ criterion exhibits significant differences. As the time series length $T$ increases, the accuracy of rank identification generally improves. For instance, in Table \ref{tab:rank_selection_6x4}, when $(r_1,r_2) = (1,1)$, the success rate under Setting I increases significantly from $(0.65,0.60)$ when $T=600$ to $(0.95,0.90)$ when $T=2000$. The identification accuracy for the low-rank combination $(1,1)$ is the highest, while the identification of the high-rank combination $(4,3)$ is the most challenging. Specifically, under Setting III with $T=600$, the success rate for \((4, 3)\) is only $(0.05,0.55)$. Obvious performance differences are observed across the three settings: Setting I performs best overall, while Setting III shows greater instability. Furthermore, the recognition effect of the dimension $(3,3)$ is generally better than that of the dimension $(6,4)$, indicating that low-dimensional problems are easier to solve. When the sample size is low, the success rate of correct selection for low-rank is higher, while that for high-rank is lower. However, as the sample size increases, the accuracy improved across all settings. These results provide an important basis for further optimizing the rank selection method and parameter settings. In addition, considering the substantial variation in the number of parameters of the MGINAR model across different dimensions, a larger sample size for dimensions $(6,4)$ is required to achieve the same proportion of correct rank selection as for dimensions $(3,3)$. Consequently, the sample sizes presented in Table \ref{tab:rank_selection_6x4} are greater than those in Table \ref{tab:rank_selection_3x3}. Other order determination methods such as AIC and BIC need further research.

\subsection{Crime data analysis}
\begin{figure}[htbp]
  \centering
  \includegraphics[scale=0.8]{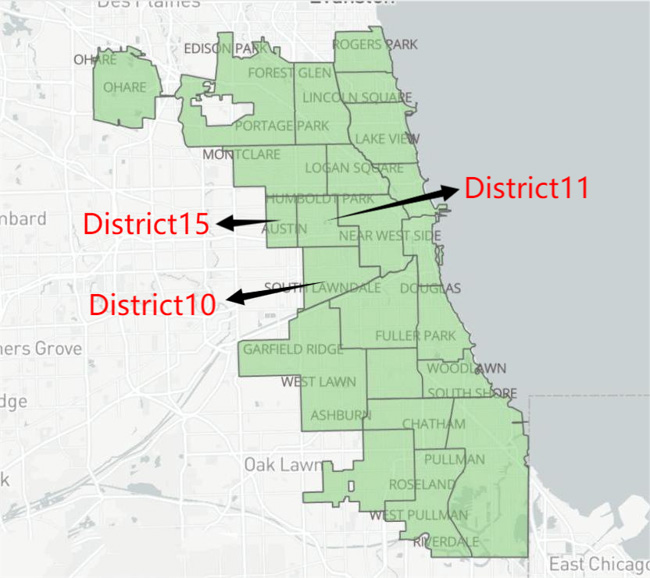}
  \caption{A map of the division of the district of Chicago.}
  \label{fig:Chicago}
\end{figure}

According to the content of criminal psychology, there are interactions among different types of crime in different districts. The interactions may be led to by the spatial proximity (criminal activities often cluster geographically, that is, the so-called "crime hotspots"), the socio-economic factors (socio-economic problems such as poverty, unemployment, low educational levels, and lack of community resources may lead to certain types of crime being more prevalent in specific districts), the criminal chain reaction (there is a certain correlation among certain types of crime, such as prostitution, gambling, and drug-related crime), and the urban crime prevention measures. To study the interactions between different districts and different types of crime, we selected three districts in Chicago (divided into district10, district11 and district15 according to police districts) as the research objects. For detailed division rules, please refer to https://www.chicagopolice.org/statistics-data/data-dashboards/sentiment-dashboard/. Not all crime data in various districts are adapted to the RRMINAR model. This is due to the fact that the estimates for parameter matrices of the RRMINAR(1) model derived from Algorithm \ref{alg:rrminar1} can contain negative values with certain data. After debugging, we found that the crime data in districts 10, 11, and 15 are adapted to the model. We selected the daily counts of three types of crime (i.e., THEFT, ASSAULT and ROBBERY) within the three districts between 2010.1.1 and 2011.2.19, with a total of 415 observations (i.e., $T=415$). The data can be downloaded for free from https://aistudio.baidu.com/datasetdetail/53627. The three types of crime data from three districts are shown in Figure \ref{fig:DQG}. We show the geographical locations of the three districts in Chicago in Figure \ref{fig:Chicago}.

\begin{figure}[!htbp]
  \centering
  \includegraphics[width=1\textwidth]{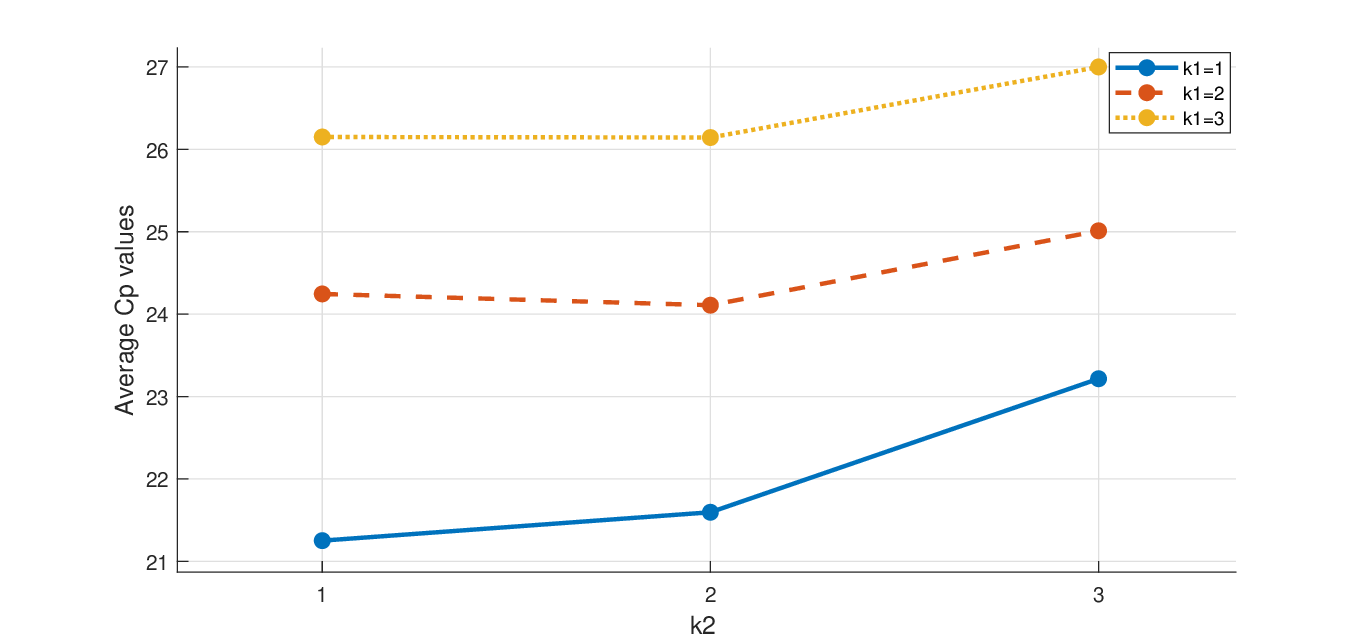}
  \caption{Results of the $C_p$ criterion for selecting the ranks of $\mathbf{A}$ and $\mathbf{B}$.}
  \label{fig:rank change}
\end{figure}

At any time, the data on crime types across different districts takes in the form of a $3\times3$ matrix. The first 355 observation matrices are selected as the training set, and the remaining 60 observation matrices are used as the test set. The ranks of $\mathbf{A}$ and $\mathbf{B}$ in the RRMINAR(1) model are selected using the $C_p$ criterion combined with non-overlapping rolling cross-validation. The training set is evenly divided into three segments, and the $C_p$ values for different rank combinations are computed on each segment respectively. Figure \ref {fig:rank change} shows the average $C_p$ values corresponding to each rank combination. Based on \ref{fig:rank change}, the ranks of $\mathbf{A}$ and $\mathbf{B}$ are set to $k_1=1$ and $k_2 = 1$, respectively. 
  
Based on the selected ranks and the observation matrices in the training dataset, the MGINAR(1), MINAR(1) and RRMINAR(1) models are established for one-step-prediction, respectively. The training and prediction performance is evaluated using the following four error metrics, with any zero denominator in $\rm E_3$ or $\rm E_4$ replaced by a value of 1:
\begin{align*}
\mathrm{E}_1 &= \sum_{t=1}^{T} \sqrt{ \sum_{i=1}^{m} \sum_{j=1}^{n} \left( \mathbf{X}_{ij,t} - \hat{\mathbf{X}}_{ij,t} \right)^2 },\\
\mathrm{E}_2 &= \sqrt{ \frac{1}{m n T} \sum_{t=1}^{T} \sum_{i=1}^{m} \sum_{j=1}^{n} \left( \mathbf{X}_{ij,t} - \hat{\mathbf{X}}_{ij,t} \right)^2 },\\
\mathrm{E}_3 &= \sqrt{\frac{1}{m nT} \sum_{t=1}^{T} \sum_{i=1}^{m} \sum_{j=1}^{n} \frac{ \left( \mathbf{X}_{ij,t} - \hat{\mathbf{X}}_{ij,t} \right)^2 }{ \mathbf{X}_{ij,t}^2 }},\\
\mathrm{E}_4 &= \frac{1}{m  n T} \sum_{t=1}^{T} \sum_{i=1}^{m} \sum_{j=1}^{n} \frac{ \left| \mathbf{X}_{ij,t} - \hat{\mathbf{X}}_{ij,t} \right| }{  \overline{\mathbf{X}}_{ij}}.
\end{align*}
Tables \ref{tab:within_sample_fitting} and \ref{tab:out_of_sample_fitting} summarize the in-sample and out-of-sample average prediction errors of the daily crime matrix data under seven models. The seven models compared are as follows:
\begin{enumerate}[label=(\roman*)]
  \item MGINAR\_row(1): Fit three MGINAR(1) models to each row of the observation matrix series, i.e., fit them to the occurrence numbers of THEFT, ROBBERY, and ASSAULT in different districts, respectively.
  \item MGINAR\_column(1): Fit three MGINAR(1) models to each column of the observational matrix series; each model corresponds to a different type of criminal event and is fitted to the occurrence numbers from districts 10, 11, and 15, respectively.
  \item iINAR(1): Fit 9 INAR(1) models, one for each combination of crime type and district.
  \item iINAR(2): Fit 9 INAR(2) models, one for each combination of crime type and district.
  \item MGINAR(1): Fit the MGINAR(1) model to $\mathrm{vec}(\mathbf{X}_{t})$.
  \item MINAR(1).ICLSE: Fit the MINAR(1) model (without rank constraint) to $\mathbf{X}_{t}$ using ICLSE.
  \item RRMINAR(1).ICLSE: Fit the RRMINAR(1) model to $\mathbf{X}_{t}$ using ICLSE.
\end{enumerate}
Note that, except for MINAR(1) and RRMINAR(1), the parameter estimates of other models may contain a small number of slightly negative elements. This makes the relevant models unable to produce predictions and renders them practically meaningless. However, the small magnitude of these negative values implies that their contribution to the model's predictions is marginal. The negative values can be corrected by taking absolute values or setting them to zero to obtain fully non-negative parameter estimates. Here, we take their absolute values.
\begin{table}[htbp]
  \centering
  \caption{Comparison of in-sample fitting effects of seven models}
  \label{tab:within_sample_fitting}
  \begin{tabular}{lccccc}
      \toprule
      & E1 & E2 & E3 & E4 & Number of parameters \\
      \midrule
      MGINAR\_row(1) & 3209.2783 & 3.1440 & 3.5182 &0.5810 & 36 \\
      MGINAR\_column(1) & 3217.9238 & 3.1436 & 3.4759 & 0.5817 & 36 \\
      iINAR(1) & 3234.4988 & 3.1634 & 3.4212 & 0.5794 & \textbf{18} \\
      iINAR(2) & 3216.2731 & 3.1530 & 3.4059 & 0.5765 & 27 \\
      MGINAR(1) & 3285.9156 & 3.1975 & 4.0421 & 0.6111 & 90 \\       
      MINAR(1).ICLSE & 3183.7184 & 3.1246 & \textbf{3.3582} & \textbf{0.5736} & 27 \\     
      RRMINAR(1).ICLSE & \textbf{3179.7788} & \textbf{3.1224} & 3.3841 & 0.5754 & 19 \\     
      \bottomrule
  \end{tabular}
\end{table}

The in-sample fitting results in Table \ref{tab:within_sample_fitting} show that the MGINAR(1) model performs poorly across all indicators. Theoretically, models with more parameters should have smaller training errors. However, the results in Table \ref{tab:within_sample_fitting} indicate that the training error of the MGINAR(1) model under a limited sample size is not the smallest. This may be due to the sample size being too small and the model being overfitted. The first three indicators of MGINAR\_column(1) are all lower than those of MGINAR\_row(1), meaning that interactions among different crime types are greater than those among different districts. iINAR(2) performs better across all indicators compared to iINAR(1), indicating that higher-order lag terms can capture more information when interactions among districts and crime types are missed. In terms of the number of parameters, both iINAR(1) and RRMINAR(1) have a similar, small number of parameters. However, the former processes the data sequence alone, losing the structure information of the data and failing to capture the interrelationships among the sequences, which lead to poor fitting performance. From a model structure perspective, both RRMINAR(1) and MINAR(1) directly handle matrix data. The former is close to and lower than MINAR(1) in terms of $\rm{E}_1$, $\rm{E}_2$, while also having an advantage in the number of parameters. Based on the analysis of both aspects, RRMINAR(1) has a better fitting effect.

Table \ref{tab:out_of_sample_fitting} presents a comparison of the out-of-sample fitting indicators of seven models. The RRMINAR(1) model has a relatively small number of parameters while its various indicators are superior to those of the MINAR(1) model. iINAR(2) captures more information and has a better prediction performance compared to iINAR(1). MGINAR\_row(1) and MGINAR\_column(1) only consider the influences between columns and rows, respectively, resulting in an incomplete capture of the data’s structural information. Consequently, their prediction performance is worse than that of RRMINAR(1) and MINAR(1), but better than iINAR(1). These results further confirm the results that the influence among rows (crime types) is greater than that among columns (districts) in Tables \ref{tab:within_sample_fitting} and \ref{tab:out_of_sample_fitting}. 
\begin{table}[htbp]
  \centering
  \caption{Comparison of the out-of-sample fitting effects of seven models}
  \label{tab:out_of_sample_fitting}
  \begin{tabular}{lccccc}
      \toprule
      & E1 & E2 & E3 & E4 & Number of parameters \\
      \midrule
      MGINAR\_row(1) & 477.3997 & 2.7132 & 4.4923 & 1.2868 & 36 \\
      MGINAR\_column(1) & 475.6551 & 2.7055 & 4.4903 & 1.2921 & 36 \\
      iINAR(1) & 521.3338 & 2.9553 & 4.7442 & 1.3969 & \textbf{18} \\
      iINAR(2) & 473.3887 & 2.7118 & 4.3672 & 1.2684 & 27 \\
      MGINAR(1) & 427.9688 & 2.4472 & 4.1385 & 1.1738 & 90 \\      
      MINAR.ICLSE(1) & 385.6819 & 2.204 & 3.7618 & 1.0313 & 27 \\  
      RRMINAR.ICLSE(1) & \textbf{382.5831} & \textbf{2.1844} & \textbf{3.7216} & \textbf{1.0281} & 19 \\    
      \bottomrule
  \end{tabular}
\end{table}

\begin{figure}[!htbp]
  \centering
  \includegraphics[width=1\textwidth]{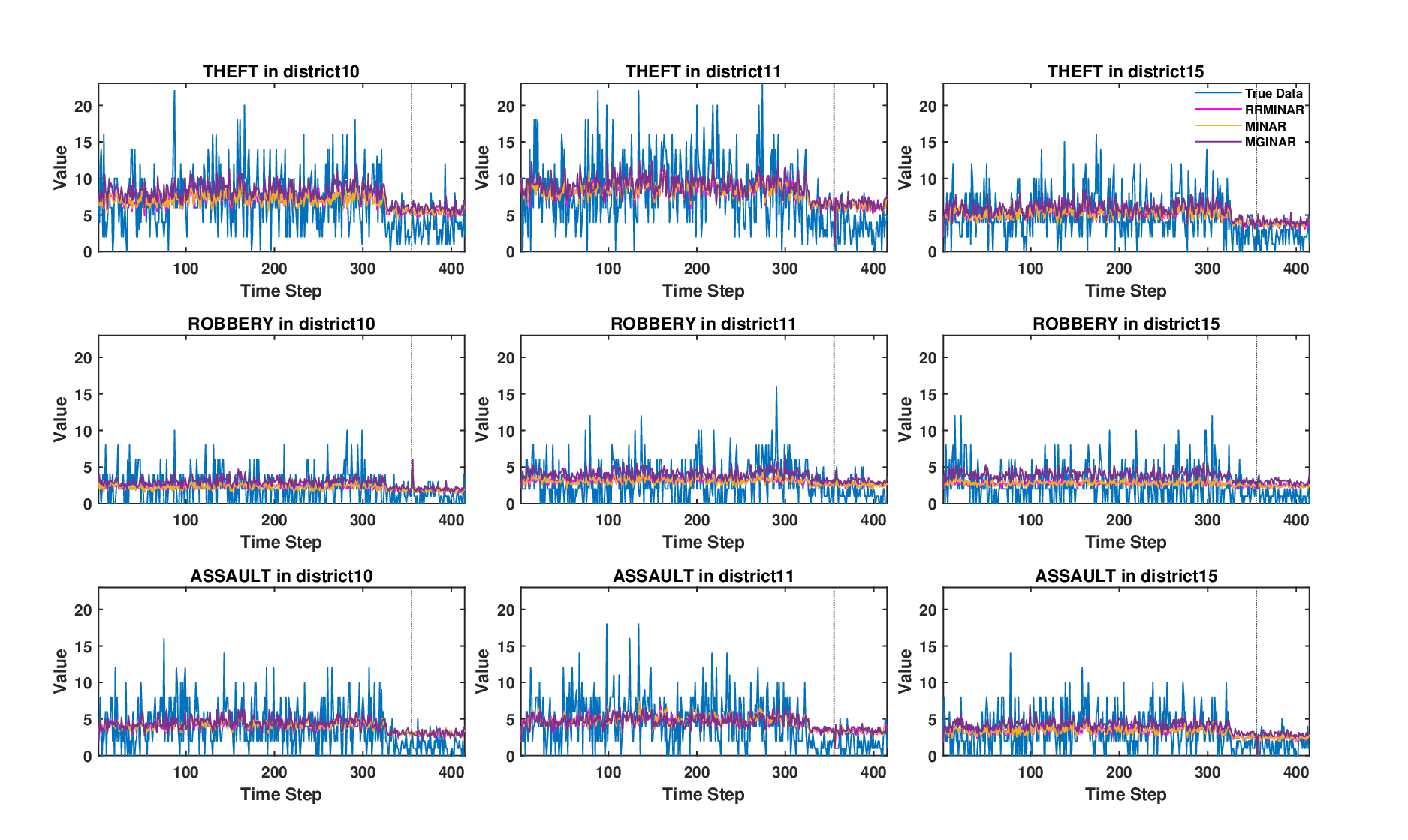}
  \caption{Comparison of original data and model predictions; out-of-sample predictions are to the right of the black dotted line.}
  \label{fig:Prediction comparison}
\end{figure}

Figure \ref{fig:Prediction comparison} shows the model fitting of each matrix sequence for the MGINAR(1), MINAR(1) and RRMINAR(1) models, as well as a comparison of their in-sample and out-of-sample prediction performance. The region to the right of the black dotted line represents the out-of-sample prediction results of the three models. It can be seen from Figure \ref{fig:Prediction comparison} that the prediction deviation of the MGINAR(1) model is relatively large. The prediction results of the RRMINAR(1) and MINAR(1) models are close but exhibit subtle differences.

Table \ref{tab:Estimation_of_Parameter_Matrices} presents the ICLSEs of $\mathbf{A}$, $\mathbf{B}$ and $\mathbf{C}$ of the RRMINAR(1) model, respectively. The estimater of the left-multiplicative coefficient matrix $\mathbf{A}$ indicates that THEFT has a stronger influence on ROBBERY and ASSAULT, while ROBBERY has a weaker influence on the other two crime types.
The estimater of right-multiplication coefficient matrix $\mathbf{B}$ indicates that the degree of influence among districts is relatively low. It can be seen from Table \ref{tab:Estimation_of_Parameter_Matrices} that district11 has a stronger influence on district10 and district15, while district10 is less affected by other districts.
These results further confirm the finding presented in Tables \ref{tab:within_sample_fitting} and \ref{tab:out_of_sample_fitting} that the influence among rows (crime types) is greater than that among columns (districts).
\begin{table}[htbp!]
  \centering
  \caption{Estimates of parameter matrices of the RRMINAR(1) model}
  \label{tab:Estimation_of_Parameter_Matrices} 
  \begin{tabular}{lcccc}
    \toprule
    Parameter matrix & \multicolumn{4}{c}{Estimates} \\
    \midrule
    \multirow{4}{*}{$\mathbf{A}$} & & THEFT & ROBBERY & ASSAULT \\
    & THEFT & 0.3719 & 0.3508 & 0.5973 \\
    & ROBBERY & 0.1005 & 0.0948 & 0.1614 \\
    & ASSAULT & 0.2501 & 0.2360 & 0.4018 \\
    \midrule
    \multirow{4}{*}{$\mathbf{B}$} & & District10 & District11 & District15 \\
    & District10 & 0.1065 & 0.1874 & 0.1861 \\
    & District11 & 0.1321 & 0.2324 & 0.2308 \\
    & District15 & 0.0883 & 0.1554 & 0.1544 \\
    \midrule
    \multirow{4}{*}{$\mathbf{C}$} & & District10 & District11 & District15 \\
    & THEFT & 4.2171 & 4.8405 & 2.7914 \\
    & ROBBERY & 1.5090 & 2.1651 & 2.1343 \\
    & ASSAULT & 2.1961 & 2.3775 & 1.7678 \\
    \bottomrule
  \end{tabular}
\end{table}

\section{Conclusion}
In this paper, the MINAR model is considered, which moves beyond the traditional vector model by directly processes matrix data. It accurately characterizes the cross-dependency of data row and column and reduces the number of model parameters for high-dimensional data. For example, for an \(m\times n\times T\) dimensional matrix time series, the traditional MGINAR(1) model involves \(m^2 n^2+mn\) parameters fitting, while MINAR(1) involves \(m^2+n^2+mn\) parameters. Howere, as the data dimension increases further, the redundant information in the data can lead to poor modeling performance of the MINAR(1) model and a sharp increase in the number of parameters. To address this issue, we propose the rank-reducing matrix integer-valued autoregressive (RRMINAR) model. The proposed model relies on an autoregressive structure based on a bilinear coefficient matrix and assumes low-rank constraints on the coefficient matrices. We develop an iterative conditional least squares estimation (ICLSE) and analyze its asymptotic properties. Compared with the MINAR model without a low-rank structure, the number of parameters involved in the RRMINAR(1) model is reduced to \(m^2+n^2-(m-k_1)^2-(n-k_2)^2+mn\) \citep{Xiao2022}, thereby improving estimation efficiency, reducing the risk of overfitting, and increasing robustness. As a result, the proposed model shows good performance in the real crime data analysis. For the proposed RRMINAR(1) model, this paper uses the $C_p$ criterion to select the appropriate rank for the coefficient matrices. Through simulation experiments, it can be observed that in low-rank scenarios, the proposed model outperforms both the MGINAR(1) and MINAR(1) models in most cases. Based on a case verification of criminal incidents in the Chicago area, the RRMINAR(1) model demonstrates better fitting and prediction performance compared to other models. However, there are still some issues worth studying in future. On the one hand, as discussed in Remark \ref{remark:3.1}, the model structure can be extended. From the crime data analysis, especially the results in Table \ref{tab:Estimation_of_Parameter_Matrices}, the proposed RRMINAR(1) model cannot fully explain the dynamic trend of the crime data, and a $p$-order RRMINAR model is required. On the other hand, the thinning operator imposes a constraint that all model parameters be positive, which limits the model's application. For instance, the model needs to be extended to scenarios such as data with negative correlations. Additionally, the rank selection method and its theoretical properties still require further investigation.

\section*{Funding}
This work was supported by the following funding sources:
\begin{itemize}
    \item National Natural Science Foundation of China (Grant No.~12201370);
    \item Fundamental Research Program of Shanxi Province (Grants No.~202203021211305 and 202203021222190);
    \item Wenying Young Scholar Talent Program of Shanxi University (Grant No.~109541073);
    \item Taiyuan University of Science and Technology Scientific Research Initial Funding (Grant No.~2022085);
    \item Reward Fund for Excellent Doctors Working in Shanxi Province (Grant No.~20242059).
\end{itemize}
\newpage
\bibliographystyle{elsarticle-num-names}
\bibliography{1}
\newpage

\appendix
\renewcommand{\thefigure}{A\arabic{figure}}
\setcounter{figure}{0}
\section{Appendix} 
\subsection{Figures}\label{fig:appendix}
\begin{figure}[htbp]
  \centering
  \includegraphics[width=1\textwidth]{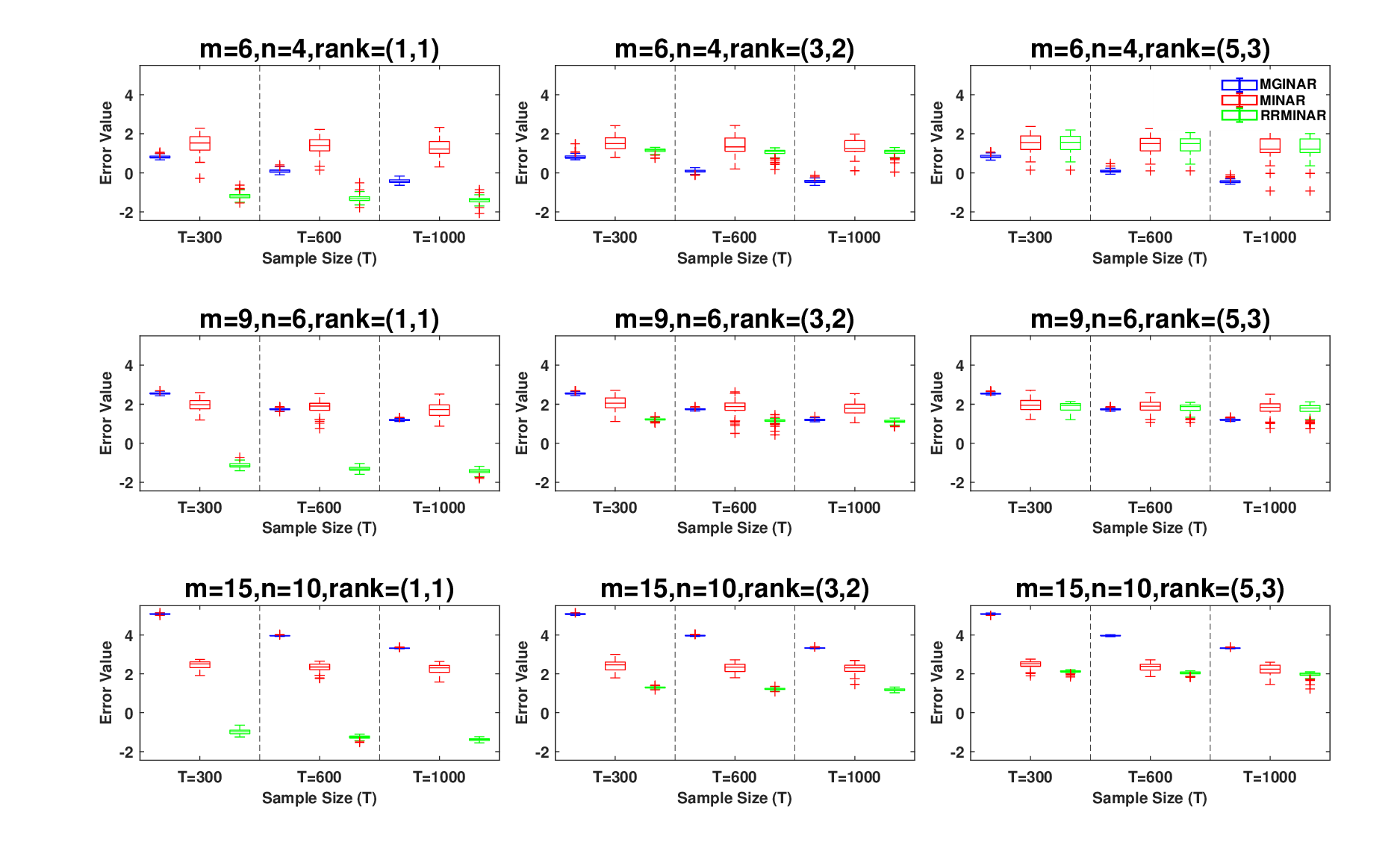}
  \caption{Comparison of the estimation errors for the MGINAR(1) (LSE), MINAR(1) (ICLSE), and RRMINAR(1) (ICLSE) models under Setting III. The three panels in each sub-figure correspond to sample sizes of 300, 600, and 1000, respectively.}
  \label{fig:set3_exp1}
\end{figure}

\begin{figure}[htbp]
  \centering
  \includegraphics[width=1\textwidth]{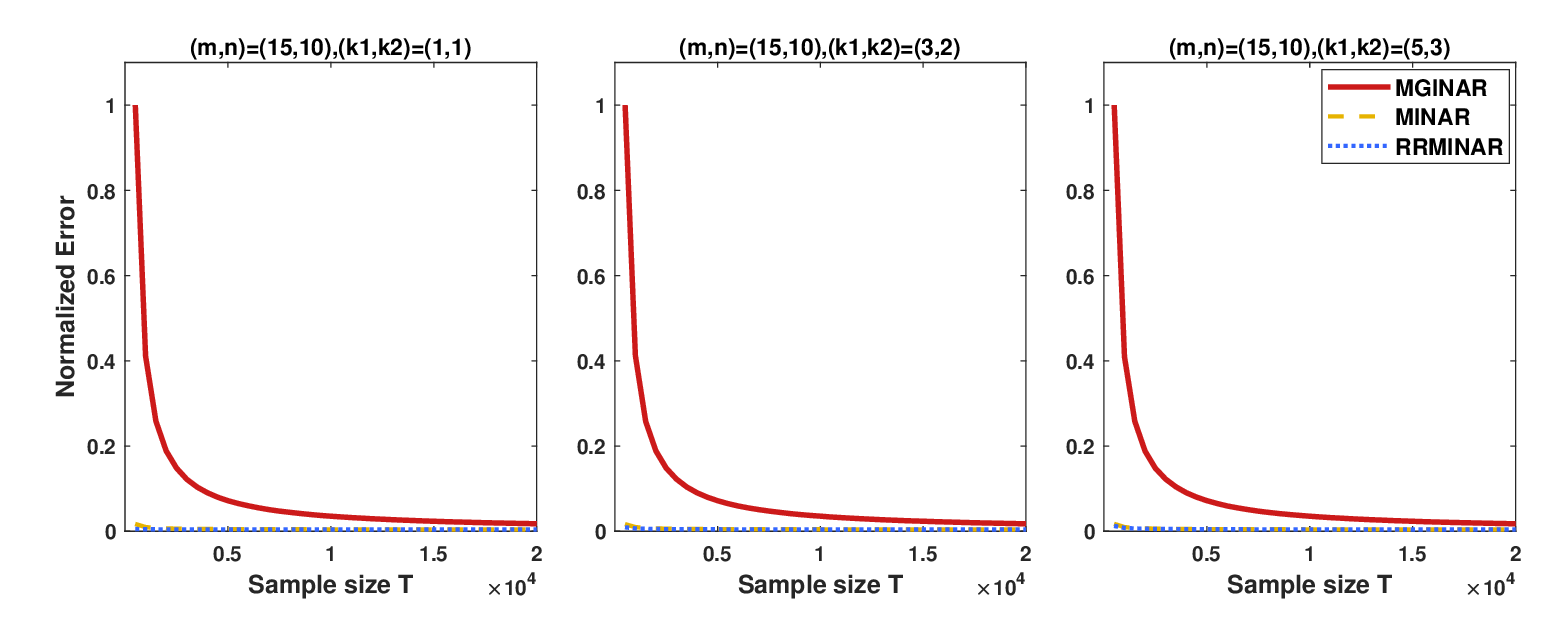}
  \caption{Convergence efficiency of coefficient matrix \(\mathbf{B}\otimes \mathbf{A}\) under Setting I.}
  \label{fig:se1_exp2}
  \end{figure}

\begin{figure}[htbp]
  \centering
  \includegraphics[width=1\textwidth]{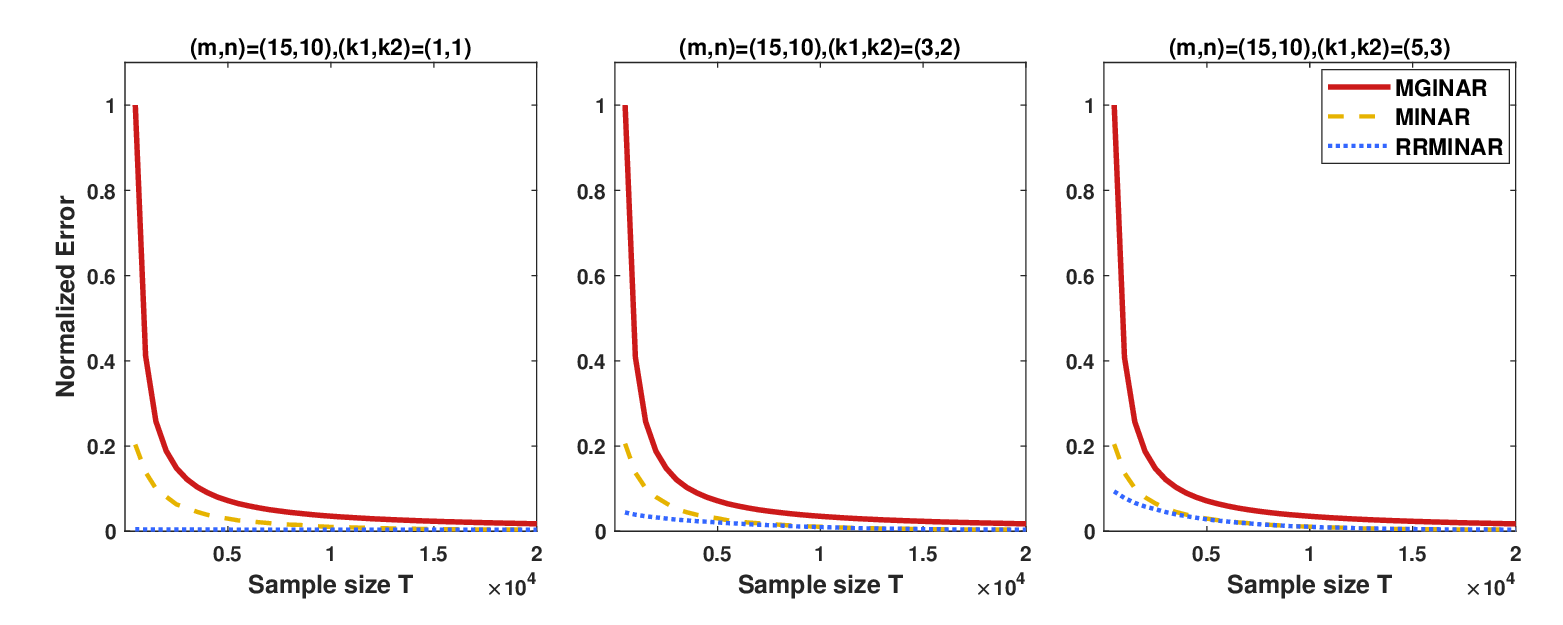}
  \caption{Convergence efficiency of coefficient matrix \(\mathbf{B}\otimes \mathbf{A}\) under Setting II.}
  \label{fig:set2_exp2}
  \end{figure}

Here, a comparison of the estimation errors for the three models—MGINAR(1), MINAR(1), and RRMINAR(1)—under Setting III is provided in Figure \ref{fig:set3_exp1}. The comparison covers different time series lengths \(T\), different dimensions \((m, n)\), and different ranks \((k_1, k_2)\). Additionally, convergence efficiency of the estimates of parameter matrices under different ranks with dimension \((15, 10)\) are presented for Settings I and II in Figures \ref{fig:se1_exp2} and \ref{fig:set2_exp2}, respectively.

\newpage

\subsection{Proof of Theorem \ref{pro:4.1}}\label{Proof of Proposition 4.1}
\begin{proof}
We first show that $\{\mathbf{\Delta}_t\}_{t\in\mathbb{Z}}$ is a matrix white noise sequence. From the property of the thinning operation and the matrix Kronecker product, we have \[\mathbb{E}\left(\mathbf{A}\circledast\mathbf{X}_{t}\circledast \mathbf{B}^{\top}\right)=\mathbf{A} \mathbb{E}(\mathbf{X}_{t})\mathbf{B}^{\top}.\]
Then
\[
\begin{aligned}
\mathbb{E} (\mathbf{\Delta}_{t}) &=\mathbb{E}\left(\mathbf{X}_{t}-\mathbf{A}\mathbf{X}_{t-1} \mathbf{B}^{\top} - \mathbf{C}\right) \\
&=\mathbb{E}\left(\mathbf{A}\circledast\mathbf{X}_{t-1}\circledast \mathbf{B}^{\top} + \mathbf{E}_t-\mathbf{A}\mathbf{X}_{t-1} \mathbf{B}^{\top} - \mathbf{C}\right) \\
&=\mathbf{A}\mathbb{E}(\mathbf{X}_{t-1})\mathbf{B}^{\top} + \mathbf{C}-\mathbf{A}\mathbb{E}(\mathbf{X}_{t-1})\mathbf{B}^{\top} - \mathbf{C} \\
&=\mathbf{0}_{m\times n}.
\end{aligned}
\]
Notice that, $\forall t_2<t_1$,
\[
\begin{aligned}
\mathbb{E}\left(\mathbf{\Delta}_{t_1}\otimes\mathbf{X}_{t_2}\right) &=\mathbb{E}\left\{(\mathbf{X}_{t_1}-\mathbf{A}\mathbf{X}_{t_1-1}\mathbf{B}^{\top} - \mathbf{C} )\otimes \mathbf{X}_{t_2}\right\} \\
&=\mathbb{E}\left[\left\{\mathbf{X}_{t_1}-\mathbb{E}\left(\mathbf{X}_{t_1} \mid \sigma\left(\mathbf{X}_{t_1-1}, \mathbf{X}_{t_1-2}, \ldots\right)\right)\right\}\otimes\mathbf{X}_{t_2}\right] \\
&=\mathbb{E}\left(\mathbf{X}_{t_1}\otimes\mathbf{X}_{t_2}\right)-\mathbb{E}\left\{\mathbb{E}\left(\mathbf{X}_{t_1}\otimes \mathbf{X}_{t_2} \mid \sigma\left(\mathbf{X}_{t_1-1}, \mathbf{X}_{t_1-2}, \ldots\right)\right)\right\} \\
&=\mathbf{0}_{m^2\times n^2}.
\end{aligned}
\]
Then $\mathbf{\Delta}_{t_1}$ is uncorrelated with the previous value $\mathbf{X}_{t_2}$. Thus we have, without loss of generality, for $t_2<t_1$,
\[
\begin{aligned}
\mathbb{E}\left(\mathbf{\Delta}_{t_1}\otimes\mathbf{\Delta}_{t_2}\right) &=\mathbb{E}\left\{\mathbf{\Delta}_{t_1}\otimes\left(\mathbf{X}_{t_2}-\mathbf{A}\mathbf{X}_{t_2-1} \mathbf{B}^{\top} - \mathbf{C}\right)\right\} \\
&=\mathbb{E}\left(\mathbf{\Delta}_{t_1}\otimes\mathbf{X}_{t_2}\right)- (\mathbf{I}_{m\times m}\otimes \mathbf{A}) \mathbb{E}\left\{\left(\mathbf{\Delta}_{t_1}\otimes\mathbf{X}_{t_2-1}\right)(\mathbf{I}_{n\times n}\otimes \mathbf{B}^{\top})\right\}-\mathbb{E}\left(\mathbf{\Delta}_{t_1}\otimes \mathbf{C}\right) \\
&=-\mathbb{E}\left(\mathbf{\Delta}_{t_1}\right)\otimes \mathbf{C} \\
&=\mathbf{0}_{m^2\times n^2}.
\end{aligned}
\]
Now we have shown that $\{\mathbf{\Delta}_t\}_{t\in\mathbb{Z}}$ is a matrix white noise sequence. Next, we derive its marginal covariance matrix. Since $\mathbb{E} (\mathbf{\Delta}_{t})=\mathbf{0}_{m\times n}$ and $\mathbb{E}\left(\mathbf{\Delta}_{t}\otimes\mathbf{X}_{t-1}\right)=\mathbf{0}_{m^2\times n^2}$, we have
\[
\begin{aligned}
&~\quad\mathbb{E}\left\{{\rm vec}(\mathbf{\Delta}_{t}){\rm vec}(\mathbf{\Delta}_{t})^{\top}\right\} \\
&=\mathbb{E}\left[{\rm vec}(\mathbf{\Delta}_{t})\left\{{\rm vec}(\mathbf{X}_{t})-{\rm vec}(\mathbf{A}\mathbf{X}_{t-1} \mathbf{B}^{\top}) - {\rm vec}(\mathbf{C})\right\}^{\top}\right] \\
&=\mathbb{E}\left\{{\rm vec}(\mathbf{\Delta}_{t}) {\rm vec}(\mathbf{X}_{t})^{\top}\right\} \\
&=\mathbb{E}\left[\left\{{\rm vec}(\mathbf{X}_{t})-{\rm vec}(\mathbf{A}\mathbf{X}_{t-1} \mathbf{B}^{\top}) - {\rm vec}(\mathbf{C})\right\}\left\{{\rm vec}(\mathbf{A}\circledast\mathbf{X}_{t-1}\circledast \mathbf{B}^{\top}) + {\rm vec}(\mathbf{E}_t)\right\}^{\top}\right] \\
&=\mathbb{E}\left[\left\{{\rm vec}(A\circledast\mathbf{X}_{t-1}\circledast \mathbf{B}^{\top}) + {\rm vec}(\mathbf{E}_t)-{\rm vec}(\mathbf{A}\mathbf{X}_{t-1} \mathbf{B}^{\top}) - {\rm vec}(\mathbf{C})\right\}\left\{{\rm vec}(\mathbf{A}\circledast\mathbf{X}_{t-1}\circledast \mathbf{B}^{\top}) + {\rm vec}(\mathbf{E}_t)\right\}^{\top}\right] \\
&={\rm Cov}\{{\rm vec}(\mathbf{E}_t)\}+\mathbb{E}\left\{{\rm vec}(\mathbf{A}\circledast\mathbf{X}_{t-1}\circledast \mathbf{B}^{\top}){\rm vec}(\mathbf{A}\circledast\mathbf{X}_{t-1}\circledast \mathbf{B}^{\top})^{\top}-{\rm vec}(\mathbf{A}\mathbf{X}_{t-1} \mathbf{B}^{\top}){\rm vec}(\mathbf{A}\circledast\mathbf{X}_{t-1}\circledast \mathbf{B}^{\top})^{\top}\right\}.
\end{aligned}
\]
For fixed $\mathbf{Y}= (y_{i,j})_{1\leq i\leq m,1\leq j\leq n}\in \mathbb{N}_{0}^{m\times n}$, we have
\[\mathbf{A}\circledast\mathbf{Y}\circledast \mathbf{B}^{\top}=\left\{\sum_{l=1}^{n}\sum_{k=1}^{m}(\beta_{j,l}\alpha_{i,k})\circ y_{k,l}\right\}_{1\leq i\leq m,1\leq j\leq n}.\]

According to Definition \ref{def:2.1} and the additivity of the Poisson distribution, the component $\sum_{l=1}^{n}\sum_{k=1}^{m}(\beta_{j,l}\alpha_{i,k})\circ y_{k,l}$ is Poisson random variable with parameter $\sum_{l=1}^{n}\sum_{k=1}^{m}(\beta_{j,l}\alpha_{i,k}) y_{k,l}$ for any $ 1\leq i \leq m,1 \leq j\leq n$. Since the thinning operator $\circ$ are mutually independent, the above components are uncorrelated. Notice that the expectation of the Poisson distribution is equal to the variance. Then the covariance matrix of ${\rm vec}(\mathbf{A}\circledast\mathbf{Y}\circledast \mathbf{B}^{\top})$ is a diagonal matrix, and the diagonal element is $\sum_{l=1}^{n}\sum_{k=1}^{m}(\beta_{j,l}\alpha_{i,k}) y_{k,l}$. Therefore, for any fixed $\mathbf{Y}$, we have 
 \[{\rm Cov}\left\{{\rm vec}(\mathbf{A}\circledast\mathbf{Y}\circledast \mathbf{B}^{\top})\right\}=\operatorname{diag}\left\{{\rm vec}(\mathbf{A}\mathbf{Y} \mathbf{B}^{\top})\right\}.\]
Notice that
\[
\begin{aligned}
&~\quad\mathbb{E}\left\{{\rm vec}(\mathbf{A}\mathbf{X}_{t-1} \mathbf{B}^{\top}){\rm vec}(\mathbf{A}\circledast\mathbf{X}_{t-1}\circledast \mathbf{B}^{\top})^{\top} \mid \sigma\left(\mathbf{X}_{t-1}\right)\right\} \\
&=\rm{vec}(\mathbf{A}\mathbf{X}_{t-1} \mathbf{B}^{\top})\mathbb{E}\left\{\rm{vec}(\mathbf{A}\circledast\mathbf{X}_{t-1}\circledast \mathbf{B}^{\top})^{\top} \mid \sigma\left(\mathbf{X}_{t-1}\right)\right\}\\
&=\mathbb{E}\left\{\rm{vec}(\mathbf{A}\mathbf{X}_{t} \mathbf{B}^{\top}) \mid \sigma\left(\mathbf{X}_{t-1}\right)\right\}\mathbb{E}\left\{\rm{vec}(\mathbf{A}\circledast\mathbf{X}_{t-1}\circledast \mathbf{B}^{\top})^{\top}\mid \sigma\left(\mathbf{X}_{t-1}\right)\right\}.\\
\end{aligned}
\]
So we can get
\[
\begin{aligned}
&~\quad\mathbb{E}\left\{{\rm vec}(\mathbf{A}\circledast\mathbf{X}_{t-1}\circledast \mathbf{B}^{\top}){\rm vec}(\mathbf{A}\circledast\mathbf{X}_{t-1}\circledast \mathbf{B}^{\top})^{\top}-{\rm vec}(\mathbf{A}\mathbf{X}_{t-1} \mathbf{B}^{\top}){\rm vec}(\mathbf{A}\circledast\mathbf{X}_{t-1}\circledast \mathbf{B}^{\top})^{\top}\right\}\\
&=\mathbb{E}\left[{\rm Cov}\left\{{\rm vec}(\mathbf{A}\circledast\mathbf{X}_{t-1}\circledast \mathbf{B}^{\top})\right\}\mid \sigma\left(\mathbf{X}_{t-1}\right)\right]\\
&=\mathbb{E}\left[\operatorname{diag}\left\{{\rm vec}(\mathbf{A}\mathbf{X}_{t-1} \mathbf{B}^{\top})\right\}\right].
\end{aligned}
\]
On the other hand, since the components of $\mathbf{E}_{t}$ are Poisson random variables and mutually independent, we have ${\rm Cov}\left\{{\rm vec}(\mathbf{E}_{t})\right\}=\operatorname{diag}\{{\rm vec}(\mathbf{C})\}$. Recall $\mathbb{E}(\mathbf{\Delta}_t) = \mathbf{0}_{m\times n}$, and let $\mathbf{1}_{a\times b} $ represent the $a\times b$ matrix of all ones. We have 
\[
\begin{aligned}
\mathbf{0}_{mn}&=\mathbb{E}\left\{\operatorname{vec}\left(\mathbf{\Delta}_t\right)\right\} \\
&=\mathbb{E}\left\{\operatorname{vec}\left(\mathbf{X}_t\right)-(\mathbf{B} \otimes \mathbf{A})  \operatorname{vec}\left(\mathbf{X}_{t-1}\right)\right\}-\operatorname{vec}(\mathbf{C}) \\
&=(\mathbf{1}_{m m \times m n}-\mathbf{B} \otimes \mathbf{A})\mathbb{E}\left\{\operatorname{vec}\left(\mathbf{X}_{t-1}\right)\right\}-\operatorname{vec}(\mathbf{C}),\\
\end{aligned}
\]
which implies $\mathbb{E}\left\{\operatorname{vec}\left(\mathbf{X}_{t-1}\right)\right\}=(\mathbf{1}_{m m \times m n}-\mathbf{B} \otimes \mathbf{A})^{-1} \operatorname{vec}(\mathbf{C})$ and
\[
\begin{aligned}
\operatorname{vec}\left\{\mathbf{A}\mathbb{E}\left(\mathbf{X}_{t-1}\right) \mathbf{B} ^{\top}\right\} = (\mathbf{B} \otimes \mathbf{A}) \operatorname{vec}\left\{\mathbb{E}\left(\mathbf{X}_{t-1}\right)\right\} = (\mathbf{B} \otimes \mathbf{A})\left(\mathbf{1}_{m m \times m n}- \mathbf{B} \otimes \mathbf{A}\right)^{-1} \operatorname{vec}(\mathbf{C}).
\end{aligned}
\]
Therefore, it holds that
\[
\begin{aligned}
&~\quad\mathbb{E}\left\{{\rm vec}(\mathbf{\Delta}_{t}){\rm vec}(\mathbf{\Delta}_{t})^{\top}\right\}\\ 
&=\operatorname{diag}\left\{{\rm vec}(\mathbf{C})\right\}+\mathbb{E}\left[\operatorname{diag}\left\{{\rm vec}(\mathbf{A}\mathbf{X}_{t-1} \mathbf{B}^{\top})\right\}\right] \\
&=\operatorname{diag}\left[{\rm vec}(\mathbf{C})+{\rm vec}\{\mathbf{A}\mathbb{E}(\mathbf{X}_{t-1}) \mathbf{B}^{\top}\}\right] \\
&=\operatorname{diag}\left\{{\rm vec}(\mathbf{C})+(\mathbf{B}\otimes \mathbf{A})\left(\mathbf{1}_{mn \times mn}- \mathbf{B}\otimes \mathbf{A}\right)^{-1}{\rm vec}(\mathbf{C})\right\} \\
&=\operatorname{diag}\left\{\left(\mathbf{1}_{mn \times mn}- \mathbf{B}\otimes \mathbf{A}\right)^{-1}{\rm vec}(\mathbf{C})\right\}.
\end{aligned}
\]
This completes the proof of Theorem 
\ref{pro:4.1}\label{Proof of Proposition 4.1}.
\end{proof}

\subsection{Proof of Theorem \ref{t:Iterative convergence}}\label{proof:iterative consistency}
\begin{proof}
To prove Theorem \ref{t:Iterative convergence}, we need the Lemmas \ref{t:4.1}-\ref{lemmaA.5}. The proofs of Lemmas \ref{t:4.1}-\ref{lemmaA.5} are given in Section \ref{proof:t:4.1}-\ref{proof of lemmaA.5}, respectively.

\begin{lemma}\label{t:4.1} 
Let $\widehat{\mathbf{U}}^{\top}_T\widehat{\mathbf{\Phi}}_T\widehat{\mathbf{V}}_T=\widehat{\mathbf{\Lambda}}_T$ and $\mathbf{U}^{\top}\mathbf{\Phi} \mathbf{V}=\mathbf{\Lambda}$ be the SVDs of $\widehat{\mathbf{\Phi}}_T\in\mathbb{R}^{m\times n}$ and $\mathbf{\Phi}\in\mathbb{R}^{m\times n}$, respectively. Assume that $\widehat{\mathbf{\Phi}}_T = \mathbf{\Phi} + \mathit{o}_\mathrm{p}(T^\mathit{-a})$ with $a>0$. Then, it holds that
\begin{itemize}
  \item[\rm(i)] $\widehat{\mathbf{\Lambda}}^2_T(i,i) = \mathbf{\Lambda}^2(i,i) + \mathit{o}_\mathrm{p}(T^\mathit{-a}),~1\leq i\leq n;$
  \item[\rm(ii)] $|\widehat{\mathbf{U}}_T\mathit{(i,j)}| = |\mathbf{U}\mathit{(i,j)}| + \mathit{o}_\mathrm{p}(T^\mathit{-a}),~1\leq i,j\leq m;$
  \item[\rm(iii)] $|\widehat{\mathbf{V}}_T\mathit{(i,j)}| = |\mathbf{V}\mathit{(i,j)}| + \mathit{o}_\mathrm{p}(T^\mathit{-a}),~1\leq i,j\leq n.$
\end{itemize}
\end{lemma}

\begin{lemma}\label{lemmaA.4}
Assume that $\hat{\mathbf{A}} =  \mathbf{A} + O_{\mathrm{p}}(a_{n})$ with invertible matrix $\mathbf{A}\in \mathbb{R}^{m\times m}$ and $a_{n}=o(1)$. Then, it holds that $\hat{\mathbf{A}}^{-1}=\mathbf{A}^{-1}+O_{\mathrm{p}}(a_{n})$.
\end{lemma}

\begin{lemma}\label{lemmaA.5}
Let \(\{\mathbf{X}_t\}_{t = 1}^{T}\) be a \(m \times n\) MINAR(1) sequence with coefficient matrices $\mathbf{A} \in \mathbb{R}^{m\times m}$, $\mathbf{B} \in \mathbb{R}^{n\times n}$ and $\mathbf{C} \in \mathbb{R}^{m\times n}$, and define \(\mathbf{\mathbf{\Delta}}_t\coloneqq\mathbf{X}_t - \mathbf{A}\mathbf{X}_{t-1}\mathbf{B}^{\top} - \mathbf{C},~t\in\mathbb{Z}\). If \(\hat{\mathbf{A}} - \mathbf{A} = O_{\mathrm{p}}(T^{-1/2})\), then we  have
\begin{itemize}
  \item[\rm(i)] \(\sum_{t = 1}^{T}(\hat{\mathbf{A}} - \mathbf{A})\mathbf{X}_t = O_{\mathrm{p}}(\sqrt{T});\)
  \item[\rm(ii)] \(\sum_{t = 1}^{T}(\hat{\mathbf{A}} - \mathbf{A})^\top(\hat{\mathbf{A}} - \mathbf{A})\mathbf{X}_t = O_{\mathrm{p}}\left(\frac{1}{\sqrt{T}}\right);\)
  \item[\rm(iii)] \(\sum_{t = 2}^{T}\mathbf{\Delta}_t^{\top}(\hat{\mathbf{A}} - \mathbf{A})^{\top}\mathbf{X}_{t-1} = o_{\mathrm{p}}(\sqrt{T}).\)
\end{itemize}
Similar results also hold for $\hat{\mathbf{B}} - \mathbf{B}$.
\end{lemma}

Define $\hat{\mathbf{M}}_1=T^{-1/2}\hat{\mathbf{S}}_{1yx} \hat{\mathbf{S}}_{1xx}^{-\frac{1}{2}}$ with $\mathbf{\hat{S}}_{1xx}=\sum_{t = 2}^{T}\mathbf{X}_{t - 1}\mathbf{\hat{B}}_{0}^{\top}\mathbf{\hat{B}}_{0}\mathbf{X}_{t - 1}^{\top}$ and 
$\mathbf{\hat{S}}_{1yx}=\sum_{t = 2}^{T}(\mathbf{X}_{t}-\mathbf{\hat{C}}_{0})\mathbf{\hat{B}}_{0}\mathbf{X}_{t - 1}^{\top}$. Recall \(\mathbf{\mathbf{\Delta}}_t=\mathbf{X}_t - \mathbf{A}\mathbf{X}_{t-1}\mathbf{B}^{\top} - \mathbf{C}\), \(\Gamma_2 = \mathbb{E}(\mathbf{X}_t\mathbf{B}^{\top}\mathbf{B}\mathbf{X}_{t}^{\top})\), \(\mathbf{\hat{B}}_{0}=\mathbf{B} + \mathit{O}_\mathrm{p}(T^{-1/2})\) and \(\mathbf{\hat{C}}_{0}=\mathbf{C} + \mathit{O}_\mathrm{p}(T^{-1/2})\). Notice that \(T^{-1} \sum_{t = 2}^{T} \Delta_{t} \mathbf{B}\mathbf{X}_{t - 1}^{\top}=O_{\rm p}(T^{-1/2})\) and \(\hat{\mathbf{S}}_{1xx}=\Gamma_2+O_{\rm p}(T^{-1/2})\). Letting \(\mathbf{\tilde{M}}_1=T^{-\frac{1}{2}}(\sum_{t = 2}^{T}\mathbf{A}\mathbf{X}_{t - 1}\mathbf{B}^{\top}+\mathbf{C} - \mathbf{\hat{C}}_{0})\mathbf{\hat{B}}_{0}\mathbf{X}_{t - 1}^{\top}\mathbf{\hat{S}}_{1xx}^{-\frac{1}{2}}\), by Lemmas \ref{lemmaA.4} and \ref{lemmaA.5}, we have
\begin{align}\label{equ:A.1}
&\mathbf{\hat{M}}_1 - \mathbf{\tilde{M}}_1 = \frac{1}{\sqrt{T}}\left\{\mathbf{\hat{S}}_{1yx}-\sum_{t = 2}^{T}(\mathbf{A}\mathbf{X}_{t - 1}\mathbf{B}^{\top}+\mathbf{C} - \mathbf{\hat{C}}_{0})\mathbf{\hat{B}}_{0}\mathbf{X}_{t - 1}\right\}\mathbf{\hat{S}}_{1xx}^{-\frac{1}{2}}\notag \\
&= \frac{1}{T}\left\{\sum_{t = 2}^{T}(\mathbf{A}\mathbf{X}_{t - 1}\mathbf{B}^{\top}+\mathbf{C} + \mathbf{\Delta}_t - \mathbf{\hat{C}}_{0})\mathbf{\hat{B}}_{0}\mathbf{X}_{t - 1} -\sum_{t = 2}^{T}(\mathbf{A}\mathbf{X}_{t - 1}\mathbf{B}^{\top}+\mathbf{C} - \mathbf{\hat{C}}_{0})\mathbf{\hat{B}}_{0}\mathbf{X}_{t - 1}\right\}\left(\frac{1}{T}\mathbf{\hat{S}}_{1xx}\right)^{-\frac{1}{2}} \notag\\
&= \left\{\frac{1}{T}\sum_{t = 2}^{T}(\mathbf{\Delta}_t\mathbf{B}\mathbf{X}_{t - 1}^{\top})+\mathit{o}_\mathrm{p}\left(\frac{1}{\sqrt{T}}\right)\right\}\left\{\Gamma_2^{-\frac{1}{2}}+\mathit{O}_\mathrm{p}\left(\frac{1}{\sqrt{T}}\right)\right\}\notag \\
&= \frac{1}{T}\sum_{t = 2}^{T}(\mathbf{\Delta}_t\mathbf{B}\mathbf{X}_{t - 1}^{\top})\Gamma_2^{-\frac{1}{2}}+\mathit{o}_\mathrm{p}\left(\frac{1}{\sqrt{T}}\right)
=\mathit{O}_\mathrm{p}\left(\frac{1}{\sqrt{T}}\right).
\end{align}
Writing \(\mathbf{M}_1=\mathbf{A}\Gamma_2^{1/2}\), then we have 
\begin{align}\label{new_equ:A.2}
  \mathbf{\tilde{M}}_{1}
  &=\frac{1}{\sqrt{T}}\sum_{t = 2}^{T}(\mathbf{A}\mathbf{X}_{t-1}\mathbf{B}^{\top}+\mathbf{C}-\mathbf{\hat{C}}_{0})\mathbf{\hat{B}}_{0}\mathbf{X}_{t-1}\mathbf{\hat{S}}_{1xx}^{-\frac{1}{2}}\notag\\
  &=\left\{\frac{1}{T}\sum_{t = 2}^{T}(\mathbf{A}\mathbf{X}_{t-1}\mathbf{B}^{\top}\mathbf{\hat{B}}_{0}\mathbf{X}_{t-1}^{\top})+\frac{1}{T}\sum_{t = 2}^{T}(\mathbf{C}-\mathbf{\hat{C}}_{0})\mathbf{\hat{B}}_{0}\mathbf{X}_{t-1}^{\top}\right\}\left(\frac{1}{T}\mathbf{\hat{S}}_{1xx}\right)^{-\frac{1}{2}}\notag\\
  &=\left\{\mathbf{A}\Gamma_2+\mathit{O}_\mathrm{p}\left(\frac{1}{\sqrt{T}}\right)\right\}\left\{\Gamma_2^{-\frac{1}{2}}+\mathit{O}_\mathrm{p}\left(\frac{1}{\sqrt{T}}\right)\right\}\notag\\
  &=\mathbf{M}_1+\mathit{O}_\mathrm{p}\left(\frac{1}{\sqrt{T}}\right).
\end{align}
Let \(\mathbf{\hat{M}}_1=\mathbf{\hat{U}_1}\mathbf{\hat{D}}_1\mathbf{\hat{V}}_1^{\top}\), \(\mathbf{\tilde{M}}_1=\mathbf{\tilde{U}_1}\mathbf{\tilde{D}}_1\mathbf{\tilde{V}}_1^{\top}\) and \(\mathbf{M}_1 = \mathbf{U}_1\mathbf{D}_1\mathbf{V}_1^{\top}\) be the SVD of \(\mathbf{\hat{M}}_1\), \(\mathbf{\tilde{M}}_1\) and \(\mathbf{M}_1\), respectively. By \eqref{equ:A.1}, \eqref{new_equ:A.2} and Lemma \ref{t:4.1}, it holds that 
\begin{align}
&\mathbf{\hat{U}_1}=\mathbf{\tilde{U}_1} + \mathit{O}_\mathrm{p}({T}^{-1/2}),\quad 
\mathbf{\hat{V}}_1=\mathbf{\tilde{V}}_1 + \mathit{O}_\mathrm{p}({T}^{-1/2}),\quad \mathbf{\hat{D}}_1=\mathbf{\tilde{D}}_1 + \mathit{O}_\mathrm{p}({T}^{-1/2});\label{new_equ:A.3}\\
&\mathbf{{U}}_1=\mathbf{\tilde{U}}_1+O_{\mathrm{p}}({T}^{-1/2}),\quad
{\mathbf{V}}_{1}=\tilde{\boldsymbol{V}}_{1}+O_{\mathrm{p}}({T}^{-1/2}),\quad
{\mathbf{D}}_{1}=\tilde{\mathbf{D}}_{1}+O_{\mathrm{p}}({T}^{-1/2}).\label{new_equ:A.4}
\end{align}
Let \(\mathbf{\hat{U}}_{1:k_1}\), \(\mathbf{\tilde{U}}_{1:k_1}\) and \(\mathbf{U}_{1:k_1}\) be the matrices composed of the first \(k_1\) normalized left-singular vectors of \(\mathbf{\hat{M}}_1\), \(\mathbf{\tilde{M}}_1\) and \(\mathbf{M}_1\), respectively; let \(\mathbf{\hat{V}}_{1:k_1}\), \(\mathbf{\tilde{V}}_{1:k_1}\) and \(\mathbf{V}_{1:k_1}\) be the matrices composed of the first \(k_1\) normalized right-singular vectors of \(\mathbf{\hat{M}}_1\), \(\mathbf{\tilde{M}}_1\) and \(\mathbf{M}_1\), respectively; and let \(\mathbf{\hat{D}}_{1:k_1}\), \(\mathbf{\tilde{D}}_{1:k_1}\) and \(\mathbf{D}_{1:k_1}\) be the matrices composed of the \(k_1\)th principal sub-matrix of \(\mathbf{\hat{M}}_1\), \(\mathbf{\tilde{M}}_1\) and \(\mathbf{M}_1\), respectively. Since \({\rm rank}(\mathbf{M}_1)=k_1\), we have \(\mathbf{M}_1 = \mathbf{A}\Gamma_2^{1/2}=\mathbf{U}_{1:k_1}\mathbf{D}_{1:k_1}\mathbf{V}_{1:k_1}^{\top}\). Using the orthogonality of the left-singular matrix, we have \(\mathbf{U}_{1:k_1}^{\top}\mathbf{U}_{1:k_1}= \mathbf{I}_{k_1}\). Therefore, it holds that
\begin{align*}
(\mathbf{I}_{m}-\mathbf{U}_{1:k_1}\mathbf{U}_{1:k_1}^{\top})\mathbf{A}\Gamma_2^{1/2}
&=(\mathbf{I}_{m}-\mathbf{U}_{1:k_1}\mathbf{U}_{1:k_1}^{\top})\mathbf{U}_{1:k_1}\mathbf{D}_{1:k_1}\mathbf{V}_{1:k_1}^{\top}\\
&=(\mathbf{U}_{1:k_1}-\mathbf{U}_{1:k_1})\mathbf{D}_{1:k_1}\mathbf{V}_{1:k_1}^{\top}=\mathbf{0}.
\end{align*}
which implies
\begin{align}\label{new_equ:A.5}
\mathbf{U}_{1:k_1}\mathbf{U}_{1:k_1}^{\top}\mathbf{A}=\mathbf{A}.
\end{align}
By \eqref{equ:A.1}-\eqref{new_equ:A.5}, we have
\begin{align*}
\hat{\mathbf{A}}_1&=\hat{\mathbf{U}}_{1:k_1} \hat{\mathbf{U}}_{1:k_1}^{\top} \hat{\mathbf{S}}_{1yx} \hat{\mathbf{S}}_{1xx}^{-1}
=\hat{\mathbf{U}}_{1:k_1} \hat{\mathbf{U}}_{1:k_1}^{\top} (\sqrt{T}\hat{\mathbf{M}}_1 \hat{\mathbf{S}}_{1xx}^{-\frac{1}{2}})\\
&=\bigg\{\mathbf{U}_{1:k_{1}}+O_{\rm p}\bigg(\frac{1}{\sqrt{T}}\bigg)\bigg\} \bigg\{\mathbf{U}_{1:k_{1}}+O_{\rm p}\bigg(\frac{1}{\sqrt{T}}\bigg)\bigg\}^{\top} \bigg\{\mathbf{M}_1 + O_{\mathrm{p}}\bigg(\frac{1}{\sqrt{T}}\bigg)\bigg\}\bigg(\frac{1}{T}\hat{\mathbf{S}}_{1xx}\bigg)^{-\frac{1}{2}}\\
&=\bigg\{\mathbf{U}_{1:k_{1}}\mathbf{U}_{1:k_{1}}^{\top}\mathbf{M}_1+O_{\rm p}\bigg(\frac{1}{\sqrt{T}}\bigg)\bigg\}\bigg\{\Gamma_2^{-1/2} + O_{\mathrm{p}}\bigg(\frac{1}{\sqrt{T}}\bigg)\bigg\}\\
&=\mathbf{U}_{1:k_{1}}\mathbf{U}_{1:k_{1}}^{\top}\mathbf{A} + O_{\mathrm{p}}\bigg(\frac{1}{\sqrt{T}}\bigg)
=\mathbf{A} + O_{\mathrm{p}}\bigg(\frac{1}{\sqrt{T}}\bigg).
\end{align*}
Similarly, we have \(\hat{\mathbf{B}}_1=\mathbf{B} + O_{\mathrm{p}}(T^{-1/2})\). Using \(\mathbf{A}\mathbf{X}_{t - 1}\mathbf{B}^{\top}+\mathbf{C} + \mathbf{\Delta}_t\) to replace \(\mathbf{X}_t\) in \eqref{equ:3.3}, we obtain 
\[\mathbf{\hat{C}}_1=\frac{1}{T}(\mathbf{A}\mathbf{X}_{t - 1}\mathbf{B}^{\top}+\mathbf{C} + \mathbf{\Delta}_t-\mathbf{\hat A}_1\mathbf{X}_{t - 1}\mathbf{\hat B}_1^{\top}),\]
which implies
\begin{align*}
  \hat{\mathbf{C}}_1&=\frac{T-1}{T}\mathbf{C}+\frac{1}{T}\sum_{t = 2}^{T}(\mathbf{\hat{A}}_1-\mathbf{A})\mathbf{X}_{t-1}\mathbf{B}^{\top}+\frac{1}{T}\sum_{t = 2}^{T}\mathbf{A}\mathbf{X}_{t-1}(\mathbf{\hat{B}}_1-\mathbf{B})^{\top} +\frac{1}{T}\sum_{t = 2}^{T}\mathbf{\Delta}_t+\mathit{o}_\mathrm{p}\left(\frac{1}{\sqrt{T}}\right)\\
  &=\mathbf{C}+\mathit{O}_\mathrm{p}\left(\frac{1}{\sqrt{T}}\right).
\end{align*}
According to the iterative process in Algorithm \ref{alg:rrminar1}, by mathematical induction, we have $\mathbf{\hat{A}}_{r}=\mathbf{A}+\mathit{O}_\mathrm{p}(T^{-1/2})$, $\mathbf{\hat{B}}_{r}=\mathbf{B}+\mathit{O}_\mathrm{p}(T^{-1/2})$ and $\mathbf{\hat{C}}_{r}=\mathbf{C}+\mathit{O}_\mathrm{p}(T^{-1/2})$ for any $r\ge 1$. This completes the proof of Theorem \ref{t:Iterative convergence}.
\end{proof}

\subsection{Proof of Theorem \ref{t:4.2}}\label{Proof of Theorem 4.2}
\begin{proof}
In the proposed Algorithm \ref{alg:rrminar1}, the estimates of $\mathbf{A}$, $\mathbf{B}$, and $\mathbf{C}$ are updated alternately in an iterative manner. Under the assumptions $\hat{\mathbf{B}}_{r+1}-\hat{\mathbf{B}}_{r}=o_\mathrm{p}(T^{-1/2})$ and $\hat{\mathbf{C}}_{r+1}-\hat{\mathbf{C}}_{r}=o_\mathrm{p}(T^{-1/2})$ for some $r>0$, Algorithm \ref{alg:rrminar1} terminates after a finite number of iterations. Without loss of generality, assume that the algorithm stops at the $r$th iteration. Then, $\hat{\mathbf{A}}_{\rm RR.LS}=\hat{\mathbf{A}}_r$ is calculated from $\hat{\mathbf{B}}_{r-1}$ and $\hat{\mathbf{C}}_{r-1}$ in step 9 of Algorithm \ref{alg:rrminar1}; $\hat{\mathbf{B}}_{\rm RR.LS}=\hat{\mathbf{B}}_{r}$ is calculated from $\hat{\mathbf{A}}_{r}$ and $\hat{\mathbf{C}}_{r-1}$ in step 16 of Algorithm \ref{alg:rrminar1}; $\hat{\mathbf{C}}_{\rm RR.LS}=\hat{\mathbf{C}}_{r}$ is computed from $\hat{\mathbf{A}}_{\rm RR.LS}=\hat{\mathbf{A}}_r$ and $\hat{\mathbf{B}}_{\rm RR.LS}=\hat{\mathbf{B}}_r$ in step 21 of Algorithm \ref{alg:rrminar1}. There are inconsistent subscripts among these iterative estimates. By introducing an error of order $o_\mathrm{p}(T^{-1/2})$, $\hat{\mathbf{B}}_{r-1}$ and $\hat{\mathbf{C}}_{r-1}$ on the right-hand side of step 9 in Algorithm \ref{alg:rrminar1} at the $r$th iteration can be replaced with $\hat{\mathbf{B}}_r$ and $\hat{\mathbf{C}}_r$, which does not affect the subsequent proof of asymptotic normality. Therefore, in the remaining proof of Theorem \ref{t:4.2}, for simplicity, we will remove the RR.LS notation.

By \eqref{new_equ:A.3} and \eqref{new_equ:A.4}, we have
\begin{align*}
&\quad \hat{\mathbf{M}}_1 - \mathbf{M}_1
= (\hat{\mathbf{U}}_1 - \mathbf{U}_1)\hat{\mathbf{D}}_1 \hat{\mathbf{V}}_1^{\top} + \mathbf{U}_1 \hat{\mathbf{D}}_1 \hat{\mathbf{V}}_1^{\top} - \mathbf{U}_1 \mathbf{D}_1 \mathbf{V}_1^{\top} \\
&= (\hat{\mathbf{U}}_1 - \mathbf{U}_1)(\hat{\mathbf{D}}_1 - \mathbf{D}_1)\hat{\mathbf{V}}_1^{\top} + (\hat{\mathbf{U}}_1 - \mathbf{U}_1)\mathbf{D}_1 \hat{\mathbf{V}}_1^{\top} + \mathbf{U}_1\hat{\mathbf{D}}_1 \hat{\mathbf{V}}_1^{\top} - \mathbf{U}_1 \mathbf{D}_1 \mathbf{V}_1^{\top} \\
&= (\hat{\mathbf{U}}_1 - \mathbf{U}_1)(\hat{\mathbf{D}}_1 - \mathbf{D}_1)\hat{\mathbf{V}}_1^{\top}+(\hat{\mathbf{U}}_1 - \mathbf{U}_1)\mathbf{D}_1 (\hat{\mathbf{V}}_1^{\top} - \mathbf{V}_1^{\top})+(\hat{\mathbf{U}}_1 - \mathbf{U}_1)\mathbf{D}_1 \mathbf{V}_1^{\top}+ \mathbf{U}_1 \hat{\mathbf{D}}_1 \hat{\mathbf{V}}_1^{\top} - \mathbf{U}_1 \mathbf{D}_1 \mathbf{V}_1^{\top} \\
&\quad+\mathbf{U}_1 \mathbf{D}_1 \hat{\mathbf{V}}_1^{\top} - \mathbf{U}_1 \mathbf{D}_1 \mathbf{V}_1^{\top} \\
&= (\hat{\mathbf{U}}_1 - \mathbf{U}_1)(\hat{\mathbf{D}}_1 - \mathbf{D}_1)\hat{\mathbf{V}}_1^{\top}+(\hat{\mathbf{U}}_1 - \mathbf{U}_1)\mathbf{D}_1 (\hat{\mathbf{V}}_1^{\top} - \mathbf{V}_1^{\top})+(\hat{\mathbf{U}}_1 - \mathbf{U}_1)\mathbf{D}_1 \mathbf{V}_1^{\top}+ \mathbf{U}_1(\hat{\mathbf{D}}_1 - \mathbf{D}_1)\hat{\mathbf{V}}_1^{\top} \\
&\quad+ \mathbf{U}_1 \mathbf{D}_1(\hat{\mathbf{V}}_1^{\top} - \mathbf{V}_1^{\top}) \\
&=  (\hat{\mathbf{U}}_1 - \mathbf{U}_1)(\hat{\mathbf{D}}_1 - \mathbf{D}_1)\hat{\mathbf{V}}_1^{\top}+(\hat{\mathbf{U}}_1 - \mathbf{U}_1)\mathbf{D}_1 (\hat{\mathbf{V}}_1^{\top} - \mathbf{V}_1^{\top})+(\hat{\mathbf{U}}_1 - \mathbf{U}_1)\mathbf{D}_1 \mathbf{V}_1^{\top}+\mathbf{U}_1(\hat{\mathbf{D}}_1 - \mathbf{D}_1)(\hat{\mathbf{V}}_1^{\top} - \mathbf{V}_1^{\top})  \\
&\quad+\mathbf{U}_1(\hat{\mathbf{D}}_1 - \mathbf{D}_1)\mathbf{V}_1^{\top} + \mathbf{U}_1 \mathbf{D}_1(\hat{\mathbf{V}}_1^{\top} - \mathbf{V}_1^{\top})\\
&=(\hat{\mathbf{U}}_1 - \mathbf{U}_1)\mathbf{D}_1 \mathbf{V}_1^{\top} + \mathbf{U}_1(\hat{\mathbf{D}}_1-\mathbf{D}_1)\mathbf{V}_1^{\top} + \mathbf{U}_1 \mathbf{D}_1(\hat{\mathbf{V}}_1^{\top} - \mathbf{V}_1^{\top}) + o_{\mathrm{p}}\bigg(\frac{1}{\sqrt{T}}\bigg).
\end{align*}
Notice that $\mathbf{D}_1$ is a diagonal matrix, and the last $m - k_{1}$ diagonal entries are zero. Then 
\begin{align*}
(\hat{\mathbf{U}}_1 - \mathbf{U}_1)\mathbf{D}_1 \mathbf{V}_1^{\top}
&=(\hat{\mathbf{U}}_{1:k_{1}} - \mathbf{U}_{1:k_{1}})\mathbf{D}_{1:k_{1}} \mathbf{V}_{1:k_{1}}^{\top},\\
\mathbf{U}_1 \mathbf{D}_1(\hat{\mathbf{V}}_1^{\top} - \mathbf{V}_1^{\top})
&=\mathbf{U}_{1:k_{1}} \mathbf{D}_{1:k_{1}}(\hat{\mathbf{V}}_{1:k_{1}}^{\top} - \mathbf{V}_{1:k_{1}}^{\top}),
\end{align*}
and
\begin{align*}
\mathbf{U}_1(\hat{\mathbf{D}}_1-\mathbf{D}_1)\mathbf{V}_1^{\top}\mathbf{V}_{1:k_{1}}\mathbf{D}_{1:k_{1}}^{-1}
&=\mathbf{U}_1(\hat{\mathbf{D}}_1-\mathbf{D}_1)[\mathbf{I}_{1:k_{1}};\mathbf{0}]\mathbf{D}_{1:k_{1}}^{-1}\\
&=\mathbf{U}_{1:k_{1}}(\hat{\mathbf{D}}_{1:k_{1}}-\mathbf{D}_{1:k_{1}})\mathbf{I}_{1:k_{1}}\mathbf{D}_{1:k_{1}}^{-1}\\
&=\mathbf{U}_{1:k_{1}}(\hat{\mathbf{D}}_{1:k_{1}}-\mathbf{D}_{1:k_{1}})\mathbf{V}_{1:k_{1}}^{\top}\mathbf{V}_{1:k_{1}}\mathbf{D}_{1:k_{1}}^{-1}.
\end{align*}
Since \((\mathbf{I}_{m}-\mathbf{U}_{1:k_{1}}\mathbf{U}_{1:k_{1}}^{\top})\mathbf{U}_{1:k_{1}}=\mathbf{0}\), we have
\begin{align}\label{equ:A.2}
  &(\mathbf{I}_{m}-\mathbf{U}_{1:k_{1}}\mathbf{U}_{1:k_{1}}^{\top})(\mathbf{\hat{U}}_{1:k_{1}}-\mathbf{U}_{1:k_{1}})\mathbf{U}_{1:k_{1}}^{\top}\notag\\
  &=(\mathbf{I}_{m}-\mathbf{U}_{1:k_{1}}\mathbf{U}_{1:k_{1}}^{\top})(\hat{\mathbf{M}}_1 - \mathbf{M}_1)\mathbf{V}_{1:k_{1}}\mathbf{D}_{1:k_{1}}^{-1}\mathbf{U}_{1:k_{1}}^{\top}+\mathit{o}_\mathrm{p}\left(\frac{1}{\sqrt{T}}\right).
\end{align}
By \eqref{new_equ:A.3} and \eqref{new_equ:A.4}, we have
\begin{align*}
\mathbf{\hat{U}}_{1:k_1}\mathbf{\hat{U}}_{1:k_1}^{\top}
&=\mathbf{U}_{1:k_1}\mathbf{U}_{1:k_1}^{\top}+\mathbf{U}_{1:k_1}(\mathbf{\hat{U}}_{1:k_1}-\mathbf{U}_{1:k_1})^{\top}+(\mathbf{\hat{U}}_{1:k_1}-\mathbf{U}_{1:k_1})\mathbf{U}_{1:k_1}^{\top}+\mathit{o}_\mathrm{p}\left(\frac{1}{\sqrt{T}}\right)\\
&=\mathbf{U}_{1:k_1}\mathbf{U}_{1:k_1}^{\top}+\mathbf{U}_{1:k_1}(\mathbf{\hat{U}}_{1:k_1}-\mathbf{U}_{1:k_1})^{\top}(\mathbf{U}_{1:k_1}\mathbf{U}_{1:k_1}^{\top}+\mathbf{I}_m-\mathbf{U}_{1:k_1}\mathbf{U}_{1:k_1}^{\top})\\
&\quad+(\mathbf{\hat{U}}_{1:k_1}-\mathbf{U}_{1:k_1})\mathbf{U}_{1:k_1}^{\top}+\mathit{o}_\mathrm{p}\left(\frac{1}{\sqrt{T}}\right)\\
&=\mathbf{U}_{1:k_1}\mathbf{U}_{1:k_1}^{\top}+\mathbf{U}_{1:k_1}(\mathbf{\hat{U}}_{1:k_1}-\mathbf{U}_{1:k_1})^{\top}(\mathbf{I}_m-\mathbf{U}_{1:k_1}\mathbf{U}_{1:k_1}^{\top})\\
&\quad+(\mathbf{I}_m-\mathbf{U}_{1:k_1}\mathbf{U}_{1:k_1}^{\top})(\mathbf{\hat{U}}_{1:k_1}-\mathbf{U}_{1:k_1})\mathbf{U}_{1:k_1}^{\top}+\mathit{o}_\mathrm{p}\left(\frac{1}{\sqrt{T}}\right).
\end{align*}
Together with \eqref{equ:A.2}, it holds that 
\begin{equation}\label{equ:A.3}
\begin{aligned}  
\mathbf{\hat{U}}_{1:k_1}\mathbf{\hat{U}}_{1:k_1}^{\top}
&=\mathbf{U}_{1:k_1}\mathbf{U}_{1:k_1}^{\top}
+(\mathbf{I}_{m}-\mathbf{U}_{1:k_{1}}\mathbf{U}_{1:k_{1}}^{\top})
(\hat{\mathbf{M}}_1 - \mathbf{M}_1)\mathbf{V}_{1:k_{1}}\mathbf{D}_{1:k_{1}}^{-1}\mathbf{U}_{1:k_{1}}^{\top}\\
  &+\mathbf{U}_{1:k_{1}}\mathbf{D}_{1:k_{1}}^{-1}\mathbf{V}_{1:k_{1}}^{\top}(\hat{\mathbf{M}}_1 - \mathbf{M}_1)^{\top}(\mathbf{I}_{m}-\mathbf{U}_{1:k_{1}}\mathbf{U}_{1:k_{1}}^{\top})+\mathit{o}_\mathrm{p}\left(\frac{1}{\sqrt{T}}\right).
\end{aligned}
\end{equation}
Notice that \(\mathbf{V}_{1:k_{1}}\mathbf{D}_{1:k_{1}}^{-1}\mathbf{U}_{1:k_{1}}^{\top}=\mathbf{M}_{1}^{+}=(\mathbf{A}\Gamma_2^{1/2})^+\). Therefore, by \eqref{equ:A.1} and \eqref{equ:A.3}, we have
\begin{equation}\label{equ:A.3-1}
\begin{aligned}  
\mathbf{\hat{U}}_{1:k_1}\mathbf{\hat{U}}_{1:k_1}^{\top}
&=\mathbf{U}_{1:k_1}\mathbf{U}_{1:k_1}^{\top}+(\mathbf{I}_{m}-\mathbf{U}_{1:k_{1}}\mathbf{U}_{1:k_{1}}^{\top})\left(\frac{1}{T}\sum_{t = 2}^{T}\mathbf{\Delta}_t\mathbf{B}\mathbf{X}_{t-1}^{\top}\right)\Gamma_2^{-1/2}\left(\mathbf{A}\Gamma_2^{1/2}\right)^+\\
&\quad+\left(\Gamma_2^{1/2}\mathbf{A}^{\top}\right)^+\Gamma_2^{-1/2}\left(\frac{1}{T}\sum_{t = 2}^{T}\mathbf{\Delta}_t\mathbf{B}\mathbf{X}_{t-1}^{\top}\right)^{\top}(\mathbf{I}_{m}-\mathbf{U}_{1:k_{1}}\mathbf{U}_{1:k_{1}}^{\top})
+\mathit{o}_\mathrm{p}\left(\frac{1}{\sqrt{T}}\right).
\end{aligned}
\end{equation}
Recall \(\mathbf{\hat{A}}=\mathbf{\hat{U}}_{1:k_1}\mathbf{\hat{U}}_{1:k_1}^{\top}\mathbf{\hat{S}}_{1yx}\mathbf{\hat{S}}_{1xx}^{-1}\) with $\mathbf{\hat{S}}_{1xx}=\sum_{t = 2}^{T}\mathbf{X}_{t - 1}\mathbf{\hat{B}}^{\top}\mathbf{\hat{B}}\mathbf{X}_{t - 1}^{\top}$ and 
$\mathbf{\hat{S}}_{1yx}=\sum_{t = 2}^{T}(\mathbf{X}_{t}-\mathbf{\hat{C}})\mathbf{\hat{B}}\mathbf{X}_{t - 1}^{\top}$. Let \(\mathbf{\hat P}_1 = \mathbf{\hat{U}}_{1:k_1}\mathbf{\hat{U}}_{1:k_1}^{\top}\) and \(\mathbf{P}_1 = \mathbf{U}_{1:k_1}\mathbf{U}_{1:k_1}^{\top}\). Then
\begin{equation*}
\begin{aligned}
&\mathbf{\hat{A}}\bigg(\sum_{t = 2}^{T}\mathbf{X}_{t-1}\mathbf{\hat{B}}^{{\top}}\mathbf{\hat{B}}\mathbf{X}_{t-1}^{\top}\bigg)
=\mathbf{\hat{U}}_{1:k_1}\mathbf{\hat{U}}_{1:k_1}^{\top}\mathbf{\hat{S}}_{1yx}\mathbf{\hat{S}}_{1xx}^{-1}\mathbf{\hat{S}}_{1xx}
=\mathbf{\hat{U}}_{1:k_1}\mathbf{\hat{U}}_{1:k_1}^{\top}\sum_{t = 2}^{T}(\mathbf{X}_t-\mathbf{\hat{C}})\mathbf{\hat{B}}\mathbf{X}_{t-1}^{\top}\\
&=\mathbf{\hat P}_1\left\{\sum_{t = 2}^{T}\mathbf{A}\mathbf{X}_{t-1}\mathbf{B}^{\top}\mathbf{\hat{B}}\mathbf{X}_{t-1}^{\top}+\sum_{t = 2}^{T}(\mathbf{C}-\mathbf{\hat{C}})\mathbf{\hat{B}}\mathbf{X}_{t-1}^{\top}+\sum_{t = 2}^{T}\mathbf{\Delta}_t\mathbf{\hat{B}}\mathbf{X}_{t-1}^{\top}\right\}.
\end{aligned}
\end{equation*}
Since \(T^{-1}\sum_{t = 2}^{T}\mathbf{\Delta}_t\mathbf{B}\mathbf{X}_{t - 1}^{\top} \stackrel{\mathrm{p}}{\to} \mathbf{0}\) and \(\mathbf{\hat P}_1=\mathbf{P}_1+O_{\rm p}(T^{-1/2})\), we have \((\mathbf{\hat P}_1-\mathbf{P}_1)\sum_{t = 2}^{T}\mathbf{\Delta}_t\mathbf{B}\mathbf{X}_{t - 1}^{\top}=o_{\rm p}(T^{1/2})\). Notice that \(\mathbf{A}^+=\mathbf{A}^{\top}(\mathbf{A}\mathbf{A}^{\top})^+\). Hence, by \eqref{new_equ:A.3}-\eqref{new_equ:A.5}, \eqref{equ:A.3-1} and Lemma \ref{lemmaA.5}, 
\begin{align}\label{equ:new1}
    &\mathbf{\hat{A}}\left(\sum_{t = 2}^{T}\mathbf{X}_{t - 1}\mathbf{\hat{B}}^{\top}\mathbf{\hat{B}}\mathbf{X}_{t - 1}^{\top}\right)-\mathbf{A}\left(\sum_{t = 2}^{T}\mathbf{X}_{t - 1}\mathbf{B}^{\top}\mathbf{\hat{B}}\mathbf{X}_{t - 1}^{\top}\right)+\mathbf{P}_1\sum_{t = 2}^{T}(\mathbf{\hat{C}}-\mathbf{C})\mathbf{B}\mathbf{X}_{t - 1}^{\top}\notag\\
    &=(\mathbf{\hat P}_1-\mathbf{I}_{m})\left(\sum_{t = 2}^{T}\mathbf{A}\mathbf{X}_{t - 1}\mathbf{B}^{\top}\mathbf{\hat{B}}\mathbf{X}_{t - 1}\right)+\mathbf{\hat P}_1\sum_{t = 2}^{T}(\mathbf{C}-\mathbf{\hat{C}}+\mathbf{\Delta}_t)\mathbf{\hat{B}}\mathbf{X}_{t - 1}+\mathbf{P}_1\sum_{t = 2}^{T}(\mathbf{\hat{C}}-\mathbf{C})\mathbf{B}\mathbf{X}_{t - 1}^{\top}\notag\\
    &=(\mathbf{I}_m-\mathbf{P}_1)\left(\sum_{t = 2}^{T}\mathbf{\Delta}_t\mathbf{B}\mathbf{X}_{t - 1}^{\top}\right)\Gamma_2^{-\frac{1}{2}}\left(\mathbf{A}\Gamma^{\frac{1}{2}}\right)^+\left(\frac{1}{T}\sum_{t = 2}^{T}\mathbf{A}\mathbf{X}_{t - 1}\mathbf{B}^{\top}\mathbf{\hat{B}}\mathbf{X}_{t - 1}\right)+\mathbf{P}_1\sum_{t = 2}^{T}\mathbf{\Delta}_t\mathbf{B}\mathbf{X}_{t - 1}^{\top}+\mathit{o}_\mathrm{p}(\sqrt{T})\notag\\
    &=(\mathbf{I}_m-\mathbf{P}_1)\left(\sum_{t = 2}^{T}\mathbf{\Delta}_t\mathbf{B}\mathbf{X}_{t - 1}^{\top}\right)\Gamma_2^{-\frac{1}{2}}\left(\mathbf{A}\Gamma^{\frac{1}{2}}\right)^+\mathbf{A}\Gamma_2+\mathbf{P}_1\sum_{t = 2}^{T}\mathbf{\Delta}_t\mathbf{B}\mathbf{X}_{t - 1}^{\top}+\mathit{o}_\mathrm{p}(\sqrt{T})\notag\\
    &=(\mathbf{I}_m-\mathbf{P}_1)\left(\sum_{t = 2}^{T}\mathbf{\Delta}_t\mathbf{B}\mathbf{X}_{t - 1}^{\top}\right)\mathbf{A}^{\top}\left(\mathbf{A}\Gamma_2\mathbf{A}^{\top}\right)^+\mathbf{A}\Gamma_2+\mathbf{P}_1\sum_{t = 2}^{T}\mathbf{\Delta}_t\mathbf{B}\mathbf{X}_{t - 1}^{\top}+\mathit{o}_\mathrm{p}(\sqrt{T}).
\end{align}
By \eqref{new_equ:A.3}-\eqref{new_equ:A.5} and Lemma \ref{lemmaA.5}, it also holds that
\begin{align*}
&\mathbf{\hat{A}}\left(\sum_{t = 2}^{T}\mathbf{X}_{t - 1}\mathbf{\hat{B}}^{\top}\mathbf{\hat{B}}\mathbf{X}_{t - 1}^{\top}\right)-\mathbf{A}\left(\sum_{t = 2}^{T}\mathbf{X}_{t - 1}\mathbf{B}^{\top}\mathbf{\hat{B}}\mathbf{X}_{t - 1}^{\top}\right)+\mathbf{P}_1\sum_{t = 2}^{T}(\mathbf{\hat{C}}-\mathbf{C})\mathbf{B}\mathbf{X}_{t - 1}^{\top} \\
&= \hat{\mathbf{A}}\left\{ \sum_{t=2}^{T} \mathbf{X}_{t - 1}\mathbf{B}^{\top}\mathbf{B}\mathbf{X}_{t - 1}^{\top} + \sum_{t=2}^{T} \mathbf{X}_{t - 1}(\hat{\mathbf{B}}-\mathbf{B})^{\top}\mathbf{B}\mathbf{X}_{t - 1}^{\top} + \sum_{t=2}^{T}\mathbf{X}_{t - 1}\mathbf{B}^{\top}(\hat{\mathbf{B}}-\mathbf{B})\mathbf{X}_{t - 1}^{\top} + O_\mathrm{p}(1) \right\} \\
&\quad - \mathbf{A}\left\{ \sum_{t=2}^{T} \mathbf{X}_{t - 1}\mathbf{B}^{\top}\mathbf{B}\mathbf{X}_{t - 1}^{\top} + \sum_{t=2}^{T} \mathbf{X}_{t - 1}\mathbf{B}^{\top}(\hat{\mathbf{B}}-\mathbf{B})\mathbf{X}_{t - 1}^{\top} \right\} + \mathbf{P}_1\sum_{t = 2}^{T}(\mathbf{\hat{C}}-\mathbf{C})\mathbf{B}\mathbf{X}_{t - 1}^{\top} \\
&= (\hat{\mathbf{A}}-\mathbf{A})\sum_{t=2}^{T} \mathbf{X}_{t - 1}\mathbf{B}^{\top}\mathbf{B}\mathbf{X}_{t - 1}^{\top} + \hat{\mathbf{A}}\bigg\{ \sum_{t=2}^{T}\mathbf{X}_{t - 1}(\hat{\mathbf{B}}-\mathbf{B})^{\top}\mathbf{B}\mathbf{X}_{t - 1}^{\top}  
 + O_\mathrm{p}(1) \bigg\}\\
&\quad   +(\hat{\mathbf{A}}-\mathbf{A})\sum_{t=2}^{T}\mathbf{X}_{t - 1}\mathbf{B}^{\top}(\hat{\mathbf{B}}-\mathbf{B})\mathbf{X}_{t - 1}^{\top} + \mathbf{P}_1\sum_{t = 2}^{T}(\mathbf{\hat{C}}-\mathbf{C})\mathbf{B}\mathbf{X}_{t - 1}^{\top} \\
&= (\hat{\mathbf{A}}-\mathbf{A})\sum_{t=2}^{T} \mathbf{X}_{t - 1}\mathbf{B}^{\top}\mathbf{B}\mathbf{X}_{t - 1}^{\top} + \mathbf{A}\sum_{t=2}^{T}\mathbf{X}_{t - 1}(\hat{\mathbf{B}}-\mathbf{B})^{\top}\mathbf{B}\mathbf{X}_{t - 1}^{\top} + \mathbf{P}_1\sum_{t = 2}^{T}(\mathbf{\hat{C}}-\mathbf{C})\mathbf{B}\mathbf{X}_{t - 1}^{\top} + O_\mathrm{p}(1).
\end{align*}
Together with \eqref{equ:new1}, we have
\begin{equation}\label{equ:A.4}
\begin{aligned}
  &(\mathbf{\mathbf{\hat{A}}}-\mathbf{A})\sum_{t = 2}^{T}\mathbf{X}_{t-1}\mathbf{B}^{\top}\mathbf{B}\mathbf{X}_{t-1}^{\top}+\mathbf{A}\sum_{t = 2}^{T}\mathbf{X}_{t-1}(\mathbf{\hat{B}}-\mathbf{B})^{\top}\mathbf{B}\mathbf{X}_{t-1}^{\top}+\mathbf{P}_1\sum_{t = 2}^{T}(\mathbf{\hat{C}}-\mathbf{C})\mathbf{B}\mathbf{X}_{t-1}^{\top}\\
  &=(\mathbf{I}_m-\mathbf{P}_1)\sum_{t = 2}^{T}(\mathbf{\Delta}_t\mathbf{B}\mathbf{X}_{t-1}^{\top})\mathbf{A}^{\top}(\mathbf{A}\Gamma_2\mathbf{A}^{\top})^+\mathbf{A}\Gamma_2+\mathbf{P}_1\sum_{t = 2}^{T}(\mathbf{\Delta}_t\mathbf{B}\mathbf{X}_{t-1}^{\top})+\mathit{o}_\mathrm{p}(\sqrt{T}).
  \end{aligned}
\end{equation}
Recall \(\Gamma_1 = \mathbb{E}(\mathbf{X}_t^{\top}\mathbf{A}^{\top}\mathbf{A}\mathbf{X}_{t})\). Parnell to \eqref{equ:A.4}, we have
\begin{align}\label{equ:A.5}
&\left(\sum_{t = 2}^{T}\mathbf{X}_{t - 1}^{\top}\mathbf{A}^{\top}\mathbf{A}\mathbf{X}_{t - 1}\right)(\hat{\mathbf{B}}-\mathbf{B})^{\top}+\sum_{t = 2}^{T}\mathbf{X}_{t - 1}^{\top}\mathbf{A}^{\top}(\mathbf{\hat{A}}-\mathbf{A})\mathbf{X}_{t - 1}\mathbf{B}^{\top}+\sum_{t = 2}^{T}\mathbf{X}_{t - 1}^{\top}\mathbf{A}^{\top}(\mathbf{\hat{C}}-\mathbf{C})\mathbf{P}_2 \notag\\
&= \mathbf{A}^{\top}(\mathbf{A}\Gamma_1\mathbf{A}^{\top})^+\mathbf{A}\Gamma_1\sum_{t = 2}^{T}(\mathbf{X}_{t - 1}^{\top}\mathbf{A}^{\top}\mathbf{\Delta}_t)(\mathbf{I}_m-\mathbf{P}_2)+\sum_{t = 2}^{T}(\mathbf{X}_{t - 1}^{\top}\mathbf{A}^{\top}\mathbf{\Delta}_t)\mathbf{P}_2 + \mathit{o}_\mathrm{p}(\sqrt{T}).
\end{align}
Using \(\mathbf{A}\mathbf{X}_{t - 1}\mathbf{B}^{\top}+\mathbf{C} + \mathbf{\Delta}_t\) to replace \(\mathbf{X}_t\) in \eqref{equ:3.3}, we obtain 
\begin{equation}\label{equ:A.6}
\begin{aligned}
  \sum_{t = 2}^{T}(\mathbf{\mathbf{\hat{A}}}-\mathbf{A})\mathbf{X}_{t-1}\mathbf{B}^{\top}+\sum_{t = 2}^{T}\mathbf{A}\mathbf{X}_{t-1}(\mathbf{\hat{B}}-\mathbf{B})^{\top}+\sum_{t = 2}^{T}(\mathbf{\hat{C}}-\mathbf{C}) =\sum_{t = 2}^{T}\mathbf{\Delta}_t+\mathit{o}_\mathrm{p}(\sqrt{T}).
  \end{aligned}
\end{equation}
Taking $\mathrm{vec}(\cdot)$ of both sides of \eqref{equ:A.4}, \eqref{equ:A.5},\eqref{equ:A.6} gives 
\begin{align*}
&\sum_{t=2}^{T}\begin{pmatrix}
(\mathbf{X}_{t-1}\mathbf{B}^{\top}\mathbf{B}\mathbf{X}_{t-1}^{\top})\otimes \mathbf{I}_m & \begin{pmatrix}\mathbf{X}_{t-1}\mathbf{B}^{\top}\end{pmatrix}\otimes(\mathbf{A}\mathbf{X}_{t-1}) & (\mathbf{X}_{t-1}\mathbf{B}^{\top})\otimes \mathbf{P}_1 \\
(\mathbf{B}\mathbf{X}_{t-1})\otimes\begin{pmatrix}\mathbf{X}_{t-1}\mathbf{A}^{\top}\end{pmatrix} & \mathbf{I}_n\otimes (\mathbf{X}_{t-1}\mathbf{A}^{\top}\mathbf{A}\mathbf{X}_{t-1}) & \mathbf{P}_2\otimes (\mathbf{X}_{t-1}\mathbf{A}^{\top}) \\
(\mathbf{B}\mathbf{X}_{t-1}^{\top})\otimes \mathbf{I}_m & \mathbf{I}_n\otimes (\mathbf{A}\mathbf{X}_{t-1}) & \mathbf{I}_n\otimes \mathbf{I}_m
\end{pmatrix}
\begin{pmatrix}
\operatorname{vec}(\mathbf{\hat{A}} - \mathbf{A}) \\
\operatorname{vec}(\mathbf{\hat{B}}^{\top} - \mathbf{B}^{\top} ) \\
\operatorname{vec}(\mathbf{\hat{C}} - \mathbf{C})
\end{pmatrix} \\
&=\sum_{t=2}^{T}\begin{pmatrix}
(\mathbf{X}_{t-1}\mathbf{B}^{\top})\otimes \mathbf{P}_1 + \left[\Gamma_2\mathbf{A}^{\top}\left(\mathbf{A}\Gamma_2\mathbf{A}^{\top}\right)^+\mathbf{A}\mathbf{X}_{t-1}\mathbf{B}^{\top}\right]\otimes(\mathbf{I} - \mathbf{P}_1) \\
\mathbf{P}_2\otimes (\mathbf{X}_{t-1}^{\top}\mathbf{A}^{\top}) + (\mathbf{I} - \mathbf{P}_2)\otimes\left[\Gamma_1\mathbf{B}^{\top}\left(\mathbf{B}\Gamma_1\mathbf{B}^{\top}\right)^+\mathbf{B}\mathbf{X}_{t-1}^{\top}\mathbf{A}^{\top}\right] \\
\mathbf{I}_n\otimes \mathbf{I}_m
\end{pmatrix}
\mathrm{vec}(\mathbf{\Delta}_t) + \mathit{o}_\mathrm{p}(\sqrt{T}),
\end{align*}
which implies
\begin{align}\label{equ:A.7}
\sum_{t=2}^{T}\mathbf{W}_{t-1}
\begin{pmatrix}
\operatorname{vec}(\mathbf{\hat{A}} - \mathbf{A}) \\
\operatorname{vec}(\mathbf{\hat{B}}^{\top} - \mathbf{B}^{\top} ) \\
\operatorname{vec}(\mathbf{\hat{C}} - \mathbf{C})
\end{pmatrix} 
=\sum_{t=2}^{T}\mathbf{Q}_{t-1}\mathrm{vec}(\mathbf{\Delta}_t) + \mathit{o}_\mathrm{p}(\sqrt{T}).
\end{align}
When the error terms are i.i.d., \(\mathbf{X}_t\) is a strictly stationary process (see \cite{Latour1997}). Therefore, by the ergodic theorem, we have 
\begin{equation*}
\frac{1}{T}\sum_{t = 2}^{T}\mathbf{W}_{t-1}\to \mathbb{E}(\mathbf{W}_t),\quad a.s.
\end{equation*}
Recall \(\mathbf{H}=\mathbb{E}(\mathbf{W}_t)+\gamma_1\gamma_1^{\top}\) with \(\gamma_1\coloneqq[\mathrm{vec}(\mathbf{A})^{\top},\mathbf{0}^{\top}]^{\top}\in \mathbb{R}^{m^2+n^2}\). Since \(\|\mathbf{A}\|_F = 1\) and \(\|\mathbf{\mathbf{\hat{A}}}\|_F = 1\), we have \(\text{vec}(\mathbf{A})^{\top}\text{vec}(\mathbf{\mathbf{\hat{A}}}-\mathbf{A})=\mathit{O}_\mathrm{p}(T^{-1})\). By Theorem \ref{t:Iterative convergence}, it holds that
\begin{align*}
\sum_{t=2}^{T}\mathbf{W}_{t-1}
\begin{pmatrix}
\operatorname{vec}(\mathbf{\hat{A}} - \mathbf{A}) \\
\operatorname{vec}(\mathbf{\hat{B}}^{\top} - \mathbf{B}^{\top} ) \\
\operatorname{vec}(\mathbf{\hat{C}} - \mathbf{C})
\end{pmatrix} 
=T\mathbf{H}
\begin{pmatrix}
\operatorname{vec}(\mathbf{\hat{A}} - \mathbf{A}) \\
\operatorname{vec}(\mathbf{\hat{B}}^{\top} - \mathbf{B}^{\top} ) \\
\operatorname{vec}(\mathbf{\hat{C}} - \mathbf{C})
\end{pmatrix}  + \mathit{o}_\mathrm{p}(\sqrt{T}),
\end{align*}
which implies
\begin{equation*}
\mathbf{H}\begin{pmatrix}\mathrm{vec}(\mathbf{\mathbf{\hat{A}}}-\mathbf{A})\\\mathrm{vec}(\mathbf{\hat{B}^{\top}}-\mathbf{B}^{\top})\\\mathrm{vec}(\mathbf{\hat{C}}-\mathbf{C})\end{pmatrix}=\frac{1}{T}\sum_{t = 2}^{T}\mathbf{Q}_{t-1} \rm{vec}(\mathbf{\Delta}_t)+\mathit{o}_\mathrm{p}\left(\frac{1}{\sqrt{T}}\right).
\end{equation*}
Similarly to the proof of Theorem 3.5 in \cite{Kirchner2017}, it holds that \(\{\mathbf{Q}_{t-1}\mathrm{vec}(\mathbf{\Delta}_t)\}_{t=2}^T\) a $(m^2+n^2+mn)$-dimensional vector martingale difference sequence. By the ergodic theorem (see Proposition 7.9 in \cite{Hamilton1994}), we have
\begin{equation*}
\frac{1}{\sqrt{T}} \sum_{t = 2}^{T}\mathbf{Q}_{t-1}\mathrm{vec}(\mathbf{\Delta}_t)\Rightarrow N(0,\mathbb{E}(\mathbf{Q}_t\Sigma_{\mathbf{\Delta}}\mathbf{Q}_t^{\top})),
\end{equation*}
which implies
\begin{equation*}
\sqrt{T}\begin{pmatrix}\mathrm{vec}(\mathbf{\hat{A}}-\mathbf{A})\\ \mathrm{vec}(\mathbf{\hat{B}}^{\top}-\mathbf{B}^{\top})\\ \mathrm{vec}(\mathbf{\hat{C}}-\mathbf{C})\end{pmatrix}\Rightarrow N(0,\mathbf{\Xi}_2),
\end{equation*}
where \(\mathbf{\Xi}_2:=\mathbf{H}^{-1}\mathbb{E}(\mathbf{Q}_{t-1}\Sigma_{\mathbf{\Delta}}\mathbf{Q}_{t-1})\mathbf{H}^{-1}\).
This completes the proof of Theorem \ref{t:4.2}.
\end{proof}

\subsection{Proof of Lemma \ref{t:4.1}}\label{proof:t:4.1}
\begin{proof}
(i) To prove Lemma \ref{t:4.1} (i),  we need Lemma \ref{pro:4.3}, whose proof is given in \cite{Tao2012}.
\begin{lemma}\label{pro:4.3} 
Letting $\mathbf{M}=(\alpha_{i,j})_{n\times n}\in \mathbb{R}_{\geq0}^{n\times n}$ and $\mathbf{N}=(\beta_{i,j})_{n\times n}\in \mathbb{R}_{\geq0}^{n\times n}$ be two Hermitian matrices, then
\[\left\|\left(\lambda_i(\mathbf{M}+\mathbf{N})-\lambda_i(\mathbf{M})\right)_{i=1}^n\right\|_{\ell_n^p} \leq\|\mathbf{N}\|_{S^p}\]
for any $1\leq p\leq \infty $, where
\(\|\left(a_i\right)_{i=1}^n\|_{\ell_n^p}=(\sum_{i=1}^n|a_i|^p)^{1 / p}\)
is the usual $\ell^p$ norm and
\(\|\mathbf{N}\|_{S^p}=\|(\lambda_i(\mathbf{N}))_{i=1}^n\|_{\ell_n^p}\) is the $p$-Schatten norm of $\mathbf{N}$. When $p=2$, the inequality can be written as
\[\sum_{i=1}^n\left|\lambda_i(\mathbf{M+N})-\lambda_i(\mathbf{M})\right|^2 \leq\|\mathbf{N}\|_F^2.\]
\end{lemma}

Let $\lambda_{i}(\hat{\mathbf{\Phi}}_T^{\top} \hat{\mathbf{\Phi}}_T)$ be the $i$th eigenvalue of $\hat{\mathbf{\Phi}}_T^{\top} \hat{\mathbf{\Phi}}_T$. Thus, $\lambda_i(\hat{\mathbf{\Phi}}_T^{\top} \hat{\mathbf{\Phi}}_T) =\hat{\Lambda}_T^2(i, i)$. Similarly, $\lambda_i({\mathbf{\Phi}}^{\top} {\mathbf{\Phi}})={\Lambda}^2(i, i)$. Notice that $\hat{\mathbf{\Phi}}_T^{\top} \hat{\mathbf{\Phi}}_T={\mathbf{\Phi}}^{\top}{\mathbf{\Phi}}+\mathit{o}_\mathrm{p}(T^{-a})$. By Lemma \ref{pro:4.3}, we have 
\[\sum_{i=1}^n\left|\lambda_i\left(\hat{\mathbf{\Phi}}_T^{\top} \hat{\mathbf{\Phi}}_T\right)-\lambda_i\left(\mathbf{\Phi}^{\top}\mathbf{\Phi}\right)\right|^2 \leqslant\left\|\hat{\mathbf{\Phi}}_T^{\top} \hat{\mathbf{\Phi}}_T-\mathbf{\Phi}^{\top} \mathbf{\Phi}\right\|_F^2,\]
where \(n\) is the number of eigenvalues of ${\mathbf{\Phi}}^{\top} \mathbf{\Phi}$ and is bounded, which implies 
\[\left|\lambda_i\left(\hat{\mathbf{\Phi}}_T^{\top} \hat{\mathbf{\Phi}}_T\right)-\lambda_i\left(\mathbf{\Phi}^{\top} \mathbf{\Phi}\right)\right|^2=\mathit{o}_\mathrm{p}\left(T^{-2 a}\right).\]
Therefore,
\[\lambda_i\left(\hat{\mathbf{\Phi}}_T^{\top} \hat{\mathbf{\Phi}}_T\right)-\lambda_i\left(\mathbf{\Phi}^{\top} \mathbf{\Phi}\right)=\mathit{o}_\mathrm{p}\left(T^{-a}\right)\]
for $ i=1,2, \ldots, n$, which implies \[\hat{\mathbf{\Lambda}}_T^2(i, i)=\mathbf{\Lambda}^2(i, i)+\mathit{o}_\mathrm{p}(T^\mathit{-a})\]
for $ i=1,2, \ldots, n$.

(ii) Let $\widehat{{\rm{\mathbf M}}}_j$ and ${\rm{\mathbf M}}_j$ be the $j$th cofactor of $\widehat{\mathbf{\Phi}}_T^{\top}\widehat{\mathbf{\Phi}}_T$ and $\mathbf{\Phi}^{\top}\mathbf{\Phi}$, respectively. By the eigenvector-eigenvalue identity (\cite{Denton2022}), we have
\[\begin{aligned}
|\widehat{\mathbf{\rm{U}}}_T\mathit{(i,j)}|^{2}&=\frac{\prod_{k=1}^{n-1}\left\{\widehat{\mathbf{\Lambda}}^2_T(i,i)-\lambda_{k}(\widehat{\rm{\mathbf M}}_{j})\right\}}{\prod_{k=1, k \neq i}^{n-1}\left\{\widehat{\mathbf{\Lambda}}^2_T(i,i)-\widehat{\mathbf{\Lambda}}^2_T(k,k)\right\}}\\
&=\frac{\prod_{k=1}^{n-1}\left\{\mathbf{\Lambda}^2(i,i)-\lambda_{k}(\rm{\mathbf M}_{j})+\mathit{o}_\mathrm{p}(\mathit{T}^\mathit{-a})\right\}}{\prod_{k=1, k \neq i}^{n-1}\left\{\mathbf{\Lambda}^2(i,i)-\mathbf{\Lambda}^2(k,k)+\mathit{o}_\mathrm{p}(\mathit{T}^\mathit{-a})\right\}}\\
&=\frac{\prod_{k=1}^{n-1}\left\{\mathbf{\Lambda}^2(i,i)-\lambda_{k}(\rm{\mathbf M}_{j})\right\}}{\prod_{k=1, k \neq i}^{n-1}\left\{\mathbf{\Lambda}^2(i,i)-\mathbf{\Lambda}^2(k,k)\right\}+\mathit{o}_\mathrm{p}(\mathit{T}^\mathit{-a})}+\mathit{o}_\mathrm{p}(\mathit{T}^\mathit{-a})\\
&=\frac{|\rm{\mathbf{U}}\mathit{(i,j)}|^{2}}{1+\mathit{o}_\mathrm{p}(\mathit{T}^\mathit{-a})}+\mathit{o}_\mathrm{p}(\mathit{T}^\mathit{-a})\\
&=|{\rm{\mathbf{U}}}\mathit{(i,j)}|^{2}+\mathit{o}_\mathrm{p}(\mathit{T}^\mathit{-a})
\end{aligned}\]
for $1\leq i,j \leq n$.
Therefore, for $1\leq i,j \leq n$, it holds that
\[
\begin{aligned}
|\widehat{\rm{\mathbf{U}}}_T\mathit{(i,j)}|-|\rm{\mathbf{U}}\mathit{(i,j)}|=&\sqrt{|\rm{\mathbf U}\mathit{(i,j)}|^{2}+\mathit{o}_\mathrm{p}(\mathit{T}^\mathit{-a})}-|\rm{\mathbf U}\mathit{(i,j)}|\\
=&\frac{\mathit{o}_\mathrm{p}(\mathit{T}^\mathit{-a})}{\sqrt{|\rm{\mathbf U}\mathit{(i,j)}|^{2}+\mathit{o}_\mathrm{p}(\mathit{T}^\mathit{-a})}+|\rm{\mathbf U}\mathit{(i,j)}|}\\
=&\mathit{o}_\mathrm{p}(\mathit{T}^\mathit{-a}).
\end{aligned}\]

(iii) Similar to the proof of (ii), we have $|\widehat{\mathrm{\mathbf V}}_T\mathit{(i,j)}| = |\mathrm{\mathbf V}\mathit{(i,j)}| + \mathit{o}_\mathrm{p}(\mathit{T}^\mathit{-a})$ for $1\leq i,j\leq n.$

We complete the proof of Lemma \ref{t:4.1}.
\end{proof}

\subsection{Proof of Lemma \ref{lemmaA.4}}\label {proof of lemma A.4}

\begin{proof}
According to the Taylor formula, if $||\mathbf{B}||_2<1$, then
\[
(\mathbf{I}+\mathbf{B})^{-1}=\mathbf{I}+\sum_{k=1}^{\infty}(-1)^{k}\mathbf{B}^{k}.\]
For any $\|\mathbf{B}\|_F=O(a_n)$ with $a_n=o(1)$, if $n$ is sufficiently large, we have
\[||\mathbf{B}||_2\le\|\mathbf{B}\|_F\leq Ca_n<1,\]
where $C\ge0$ is some constant. Notice that $\|\mathbf{B}^{k}\|_F\le\|\mathbf{B}\|_F^{k}\le (Ca_n)^{k}$ for $k=1,2,3,\ldots$. Therefore, for a sufficiently large $n$, it holds that
\begin{align*}
\left\|\sum_{k=1}^{\infty}(-1)^{k}\mathbf{B}^{k}\right\|_F 
    \leq \sum_{k=1}^{\infty}\|\mathbf{B}^{k}\|_F 
    \leq \sum_{k=1}^{\infty}(Ca_n)^{k} 
    = \frac{Ca_n}{1-Ca_n} 
    = O(a_n),
\end{align*}
which implies
\[(\mathbf{I}+\mathbf{B})^{-1}=\mathbf{I}+O(a_n)\]
for any $\|\mathbf{B}\|_F=O(a_n)$ with $a_n=o(1)$. Hence, for any invertible matrix 
$A$, we have
\begin{align*}
\{\mathbf{A} + O(a_n)\}^{-1} 
= \mathbf{A}^{-1}\{\mathbf{I} + \mathbf{A}^{-1} O(a_n)\}^{-1} 
= \mathbf{A}^{-1}\{\mathbf{I} + O(a_n)\}^{-1} 
= \mathbf{A}^{-1} + O(a_n),
\end{align*}
which implies $\hat{\mathbf{A}}^{-1}=\mathbf{A}^{-1}+O_{\mathrm{p}}(a_{n})$ when $\hat{\mathbf{A}} =  \mathbf{A} + O_{\mathrm{p}}(a_{n})$ with invertible matrix $\mathbf{A}\in \mathbb{R}^{m\times m}$ and $a_{n}=o(1)$. This completes the proof of Lemma \ref{lemmaA.4}.
\end{proof}

\subsection{Proof of Lemma \ref{lemmaA.5}}\label{proof of lemmaA.5}
\begin{proof}
(i) Since \(\hat{\mathbf{A}} - \mathbf{A} = O_{\mathrm{p}}\left(T^{-1/2}\right)\), we have, $\forall \epsilon>0$, $\exists {M}_{\epsilon}>0$ such that \(\mathbb{P}(\|\hat{\mathbf{A}} - \mathbf{A}\|_F \leq {M}_{\epsilon}T^{-1/2}) > 1-\epsilon\) for sufficiently large $T>0$. Write \(\boldsymbol{\mu} = \mathbb{E}(\mathbf{X}_t)\). Notice that \(T^{-1}\sum_{t = 1}^{T}\mathbf{X}_t \stackrel{\mathrm{p}}{\to} \boldsymbol{\mu}\). Therefore, for any $\epsilon>0$ and $C>M_{\epsilon}\|\boldsymbol{\mu}\|_{F}$, we have
\begin{align*}
\mathbb{P}\left\{\left\|\sum_{t = 1}^{T}(\hat{\mathbf{A}} - \mathbf{A})\mathbf{X}_t\right\|_F \leq C\sqrt{T}\right\} 
& \geq \mathbb{P}\left\{\left\|\sum_{t = 1}^{T}(\hat{\mathbf{A}} - \mathbf{A})\mathbf{X}_t\right\|_F \leq C\sqrt{T}, \|\hat{\mathbf{A}} - \mathbf{A}\|_F \leq \frac{M_{\epsilon}}{\sqrt{T}}\right\} \\
& \geq \mathbb{P}\left(\frac{M_{\epsilon}}{\sqrt{T}}\left\|\sum_{t = 1}^{T}\mathbf{X}_t\right\|_F \leq C\sqrt{T}, \|\hat{\mathbf{A}} - \mathbf{A}\|_F \leq \frac{M_{\epsilon}}{\sqrt{T}}\right) \\
& \geq \mathbb{P}\left(\left\|\frac{1}{T}\sum_{t = 1}^{T}\mathbf{X}_t\right\|_F \leq \frac{C}{M_{\epsilon}}\right) - \mathbb{P}\left(\|\hat{\mathbf{A}} - \mathbf{A}\|_F > \frac{M_{\epsilon}}{\sqrt{T}}\right) \\
& \ge \mathbb{P}\left(\left\|\frac{1}{T}\sum_{t = 1}^{T}\mathbf{X}_t- \boldsymbol{\mu}\right\|_F \leq \frac{C}{M_{\epsilon}}- \|\boldsymbol{\mu}\|_{F}\right)-\epsilon\\
&= 1-o(1)-\epsilon,
\end{align*}
which implies \(\sum_{t = 1}^{T}(\hat{\mathbf{A}} - \mathbf{A})\mathbf{X}_t = O_{\mathrm{p}}(\sqrt{T})\).

\noindent
(ii) Parallel to the proof of (i), we have \(\sum_{t = 1}^{T}(\hat{\mathbf{A}} - \mathbf{A})^\top(\hat{\mathbf{A}} - \mathbf{A})\mathbf{X}_t = O_{\mathrm{p}}(T^{-1/2})\).

\noindent
(iii) Rewrite the expression in the form of element-by-element product:
\begin{align*}
\sum_{t=2}^T \Delta_t^{\top} (\mathbf{\hat{A}} - \mathbf{A})^{\top} \mathbf{X}_{t-1} = \Bigg(\sum_{t=2}^T \sum_{j=1}^m \sum_{k=1}^m \{\Delta_{j,i,t} (\mathbf{\hat{A}}_{j,k} - \mathbf{A}_{j,k}) \mathbf{X}_{k,l,t-1}\}\Bigg)_{ 1\leq i,l\leq n}. 
\end{align*}
By Theorem \ref{pro:4.1}, $\Delta_t$ and $ \mathbf{X}_{t_1}$ are uncorrelated for any $t_1 < t$, and $\mathbb{E}(\Delta_t) = \mathbf0$. Thus, we have
\[\frac{1}{T} \sum_{t=2}^T \sum_{j=1}^m \sum_{k=1}^m (\Delta_{j,i,t} \mathbf{X}_{k,l,t-1}) \stackrel{\rm p}{\to} 0.\]
Since $\mathbf{\hat{A}} - \mathbf{A} = O_\mathrm{p}(T^{-1/2})$, for any $\epsilon>0$ and $C>0$, it holds that
\begin{align*}
&\mathbb{P}\Bigg[\sum_{t=2}^T \sum_{j=1}^m \sum_{k=1}^m \{\Delta_{j,i,t} (\mathbf{\hat{A}}_{j,k} - \mathbf{A}_{j,k}) \mathbf{X}_{k,l,t-1}\} \leq C\sqrt{T}\Bigg] \\
&\geq \mathbb{P}\Bigg\{\frac{M_{\epsilon}}{\sqrt{T}} \sum_{t=2}^T \sum_{j=1}^m \sum_{k=1}^m (\Delta_{j,i,t} \mathbf{X}_{k,l,t-1}) \leq C\sqrt{T},\ \|\mathbf{\hat{A}} - \mathbf{A}\|_F \leq \frac{M_{\epsilon}}{\sqrt{T}}\Bigg\} \\
&\geq \mathbb{P}\Bigg\{\frac{1}{T} \sum_{t=2}^T \sum_{j=1}^m \sum_{k=1}^m (\Delta_{j,i,t} \mathbf{X}_{k,l,t-1}) \leq \frac{C}{M_{\epsilon}}\Bigg\} - \mathbb{P}\Bigg(\|\mathbf{\hat{A}} - \mathbf{A}\|_F \leq \frac{M_{\epsilon}}{\sqrt{T}}\Bigg) \\
&\ge 1-o(1)-\epsilon,
\end{align*}
where $M_{\epsilon}>0$, $1\leq i,l\leq n$, which implies
\[
\sum_{t=2}^T \Delta_t^{\top} (\mathbf{\hat{A}} - \mathbf{A})^{\top} \mathbf{X}_{t-1} = o_\mathrm{p}(\sqrt{T}).
\]

Similar results also hold for $\hat{\mathbf{B}} - \mathbf{B}$. We complete the proof of Lemma \ref{lemmaA.5}.
\end{proof}

\end{document}